\documentclass[11pt,a4paper]{amsart}

\usepackage{verbatim}
\usepackage[dvips]{color}
\usepackage{amssymb}
\usepackage{mdframed}
\usepackage[active]{srcltx}
\usepackage[breakable]{tcolorbox}
\tcbuselibrary{many}
\tcbset{
    enhanced,
    breakable,
    attach boxed title to top left={
        xshift=0.4cm,
        yshift= -3.5mm, 
    },
    top=4mm,
    coltitle=black,
    beforeafter skip=\baselineskip,
}
\usepackage[colorlinks,pdfpagelabels,pdfstartview = FitH,bookmarksopen = true,bookmarksnumbered = true,linkcolor = blue,plainpages = false,hypertexnames = false,citecolor = red] {hyperref}

\def\eqn#1$$#2$${\begin{equation}\label#1#2\end{equation}}
\def\charfn_#1{{\raise1.2pt\hbox{$\chi_{\kern-1pt\lower3pt\hbox{{$\scriptstyle#1$}}}$}}}
\newcommand{\rif}[1]{(\ref{#1})}
\newcommand{\trif}[1] {\textnormal{\rif{#1}}}

\numberwithin{equation}{section}
\newenvironment{boxy}[1]{%
    \begin{tcolorbox}[title={#1}]%
    }{
    \end{tcolorbox}
}

\newcommand{\eps}{\varepsilon}
\newcommand{\hhh}{\textnormal{\texttt{h}}}

\newcommand{\mf}[1]{\mathfrak{#1}}

\def\dx{\,{\rm d}x}
\def\dy{\,{\rm d}y}
\def\d{{\rm d}}
\newcommand{\kk}{\kappa}
\newcommand{\rr}{\varrho}
\newcommand{\sss}{\mathfrak{s}}

\newcommand{\data}{\textnormal{\texttt{data}}}
\newcommand\ia{i_{\textnormal{a}}}
\newcommand\sa{s_{\textnormal{a}}}

\def\dist{\operatorname{dist}}

\newcommand{\divo}{\textnormal{div}}

\newcommand{\ssf}{\mathfrak{s}}

\def\en{\mathbb N}
\def\er{\mathbb R}

\newcommand{\ao}{\alpha_0}

\newcommand{\tia}{\tilde \alpha}

\def\loc{\operatorname{loc}}

\newcommand\ccc{\mathfrak{c}}
\newcommand\aaa{\mathfrak{a}}
\def\deb{\rightharpoonup}
\newcommand\ppp{\mathfrak{p}}

\newcommand\ssss{\mathfrak{s}}

\newcommand{\zzero}{0_{\er^n}}

\newcommand{\mint}{\mathop{\int\hskip -1,05em -\, \!\!\!}\nolimits}

\def\vintslides_#1{\mathchoice%
          {\mathop{\kern 0.1em\vrule width 0.5em height 0.697ex depth -0.581ex
                  \kern -0.6em \intop}\nolimits_{\kern -0.4em#1}}%
          {\mathop{\kern 0.1em\vrule width 0.3em height 0.697ex depth -0.604ex
                  \kern -0.4em \intop}\nolimits_{#1}}%
          {\mathop{\kern 0.1em\vrule width 0.3em height 0.697ex depth -0.604ex
                  \kern -0.4em \intop}\nolimits_{#1}}%
          {\mathop{\kern 0.1em\vrule width 0.3em height 0.697ex depth -0.604ex
                  \kern -0.4em \intop}\nolimits_{#1}}}

\newcommand{\aveint}[2]{\mathchoice%
          {\mathop{\kern 0.2em\vrule width 0.6em height 0.69678ex depth -0.58065ex
                  \kern -0.8em \intop}\nolimits_{\kern -0.45em#1}^{#2}}%
          {\mathop{\kern 0.1em\vrule width 0.5em height 0.69678ex depth -0.60387ex
                  \kern -0.6em \intop}\nolimits_{#1}^{#2}}%
          {\mathop{\kern 0.1em\vrule width 0.5em height 0.69678ex depth -0.60387ex
                  \kern -0.6em \intop}\nolimits_{#1}^{#2}}%
          {\mathop{\kern 0.1em\vrule width 0.5em height 0.69678ex depth -0.60387ex
                  \kern -0.6em \intop}\nolimits_{#1}^{#2}}}
                  
\newtoks\by
\newtoks\paper
\newtoks\book
\newtoks\jour
\newtoks\yr
\newtoks\pages
\newtoks\vol
\newtoks\publ

\def\name[#1, #2]{#1 #2}
\def\ota{{\hbox{\bf ???}}}
\def\cLear{\by=\ota\paper=\ota\book=\ota\jour=\ota\yr=\ota
\pages=\ota\vol=\ota\publ=\ota}
\def\enD\ellaper{\the\by, \textit{\the\paper},
{\the\jour} \textbf{\the\vol} (\the\yr), \the\pages.\cLear}
\def\endbook{\the\by, \textit{\the\book},
\the\publ, \the\yr.\cLear}
\def\enD\ellap{\the\by, \textit{\the\paper}, \the\jour.\cLear}
\def\enD\ellroc{\the\by, \textit{\the\paper}, \the\book, \the\publ,
\the\yr, \the\pages.\cLear}

\newtheorem{theorem}{Theorem}[section]

\newtheorem{commentary}[theorem]{Commentary}
\newtheorem{corollary}[theorem]{Corollary}
\newtheorem{lemma}[theorem]{Lemma}

\theoremstyle{definition}
\newtheorem{definition}[theorem]{Definition}
\newtheorem{remark}[theorem]{Remark}

\title[Nonlinear potential theoretic methods]{Nonlinear potential theoretic methods in nonuniformly ellliptic problems}

\author[Mingione]{Giuseppe Mingione}  \address{Giuseppe Mingione\\Dipartimento SMFI, Universit\`a di Parma, Parco Area delle Scienze 53/a, Campus, 43124 Parma, Italy} \email{\url{giuseppe.mingione@unipr.it}}

\begin{document}

\begin{abstract}
Nonuniform ellipticity is a classical topic in the theory of partial differential equations. While several results in regularity theory have been adding up over decades, many basic issues, as for instance the validity of Schauder theory and sharp dependence of regularity upon data, remained opened for a while. In these notes we give an overview of recent results and techniques about the topic, that, via a novel use of nonlinear potential theoretic methods, allow to answer several of the above questions.    \end{abstract}


\keywords{Regularity, Nonlinear Potentials, Uniform Ellipticity, Nonuniform Ellipticity, Nonstandard growth conditions. \vspace{1mm}}

\maketitle

\tableofcontents

\section{About the content}\label{sezione1}
The aim of these notes is twofold, that is
\begin{itemize}
\item To briefly survey a few recent results in the regularity theory of minimizers of variational integrals of the type 
\eqn{genF}
$$
W^{1,1}_{\loc}(\Omega,\er^N) \ni w \mapsto \mathcal F(w,\Omega):=  \int_{\Omega}  \left[ F(x,Dw)-\mu \cdot w \right]\dx
$$ 
and of solutions to elliptic equations of the type
\eqn{risolvi}
$$
-\divo\, a(x,Du) = \mu \qquad \mbox{in}\  \Omega \subset \er^n\,.
$$
In both cases, and for the rest of the paper, $\Omega \subset \er^n$ denotes a bounded open subset and $n\geq 2$ unless otherwise specified.  We shall refer to the scalar case when $N=1$ and to the vectorial one otherwise. In fact, unless otherwise specified, we shall always deal with the scalar case $N=1$. In particular, this will always happen when dealing with \rif{risolvi}, i.e., we shall always deal with equations but in the case of Theorem \ref{sample6}.

As from the title, the emphasis will be put on {\em nonuniformly elliptic problems}. 
\item In connection to the previous point, we shall try to give a streamlined presentation of a few facts and techniques based on the use of Nonlinear Potential Theory in the setting of nonuniformly elliptic problems; see Section \ref{rinormalizza}. In this respect, we shall sacrifice generality to readability. We shall indeed present a series of cases and results treated via a unified approach using a certain family of nonlinear potentials. Several of such cases can be then treated in a more efficient via specific variations of the general method presented, but still based on the use of nonlinear potentials. These variations will be then described in the commentaries of Section \ref{rinormalizza} and references will be given.
\end{itemize}

In \rif{genF} the integrand $F\colon \Omega \times \er^{N\times n} \to \er$ will be assumed to be non-negative and Carath\'eodory regular. Moreover, again when dealing with functionals \rif{genF}, we shall assume (at least) that
\eqn{minimamu}
$$
\mu\in L^n(\Omega, \er^N)
$$
so that $\mu \cdot w$ is locally integrable whenever $w\in W^{1,1}_{\loc}(\Omega, \er^N)$ as an effect of Sobolev embedding theorem. As for \rif{risolvi}, the vector field $a\colon \Omega \times \er^n \to \er^n$ will always be at least Carath\'eodory regular and $\mu$ will in the most general case be a Borel measure. Minimizers of the functional in \rif{genF} and solutions to \rif{risolvi} are defined as follows, respectively:
\begin{definition}\label{defi-min} A map $u \in W^{1,1}_{\loc}(\Omega,\er^N)$ is a (local) minimizer of the functional $\mathcal F$ in \eqref{genF} with $\mu \in L^n_{\loc}(\Omega,\er^N)$ if, for every open subset $\tilde \Omega\Subset \Omega$, we have $F(\cdot, Du)\in L^1(\tilde \Omega)$ and if $\mathcal F(u,\tilde \Omega)\leq \mathcal F(w,\tilde \Omega)$ holds for every competitor $w \in u + W^{1,1}_0(\tilde \Omega, \er^N)$. 
\end{definition}
\begin{definition}\label{defisol0} A function $u \in W^{1,1}_{\loc}(\Omega)$ is a distributional (local) solution to \rif{risolvi} iff $a(\cdot, Du) \in L^1_{\loc}(\Omega, \er^n)$, $\mu$ is a Borel measure with locally finite total mass in $\Omega$ and 
\eqn{distrisol}
$$
\int_{\Omega} a(x,Du)\cdot D\varphi \dx = \int_{\Omega} \varphi\, \d\mu
$$
holds for every $\varphi \in C^{\infty}_{0}(\Omega)$. 
\end{definition}
These are not yet so-called energy solutions, for which we refer to Definition \ref{defisol}, and are actually called very weak solutions - see discussion after Theorem \ref{c1} in Section \ref{piovra1}. In these notes we shall use no more than  integrability properties of $\mu$. As all the results presented here will be interior regularity results, without loss of generality we shall always assume that the integrability properties of $\mu$ in question will be global, as for instance in \rif{minimamu}. Eventually letting $\mu\equiv 0$ outside $\Omega$, we shall assume that $\mu$ is defined, with the same integrability properties, on the whole $\er^n$.  

This paper finds its origins in a series of lectures given at the 2022 CIME school and, as mentioned before, they are mostly concerned with nonuniformly elliptic problems. For nonlinear potential theoretic results and related Nonlinear Calder\'on-Zygmund theory in the uniformly elliptic setting the reader might like to try the notes from the previous 2016 CIME lectures \cite{mincime}.

\subsection{Notation}\label{notazione} 
By $c$ we shall denote generic constants larger than $1$; such constants will change their specific value in different occurrences, although used in the same setting. What will really matter will be their dependence on various relevant parameters. This will be usually specified in parentheses. A similar role will be played by constants denoted by $\vartheta$, $\tilde c$, $\kk$ and the like; this will be clear by the context and their exact values will not be important, but their dependence on parameters will. We shall use the symbol $\lesssim$ to mean an inequality that occurs up to a universal constant, i.e., a constant usually depending on $n$, or on a fixed number of parameters whose precise value is not relevant in the context. For instance, $a\lesssim b$ means that there exists a constant $c\geq 1$ such that $a \leq c b$. In the case $c$ depends itself on parameters, for instance $\kk_1, \kk_2$, and it is useful to specify such a dependence, we shall denote $a\lesssim_{\kk_1, \kk_2} b$. Finally, $a \approx b$ means that both $a\lesssim b$ and $b \lesssim a$ occur. In a similar way, writing $a \equiv_{\kk_1, \kk_2}b$ means that $a =cb$ for a constant $c$ depending on $\kk_1, \kk_2$ and so on. In the following $$ B_{r}(x_0):=\{x \in \er^n \, : \,  | x-x_0|< r\}$$ denotes the open ball with center $x_0$ and radius $r>0$. When not important, or when it will be clear from the context, we shall omit denoting the center as follows: $B_r \equiv B_{r}(x_0)$.
As usual, the Sobolev embedding exponent $2^*$ is defined as 
\eqn{sobby0}
$$
2^*:=\begin{cases}
\ \mbox{$\frac{2n}{n-2}$ if $n>2$}\\
\ \mbox{any number $>2$ if $n=2$}\,.
\end{cases}
$$
With $\mathcal B \subset \er^{n}$ being a measurable subset with positive measure, and with 
$f \colon \mathcal B \to \er^{k}$, $k\geq 1$, being an integrable map, we shall denote by  $$
   (f)_{\mathcal B} \equiv \mint_{\mathcal B}  f \, dx  := \frac{1}{|\mathcal B|}\int_{\mathcal B}  f(x) \, d x
$$
its integral average; here $|\mathcal B|$ denotes the Lebesgue measure of $\mathcal B$. Needless to say, we shall automatically identify $L^1$-functions $\mu$ with Borel measures, thereby denoting 
$$
|\mu|(\mathcal B) = \int_{\mathcal B} |\mu|\, dx \qquad \mbox{for every measurable subset} \ \mathcal B \Subset \Omega\;.
$$
In denoting several function spaces like $L^t(\Omega), W^{1,t}(\Omega)$, we will denote the vector valued version by $L^t(\Omega,\er^N), W^{1,p}(\Omega,\er^N)$ in the case the maps considered take values in $\er^N$, $N\in \en$. When clear from the contest, we will also abbreviate $L^t(\Omega,\er^N) , W^{1,t}(\Omega,\er^N)\equiv 
L^t(\Omega) , W^{1,t}(\Omega)$ and so on. In the rest of the paper $p, q, \nu, L$ and $s$ will denote numbers such that
\eqn{assunti}
$$
1 < p \leq q\,,\quad 0 < \nu \leq L\,, \quad  0 \leq s \leq 1\,.
$$
In relation to \rif{assunti}, we set
$$
\data = (n,p,q,\nu, L)\,,
$$
where $n$ is the ambient dimension of the equations or the integrals we are considering, as indicated after \rif{genF}-\rif{risolvi}. This notation will still be used when $p=q$ and in this case $\data = (n,p,\nu, L)$. Finally, we shall also denote 
\eqn{defiH}
$$
H_{s}(z):=|z|^{2}+s^{2}\,, \quad z \in \er^n\,.
$$


\section{Notions of nonuniform ellipticity}
In this section we discuss a few definitions of uniform and nonuniform ellipticity for functionals and equations of the type in \rif{genF} and \rif{risolvi}, respectively. Major emphasis will be put on nonautonomous problems, that is, problems with coefficients. 

\subsection{Classical (pointwise) definitions}\label{unisec}
Let us start with an elliptic equation in divergence form of the type in \rif{risolvi}, where the vector field $a\colon \Omega \times \er^n \to \er^n$ is initially assumed to be Carath\'eodory regular together with its derivatives $\partial_z a(\cdot)$ being Carath\'eodory regular in $\Omega\times \er^n\setminus\{0_{\er^n}\}$. Moreover, in the following we shall always assume that $z \mapsto a(\cdot, z)\in C^{0}(\er^n)\cap C^{1}(\er^n\setminus\{0_{\er^n}\})$. 
In order to give a quantitative description of its ellipticity properties, we shall consider two non-negative, Carath\'eodory functions
$g_1, g_2\colon \Omega \times  (0, \infty) \to (0,\infty)$
that we are going to use as lower and upper bounds on the eigenvalues of $z \mapsto \partial_za (x,z)$, respectively. Specifically, we shall assume that
\eqn{doppiob}
$$
\begin{cases}
\ g_1(x,|z|) \mathbb{I}_{\rm d} \leq \partial_z a(x,z)\\
\ |\partial_z a(x,z)|\leq g_2(x,|z|)
\end{cases}
$$
holds whenever $z \in \er^n\setminus \{\zzero\}$ and $x\in \Omega$, with $ \mathbb{I}_{\rm d}$ denoting the identity.  Uniform ellipticity of the vector field $a(\cdot)$ (or, equivalently, of the equation \rif{risolvi}) then requires that 
\begin{boxy}{The Ellipticity Ratio}
\eqn{ratio1}
$$
\mathcal R_{a} (x,z):=\frac{g_2(x,|z|)}{g_1(x,|z|)}
$$
\end{boxy}
\noindent remains uniformly bounded, i.e., 
\eqn{limitata}
$$
\sup_{x\in \Omega, |z|\not=0}\mathcal R_{a} (x,z)< \infty\,.
$$
When $\partial_z a(\cdot)$ is symmetric this amounts to require that the ratio
$$
\frac{\mbox{highest eigenvalue of}\ \partial_{z}  a(x,z)}{\mbox{lowest eigenvalue of}\  \partial_{z} a(x,z)}
$$
is still uniformly bounded in the above range of $x$ and $z$; in this case the quantity in the last display provides an alternative, more tailored definition of ellipticity ratio. In the same line of thought, when considering integral functionals of the type in \rif{genF}, we shall say that the integrand $F(\cdot)$ (or, equivalently, the functional $\mathcal F$) is uniformly elliptic if so is $\partial_zF(\cdot)$. That is, if and only if its Euler-Lagrange equation
$$
-\divo\, \partial_zF(x,Du)=\mu 
$$
is uniformly elliptic. 
In this case we shall often adopt the shortened notation $\mathcal R_{F}\equiv \mathcal R_{\partial_{z}F}$. Sometimes, we shall denote $\mathcal R \equiv \mathcal R_{F}$ or $\mathcal R \equiv \mathcal R_{a}$ when no ambiguity shall arise on the identity of $F(\cdot)$ or $a(\cdot)$. 

On several occasions it is also sufficient to adopt a slightly 
different notion of uniform ellipticity, namely, instead of \rif{limitata} one might require that 
\eqn{limitata2}
$$
\sup_{x\in \Omega, |z|\geq 1}\mathcal R_{a} (x,z)< \infty
$$
or even 
\eqn{limitata3}
$$
\limsup_{|z|\to \infty}\, \mathcal R_{a} (x,z)< \infty, \ \  \mbox{uniformly with respect to $x\in \Omega$}\,.
$$
This happens especially when one is interested in proving local Lipschitz regularity of solutions so that only the behaviour of the operator for large values of $|z|$ (that is, of the gradient) matters. In this case \rif{limitata2}-\rif{limitata3} become more suitable definitions. 

On the contrary, we shall say that the equation in \rif{risolvi} is nonuniformly elliptic if, for at least one point $x\in \Omega$, it happens that
\eqn{classe}
$$
\sup_{|z|>0}\mathcal R_{a} (x,z)= \infty\,.
$$
The importance of conditions like \rif{limitata} or \rif{limitata3} stems from the fact that when trying to prove Lipschitz bounds for solutions to equations as \rif{risolvi}, the quantity $\mathcal R_{a} (x,Du(x))$ appears in all the crucial integral estimates. Examples of this will appear in \rif{para202} and \rif{similar}. Assuming \rif{limitata}, and therefore the a priopri boundedness of $\mathcal R_{a} (x,Du(x))$ no matter how large the size of $|Du(x)|$ is, allows for a whole range of estimates and results. These must be then rediscussed in the nonuniformly elliptic case \rif{classe}, up to the stage that a too fast growth of 
\eqn{cresci}
$$z \mapsto \mathcal R_{a} (x,z)$$ 
with respect to $|z|$ implies the existence of non-Lipschitz, and even unbounded solutions \cite{gia,M2}. The whole point in the regularity theory of nonuniformly elliptic equations is then to find suitable growth assumptions on the map in \rif{cresci} allowing to compensate the potential blow-up of the gradient of solutions, and therefore of $\mathcal R_{a} (x,Du(x))$, eventually leading to local Lipschitz regularity of solutions. Once this is known, the growth of \rif{cresci} becomes irrelevant as $Du$ remains bounded and so (typically) $\mathcal R_{a} (x,Du(x))$ does anyway. In this situation higher regularity of solutions can be often obtained adapting the methods developed for the uniformly elliptic case. Imposing a growth condition on \rif{cresci} is a classical and natural approach in the literature \cite{IO, LUcpam, M1, M2, M3, serrin, tru1967} and we shall also follow it. Examples of this are given in Section \ref{gapbounds} below. 
\subsection{Examples of uniform operators}
Beside the Laplacean operator, obviously, and linear variations including coefficients, the next popular uniformly elliptic operator is  
\eqn{modellop}
$$
\begin{cases}
\, -\divo\, (\ccc(x)(|Du|^{2}+s^2)^{(p-2)/2}Du)=0\\
\, 0 < \nu \leq \ccc(x)\leq L, \quad 
s\in [0,1], \quad p>1
\end{cases}
$$
which is the Euler-Lagrange equation of the functional
$$
w \mapsto \int_{\Omega} (|Dw|^{2}+s^2)^{p/2}\dx \,.
$$
In this case it is
\eqn{modellino}
$$
\begin{cases}
\, a(x,z)=\ccc(x)(|z|^{2}+s^2)^{(p-2)/2}z\\
\,   \partial_z a(x,z)\approx_{p} \ccc(x)(|z|^{2}+s^2)^{(p-2)/2} \mathbb{I}_{\rm d}\\
\, g_1(x,|z|)\approx g_2(x,|z|) \approx \ccc(x) (|z|^{2}+s^2)^{(p-2)/2} 
\end{cases}
$$
and 
$$
\mathcal R_{a} (x,z) \lesssim \frac{L}{\nu} \frac{\max\{p-1,1\}}{\min\{p-1,1\}}\,.
$$
The choice $s=0$ and $\ccc(\cdot)\equiv 1$ leads to 

\begin{boxy}{The $p$-Laplacean operator}
\eqn{chiefp}
$$
-\Delta_p u = -\divo\, (|Du|^{p-2}Du)=0\,.
$$
\end{boxy}
\noindent See \cite{Ur, Uh, Dnonana, manth1, manth2} for basic regularity in this last case; conditions on coefficients can be also relaxed, see for instance \cite{balcina} and related references. The one in \rif{modellop} belongs to the family of operators with standard polynomial $p$-growth. These have been treated at length in the literature where general equations of the type \rif{risolvi} have been considered  under the assumptions 
\eqn{ass1}
$$
\left\{
\begin{array}{c}
|a(x,z)|+|\partial_{z} a(x,z)|(|z|^{2}+s^2)^{1/2} \leq L(|z|^{2}+s^2)^{(p-1)/2}\\ [8 pt]
\nu(|z|^{2}+s^2)^{(p-2)/2}|\xi|^2\leq  \partial_z a(x,z)\xi \cdot \xi  
 \end{array}\right.
 $$
for every choice of $x\in \Omega$ and $z, \xi \in \er^n$, $|z|\not=0$; as already specified in \rif{assunti}, it is $0< \nu \leq L$ and $p>1$. These assumptions are classical after the work of Ladyzhenskaya \& Uraltseva \cite{LU} and are modelled on \rif{modellino}. In this case we have 
$$
\mathcal R_{a} (x,z) \lesssim_{p} \frac{L}{\nu}\,.
$$
Sometimes conditions \rif{ass1} are replaced by the weaker
\eqn{assumi}
$$
\nu |z|^p- L\leq a(x,z)\cdot z\,, \quad  |a(x,z)|\leq L |z|^{p-1}+ L\,,
$$
where $a(\cdot)$ is only supposed to be Carath\'eodory regular. The growth conditions \rif{ass1}$_1$ fix the concept of (local) {\em energy solutions} as distributional solutions that belong to the Sobolev  space $W^{1,p}_{\loc}(\Omega)$, that is
\begin{definition}\label{defisol} A function $u$ is an energy solution to \trif{risolvi} under assumptions \trif{assumi} and $\mu \in W^{-1, p'}(\Omega)$ (the dual of $W^{1,p}_0(\Omega)$), iff $u \in W^{1,p}_{\loc}(\Omega)$ and it is a distributional solution in the sense of Definition \ref{defisol0}. In particular \trif{distrisol} 
holds for every $\varphi \in W^{1,p}(\Omega)$ which is compactly supported in $\Omega$. 
\end{definition}
Notice that, as already mentioned, conditions \rif{assumi} are implied by \rif{ass1} (modulo changing the values of $\nu, L$), see also \cite{AKM}. This is the type of solution which is usually considered in this setting. Uniformly elliptic problems are indeed often characterized by the fact of having a natural and robust notion of energy solutions, which is indeed fixed by the growth conditions of the vector field $a(\cdot)$, as in the present setting. 

Uniformly elliptic equations are not necessarily of $p$-polynomial type as in \rif{ass1}, but growth conditions intrinsically determined by certain nonlinear functions can be considered too. We are thinking of operators of the type 
\eqn{uni-ell0}
$$ 
- \divo(\tilde a(|Du|)Du)=0
$$
under assumptions 
\eqn{uni-ell}
$$
\left\{
\begin{array}{c}
\displaystyle -1< \ia \leq \frac{\tilde a'(t)t}{\tilde a(t)}\leq \sa< \infty \quad \mbox{for every $t>0$}\\ [10 pt]
\mbox{$\tilde a\colon (0, \infty) \to [ 0, \infty)$ is of class $C^1_{\loc}(0, \infty)$}\,.
\end{array}
\right.
$$ 
In this case we find 
$$\mathcal R_{a} (z)\lesssim \frac{\max\{\sa+1,1\}}{\min\{\ia+1,1\}}\,.$$
The one in \rif{uni-ell0} is the Euler-Lagrange equation of the functional 
\eqn{defiA}
$$
w \mapsto \int_{\Omega} A(|Dw|)\dx,  \qquad 
A(t) := \int_0^t \tilde a(y)y\, dy\,.
$$
Notice that \rif{modellino} falls in the realm of \rif{uni-ell0} with 
$\tilde a(t)\equiv (t^2+s^2)^{(p-2)/2}$ so that $\ia=\sa=p-2$ and $A(t)\equiv (t^2+s^2)^{p/2}/p$. 
These problems are naturally well-posed in the Orlicz space $W^{1, A(\cdot)}_{\loc}(\Omega)$, defined as the set of $W^{1,1}_{\loc}(\Omega)$-regular functions $w$ such that 
 $A(|Dw|)\in L^1_{\loc}(\Omega)$. In turn, energy solutions to \rif{uni-ell0} are defined as distributional solutions belonging to $W^{1, A(\cdot)}_{\loc}(\Omega)$. An instance of a regularity result valid for energy solutions in this setting is presented in Theorem \ref{sample6}.  A typical example of a functional as in \rif{defiA}  is provided by 
 \eqn{loggy}
 $$
 w \mapsto \int_{\Omega} |Dw|^p\log(1+|Dw|)\dx\,, \qquad p>1\,.
 $$
For such problems, and related regularity theory of solutions, references we recommend are for instance \cite{CMPDE, CMarma, CMjems, CMsurvey, CMjmpa, diening, HO1, HO2, liebe1, stroffa}. Observe that in the limiting case $p=1$ the functional in \rif{loggy} ceases to be uniformly elliptic. See Section \ref{esempi} below.

\subsection{Examples of nonuniform operators}\label{esempi}
Functionals as
$$
w \mapsto \mathcal F_j(w, \Omega) :=\int_{\Omega} F_j(x, Dw)\dx 
$$
for $j \in \{1, \ldots, 4\}$, 
where
$$
\begin{cases}
\ \displaystyle  F_1(x,Dw) := |Dw|^p + \sum_{k=1}^m \aaa_k(x)|D_kw|^{p_i}, \ \  1 < p \leq p_k, \ \ 0 \leq \aaa_k(\cdot) \leq L \\ \\
\ \displaystyle  F_2(x,Dw)  :=    \exp(|Dw|^p),  \quad p\geq 1\\ \\ 
\ \displaystyle  F_3(x,Dw)  :=   \exp(\exp((|Dw|^p))), \quad p\geq 1  \\ \\  
\  \displaystyle 
F_4(x,Dw)  :=  \ccc(x)|Dw|\log(1+|Dw|), \quad 0 < \nu \leq \ccc(x) \leq L\,,
\end{cases}
$$
are all nonuniformly elliptic. In fact, the best upper bounds we can find on the ellipticity ratios, defined in \rif{ratio1}, are 
\eqn{ratios}
$$
\begin{cases}
\ \displaystyle\mathcal R_{F_1}(z)\lesssim_{a_i(\cdot)} |z|^{q-p}+1, \quad q:= \max\{p_i\}\\
\ \displaystyle \mathcal R_{F_2}(z)\lesssim |z|^{p} +1\\
\ \displaystyle \mathcal R_{F_3}(z)\lesssim |z|^{p}\exp(|z|^p) +1
\\
\ \displaystyle \mathcal R_{F_4}(z)\lesssim_{L/\nu} \log(1+|z|)+1\,.
\end{cases}
$$
See also \cite{seminale, ciccio} for the computations relevant to \rif{ratios}. 
$\mathcal F_1$ is a typical example of functional with so-called $(p,q)$-growth conditions. The ellipticity ratio grows as $|Du|^{q-p}$, and therefore polynomially with respect to the gradient of a minimizer $u$. We shall come back on such functionals later on, in Section \ref{gapbounds}. Integrands where each partial derivative is penalized with its own exponent are called anisotropic integrands and are amongst the most intensively studied within the realm of nonuniform operators; see for instance \cite{BB, M1, UU}. Functionals $\mathcal F_2$ and $\mathcal F_3$ are relevant examples of functionals with fast, non-polynomial growth conditions, considered for instance in \cite{seminale, ciccio, dima, DE, evans, liebe2, M4, M5}; see Section \ref{expsec}. In the case of $\mathcal F_3$ the ellipticity ratio still grows polynomially. Finally, $\mathcal F_4$ is an example of a functional with so-called nearly linear growth. See Section \ref{logsec} for more. 

Nonuniform ellipticity is a classical topic in the theory of PDE. In the polynomial growth case authors who worked on the topic are Ladyzhenskaya \& Ural'tseva \cite{LU, LUcpam}, Hartman \& Stampacchia \cite{hast}, Trudinger \cite{tru1967, tru}, Ivočkina \& A.\,P. Oskolkov \cite{IO}, Oskolkov \cite{osk}, Serrin \cite{serrin}, A.\,V. Ivanov \cite{ivanov0, ivanov1, ivanov2},  Leon Simon \cite{simon0, Simongrad}, Ural'tseva \& Urdaletova \cite{UU}, Lieberman \cite{liebe0}, just to mention a few. We shall describe more recent contributions later on, in Section \ref{gapbounds}.

\subsection{Soft nonuniform ellipticity \cite{ciccio}}\label{softy2} Uniform ellipticity is a way to express that the eigenvalues of a given operator, that at this stage depend on the solution itself as described in \rif{doppiob}, self-rebalance. Ultimately, they mimic the situation of the classical Poisson equation. Of course substantiating this last assertion requires the full strength of the entire nonlinear regularity theory, starting from De Giorgi-Nash-Moser. Passing from the linear to the nonlinear case requires building a completely different world of ideas. Nevertheless, the final outcome is that things like 
\eqn{elenco}
$$
\begin{cases}
\,  \mbox{Schauder type estimates  \cite{gg2,gg3,manth1,manth2,liebe1}}\\
\,  \mbox{Calder\'on-Zygmund estimates}\  (\mbox{see \cite{mincime} for an overview})\\
\,  \mbox{Potential estimates (see Section \ref{potsec2})}
\end{cases}
$$
work more or less exactly as in the linear case. All in all, reducing the nonlinear theory to the linear one in terms of results is a gigantic achievement. These results are actually concerned with the case of equations involving vector fields satisfying \rif{ass1} or \rif{uni-ell0}-\rif{uni-ell}. It is therefore reasonable to argue that, whenever condition \rif{limitata} is satisfied, then such parts of regularity theory would still hold. Surprisingly enough, things turn out to be different. For this we consider the so-called 
\begin{boxy}{The Double Phase Functional} 
\eqn{doppio}
$$
\begin{cases}
\ \displaystyle w\mapsto \mathcal{D}(w, \Omega):= \int_{\Omega}\mathcal H(x,Dw) \dx\\ 
\  \displaystyle \mathcal H(x,z):= |z|^p+\aaa(x)|z|^{q} \\ 
 \ 1 <p \leq  q, \    0 \leq \aaa(\cdot) \in C^{0, \alpha}(\Omega), \  \alpha \in (0,1]
 \end{cases}
$$
\end{boxy}
\noindent with corresponding Euler-Lagrange equation being 
\eqn{doppioeq}
$$
-\divo\, (|Du|^{p-2}Du+(q/p)\aaa(x)|Du|^{q-2}Du)=0\,.
$$
Observe that in this case any local minimizer obviously belongs to $W^{1,p}_{\loc}(\Omega)$, just because it makes the functional $\mathcal D$ locally finite. This kind of structure was introduced by Zhikov \cite{Z2, Z3} in the setting of Homogenization of strongly anisotropic materials \cite{Z1, Z4}, where the geometry of the double composite, with hardening exponents $p$ and $q$, is described by the zero set $\{\aaa(\cdot)=0\}$. A systematic study of qualitative properties of minima of the functional in \rif{doppio} was made in \cite{BCM, CMa1}, after the first higher integrability results were obtained in \cite{sharp}. Here we are mostly interested in the counterexamples developed in \cite{sharp, FMM}. Let us preliminary note that the functional in \rif{doppio} is uniformly elliptic in the sense of Section \ref{unisec}. Indeed, a simple computation reveals that, with respect to the notation in \rif{doppiob}, we have 
$$
\begin{cases}
g_1(x, t) \approx \min\{p-1,1\}t^{p-2} +\min\{q-1,1\}\aaa(x)t^{q-2}\\
g_2(x, t) \approx \max\{p-1,1\}t^{p-2} +\max\{q-1,1\}\aaa(x)t^{q-2}
\end{cases}
$$
so that
\eqn{limitato}
$$
\mathcal{R}_{\mathcal{H}}(x, z) \lesssim \frac{ \max\{p-1,1\}}{\min\{p-1,1\}}+
\frac{ \max\{q-1,1\}}{ \min\{q-1,1\}}\,.
$$
It follows that the equation in \rif{doppioeq} is uniformly elliptic in the classical sense of Section \ref{unisec}, \rif{limitata}. Nevertheless the regularity features \rif{elenco} generally fail for solutions to \rif{doppioeq} and local minimizers of $\mathcal D$ in \rif{doppio}. It is in fact possible to find non-negative H\"older continuous functions $\aaa(\cdot)$ leading to the existence of minimizers that do not even belong to $W^{1,q}_{\loc}$ \cite{sharp} and develop singularities on fractals with maximal Hausdorff dimension $n-p$ \cite{balci2, balci33, FMM}; see also the two-dimensional constructions by Zhikov \cite{Z2, Z3}. This can be realised taking any choice of $(p,q)$ such that 
\eqn{biba}
$$1<p<n< n+\alpha<q\,,$$ with a suitably constructed function $\aaa(\cdot)\in C^{0, \alpha}(\Omega)$. Minima found in \cite{sharp, FMM} are globally bounded. Note that \rif{biba} implies
\eqn{twist}
$$
\frac qp > 1+ \frac{\alpha}{n} \quad \mbox{and} \quad q-p>\alpha\,.
$$
New and interesting counterexamples are in \cite{ balci2, balci33, balci, balci3}. All this happens in presence of classical uniform ellipticity \rif{limitato}, which is somehow counterintuitive according to the narrative drawn at the beginning of this section. The key to understand the failure of properties \rif{elenco} is to introduce another type of ellipticity ratio, which is bound to detect milder forms of nonuniform ellipticity. 
Following \cite{ciccio}, with $a(\cdot)$ being a vector field as considered in Section \ref{unisec}, we define 
\begin{boxy}{The Nonlocal Ellipticity Ratio}
\eqn{nlratio1}
$$\mathcal R_{a}(z,B):= \frac{\sup_{x\in B} g_2(x, |z|)}{\inf_{x\in B} g_1(x, |z|)} \qquad \mbox{for any ball $B \subset \Omega$}$$
\end{boxy}
\noindent or, equivalently, when $\partial_z a(\cdot)$ is symmetric
$$
\mathcal R_{a} (z, B):=  \frac{\sup_{x\in B}\, \mbox{highest eigenvalue of}\ \partial_{z}  a(x,z)}{\inf_{x\in B}\, \mbox{ lowest eigenvalue of}\  \partial_{z} a(x,z)}\,.
$$
Of course, it is 
$$
\mathcal R_{a} (x,z)\leq \mathcal R_{a} (z, B) \quad \mbox{for every $x\in B$}\,.
$$ 
Similar definitions can be given as in Section \ref{unisec} concerning an integrand $F(\cdot)$ rather than a vector field. Needless to say, such different concepts coincide in the autonomous case, i.e., when the vector field and/or the integrand do not exhibit any direct dependence on $x$. We shall then say that a vector field $a\colon \Omega \times \er^n \to \er^n$ is {\em softly nonuniformly elliptic} if it is uniformly elliptic in the classical sense of \rif{limitata} but, for at least one ball $B$, it happens that $ R_{a} (z, B)$ is unbounded, i.e., 
\eqn{limitatasoft}
$$
\sup_{x\in \Omega, |z|\not=0}\mathcal R_{a} (x,z)< \infty \qquad \mbox{and}\qquad \sup_{|z|\not=0}\mathcal R_{a} (z,B)= \infty\,.
$$
Softly nonuniformly elliptic integrands $F(\cdot)$ can be defined in a similar manner; these are such that $\partial_zF(\cdot)$ is softly nonuniformly elliptic. Back to the functional $\mathcal D$ in \rif{doppio}, we observe that
\eqn{limitazione}
$$\mathcal R_{\mathcal H} (z, B)\approx \|\aaa\|_{L^\infty(B)}|z|^{q-p}+1\lesssim |B|^{\alpha/n}|z|^{q-p}+1$$
holds whenever $\{\aaa(\cdot)\equiv 0\}\cap B$ is nonempty. This means that \rif{limitatasoft} occurs provided $\aaa(\cdot)$ does not vanish almost everywhere. Using this second type of ellipticity ratio helps explaining the occurrence of irregularity phenomena in classical uniformly elliptic problems; these are evidently generated by what we now call soft nonuniform ellipticity. Moreover, being a nonlocal quantity, ratios of the type in \rif{nlratio1} are more suitable to control estimates in nonautonomous problems, as these typically involve integral estimates on balls. In this respect, \rif{twist} reveals that the distance between $q$ and $p$ is too large with respect to $\alpha$ and $n$, therefore by \rif{limitazione} the ratio $\mathcal R_{\mathcal H} (z, B)$ grows fast enough to generate irregular minimizers according to the discussion made at the end of Section \ref{unisec}. Due to \rif{twist}, the factor $ |B|^{\alpha/n}$ in \rif{limitazione} is not sufficiently small to rebalance the growth with respect to $|z|$ on small balls $B$. Another instance of the same phenomena is provided by 
\begin{boxy}{The Variable Exponent Functional}
\eqn{variabile}
$$
w \mapsto \mathcal V(w, \Omega) := \int_{\Omega} |Dw|^{\ppp(x)}\dx$$
\end{boxy}
\noindent where $\ppp\colon \Omega \mapsto (1, \infty)$ is a continuous function. 
In this case we have that
$$\mathcal R (x,z)\lesssim \frac{\max\{q-1,1\}}{\min\{p-1, 1\}}, 
\qquad q:= \sup \ppp(x),\  p:=\inf \ppp(x)\,,$$
and therefore the integrand $F(x,z)\equiv  |z|^{\ppp(x)}$ is uniformly elliptic in the classical, pointwise sense of Section \ref{unisec}. On the other hand, with $B\Subset \Omega$ being a ball, we note that 
\eqn{limitazione2}
$$
\begin{cases}
\mathcal R (z,B) \approx  |z|^{\ppp_+(B)-\ppp_-(B)} +1\\
\ppp_+(B) := \sup_{B} \ppp(x), \quad \ppp_-(B) := \inf_{B} \ppp(x)\,.
\end{cases}
$$
We therefore conclude that also the integrand in question is softly nonuniformly elliptic as long as $\ppp(\cdot)$ is a non-constant function. Similarly to the double phase case, counterexamples to regularity of minima emerge when $\omega(\cdot)$, the modulus of continuity of $\ppp(\cdot)$ i.e.,
$$
|\ppp(x)-\ppp(y)|\leq \omega(|x-y|) \qquad \forall \, x, y \in \Omega\,,
$$ fails to satisfy certain explicit decay properties. More precisely, if 
\eqn{limitazione3}
$$
\limsup_{\varrho \to 0} \omega(\varrho) \log \frac{1}{\varrho} =\infty,
$$
then examples of discontinuous minimizers are found by Zhikov \cite{Z2, Z3}; see also \cite{balci}. Indeed, notice that \rif{limitazione2} implies 
$$
\mathcal R (z,B) \approx  |z|^{\omega\left(c|B|^{1/n}\right)} +1\,, \qquad c\equiv c (n)
$$
so that the oscillations of $\ppp(\cdot)$ allow for a fast growth of the ellipticity ratio with respect to the size of the balls.
The situation resembles the one for the double phase functional, where a balance between growth conditions with respect to the gradient and continuity rates of coefficients is necessary and sufficient for regularity of minima.  
\vspace{3mm}
\begin{boxy}{The key takeaway}
There can be different notions of nonuniform ellipticity. Next to the classical, pointwise one, there is a softer notion aimed at explaining, for certain nonautonomous problems, lack of regularity of solutions  although in presence of classical uniform ellipticity. Summarizing, we adopt the following classification: 
\vspace{3mm}\begin{mdframed}
\begin{itemize}
\item Nonuniform ellipticity $\Longleftrightarrow z \mapsto \mathcal R (x,z)$ is unbounded  for at least one $x$.
\item Soft nonuniform ellipticity $\Longleftrightarrow z \mapsto \mathcal R (x,z)$ is uniformly bounded (with respect to $x$) but $z \mapsto \mathcal R (z,B)$ is unbounded for at least one ball $B$.
\item Uniform ellipticity $\Longleftrightarrow z \mapsto \mathcal R (z,B)$ is uniformly bounded for every ball $B$.
\end{itemize}
\end{mdframed}
\end{boxy}

\vspace{3mm}

\section{Regularity and soft nonuniform ellipticity}\label{softy} Looking at the double phase functional $\mathcal D$ the twist is that when \rif{twist} are violated regularity of minima returns and the basic features of uniformly elliptic problems hold. More precisely we have the following two theorems:
\vspace{3mm}\begin{mdframed}
\begin{theorem}[Schauder type \cite{BCM, CMa1, CMa2}]\label{doubleschauder} Let $u\in W^{1,1}_{\loc}(\Omega)$ be a local minimizer of the functional $\mathcal D$ in \trif{doppio}. If either 
\eqn{bound1}
$$
\frac qp \leq 1+ \frac{\alpha}{n}
$$
or
\eqn{twist2}
$$
u \in L^{\infty}_{\loc}(\Omega) \ \ \mbox{and}\ \  q \leq p+\alpha
$$
holds, then $Du$ is locally H\"older continuous in $\Omega$. 
\end{theorem}
\end{mdframed}
\vspace{15mm}
\begin{mdframed}
\begin{theorem}[Calder\'on-Zygmund type \cite{CM, dm0}]\label{czth} Let $u\in W^{1,1}_{\loc}(\Omega)$ be a distributional solution to 
$$
\divo\, (|Du|^{p-2}Du+\aaa(x)|Du|^{q-2}Du)=\divo\, (|F|^{p-2}F+\aaa(x)|F|^{q-2}F)
$$
such that $\mathcal H(\cdot, Du), \mathcal H(\cdot, F)\in L^1_{\loc}(\Omega)$. If \trif{bound1} holds, then 
\eqn{assertion}
$$
\mathcal H(\cdot, F)\in L^\gamma_{\loc}(\Omega) \Longrightarrow \mathcal H(\cdot, Du)\in L^\gamma_{\loc}(\Omega)\quad \mbox{for every}\ 
\gamma \geq 1\;.
$$
\end{theorem}
\end{mdframed}\vspace{3mm}
Note that, when $\aaa(\cdot)\equiv 0$, assertion \rif{assertion} is a classical fact from Nonlinear Calder\'on-Zygmund theory; see \cite{mincime} for a survey and a panorama of results. Furthermore, when $p=2$ from \rif{assertion} we recover a classical linear result of Calder\'on-Zygmund type. Note that the starting assumption $\mathcal H(\cdot, Du)\in L^1_{\loc}(\Omega)$ is necessary already when  $\aaa(\cdot)\equiv 0$, as recently shown in \cite{colombotione}; see the comments after Theorem \ref{c1}. Conditions \rif{twist} reveal that the occurrence of regularity/irregularity is linked to a subtle interaction between the so called gap $q/p$, the rate of H\"older of coefficients $\alpha$ and the ambient dimension $n$, so that a large $\alpha$ is able to compensate the growth of $\mathcal R_{\mathcal H}(z,B)$ with respect to $|z|$, as it is clear from \rif{limitazione}. This time, $ |B|^{\alpha/n}$ is small enough to rebalance the growth with respect to $|z|$ on shrinking balls $B$.

 Similar considerations can be made for so-called multiphase variational integrals, that is, functionals of the type
\eqn{multino}
$$
\begin{cases}
\displaystyle \  w \mapsto \int_{\Omega}\bigl(|Dw|^p + \sum_{k=1}^{m} \aaa_k(x)|Dw|^{q_k}\bigl) \dx\\
\displaystyle \ 0 \leq \aaa_k(\cdot) \in C^{0, \alpha_k}_{\loc}(\Omega), \quad  1 < p \leq q_1 \leq \ldots \leq  q_m
\end{cases}
$$
for arbitrary number $m$ of phases. In this case optimal conditions ensuring regularity are given by
\eqn{boundulti}
$$
\frac {q_k}p \leq 1+ \frac{\alpha_k}{n}
$$
for every $k\leq m$. 
Treated for the first time in \cite{crismulti0}, where this last fact was first proved, multiphase functionals have been already considered at length in the literature, see \cite{byun0, byun1} and related references. In particular, an optimal Calder\'on-Zygmund theory of the type reported in Theorem \ref{czth} has been established in \cite{byun1, crismulti}. Finally, looking at the variable exponent functional $\mathcal V$ in \rif{variabile}, the condition 
\eqn{condlog}
$$
\limsup_{\varrho \to 0} \, \omega(\varrho) \log \frac{1}{\varrho} <\infty\,,
$$
which is complementary to \rif{limitazione3}, 
implies the local H\"older continuity of minima; moreover, assuming that $\omega(\varrho)\lesssim \varrho^{\beta}$ for some $\beta>0$, implies the local gradient H\"older continuity of minima (maximal regularity and Schauder type result). Finally, assuming that the limit in \rif{condlog} vanishes allows to construct a natural, intrinsic nonlinear Calder\'on-Zygmund theory. For an overview of regularity results on the variable exponent functional $\mathcal V$ in \rif{variabile} we refer for instance to \cite[Section 7]{Dark}. 

Condition \rif{condlog} plays a role similar to those that \rif{bound1} and \rif{boundulti} have in the case of double and multiphase integrals, respectively. Actually, all these situations can be unified. This has been done in a recent series of papers by H\"asto \& Ok \cite{HO1, HO2, HO3}. The authors consider general softly nonuniformly elliptic functionals and equations and prove local H\"older continuity of the gradients under a general condition on the regularity of the partial map $x\mapsto a(x, z)$. This condition gives back those known for the functionals $\mathcal D$ and $\mathcal V$.

The constructions implying the existence of irregular minimizers found in \cite{sharp, FMM, Z2, Z3} are linked to the occurrence of the so-called Lavrentiev phenomenon. This means that, for a suitable boundary datum $u_0 \in W^{1,\infty}(B_1(0))$
\begin{boxy}{The Lavrentiev Phenomenon}
\begin{flalign}
\notag &\inf_{w\in u_0+W^{1,p}_0(B_1(0))} \mathcal D(w, B_1(0)) \\
& \qquad \qquad < \inf_{w\in u_0+W^{1,p}_0(B_1(0))\cap W^{1,q}_{\loc}(B_1(0))} \mathcal D(w,B_1(0)) \label{lavvi}
\end{flalign}
\end{boxy}
\noindent holds, which is an obvious obstruction to regularity of minima. In fact, \rif{lavvi} implies that minima cannot belong to $W^{1,q}$ locally, so that even the basic bootstrap of regularity $W^{1,p}_{\loc} \to W^{1,q}_{\loc}$, fails. An interesting catch, pointed out in the approach originally introduced in \cite{sharp}, is that the very same condition on $q/p$ implying regularity of minima, that is \rif{bound1}, implies the absence of Lavrentiev phenomenon for functionals 
\eqn{consider}
$$
\begin{cases}
\,  \displaystyle w \mapsto   \int_{\Omega}  F(x,Dw)\, dx\\
\, F(x, z) \approx \mathcal H(x,z) +1
\end{cases}
$$
without any further assumption on $F(\cdot)$. Not even convexity $z \mapsto F(\cdot, z)$  is necessary assuming \rif{consider}$_2$. Specifically, in \cite{sharp} it is proved that for every ball $B\Subset \Omega$ and $w\in W^{1,1}_{\loc}(\Omega)$ such that $F(\cdot, Dw)\in L^{1}_{\loc}(\Omega)$, it happens that
\begin{boxy}{Approximation-in-energy}
\eqn{energyconv}
$$
\begin{cases}
\, \mbox{there exists $\{w_k\}_{k}\subset C^{\infty}(\Omega)$ such that $w_k\to w$ strongly in $W^{1,p}(B)$}\\
\, \mbox{and $F(\cdot,Dw_k)\to F(\cdot, Dw)$ in $L^1(B)$}
\end{cases}
$$
\end{boxy}
\noindent provided \rif{bound1} and \rif{consider} hold. In fact, somehow reversing the arrow, when considering a functional as in \rif{consider}$_1$, not necessarily satisfying \rif{consider}$_2$, assuming \rif{energyconv} allows to prove higher integrability of minima under suitable $(p,q)$-convexity assumptions on the integrand $F(\cdot)$ provided \rif{bound1} takes place. See \cite{sharp} for details. 

After the initial approach introduced in \cite{sharp}, more connections between the occurrence of the approximation-in-energy property \rif{energyconv} and regularity of minima have been established, together with more results on the absence of Lavrentiev phenomenon \cite{BCM, buli, hasto1, hasto3, hasto4, HO1, HO2, HO3, koch1, koch2}. The analysis usually passes through the analysis of the Lavrentiev gap functional, see Remark \ref{lavgap}. For such topics see also the recent survey \cite{mira}. For connections to related functions spaces and their abstract properties see the monograph \cite{hasto2}. 
 
In the last years there has been a large interest in regularity for soft uniformly ellptic problems, especially when considering double phase operators. Beyond the H\"asto \& Ok papers already cited above, we mention \cite{byun2, byun3} for different extensions. An interesting direction, which is certainly worth developing, concerns the case of manifold valued problems, considered for instance in \cite{chdeko, crispx, DM}. In particular, in \cite{DM}, the problem to minimize the double phase functional $\mathcal D$ amongst maps with values into the sphere has been considered. Partial regularity results and estimates of the size of the singular sets have been obtained using certain intrinsic Hausdorff measures. These measures have eventually proven to be relevant in various issues, as for instance removability of singularities \cite{chde} and in symmetry problems \cite{biagi1, biagi2}. Double phase degeneracies also appear in the setting of fully nonlinear problems, as first considered in \cite{crispro} and then in \cite{dasilva1, dasilva2, dasilva3, crisintner}. After the initial contribution in \cite{depala}, nonlocal double phase problems were considered in \cite{byun6, byun7, mengesha}.

\begin{boxy}{The key takeaway}
For softly uniformly elliptic problems with polynomial growth of the eigenvalues, a delicate balance between the growth with respect to the gradient and the modulus of continuity of coefficients is necessary and sufficient to guarantee regularity of solutions. 
\end{boxy}


\section{Nonlinear potentials and a priori estimates}
\subsection{Potentials, Lorentz spaces and iterations}
Let us recall a few basic definitions concerning linear and nonlinear potentials. 
\begin{definition}[Riesz and Havin-Mazya-Wolff potentials]\label{defi-riesz}
Let $\mu$ be a Borel measure with (locally) finite total mass defined on the open subset $\Omega\subset \er^n$, $n\geq 2$, and let $B_r(x_0)\Subset \Omega$ be a ball.  
\begin{itemize}
\item
The (truncated) Riesz potential ${\bf I}_{\beta}^\mu$ is defined by
$$
{\bf I}_{\beta}^\mu(x_0,r):= \int_0^r \frac{|\mu|(B_{\varrho}(x_0))}{\varrho^{n-\beta}}\, \frac{d\varrho}{\varrho}\,, \qquad \beta >0\,.
$$
\item The (nonlinear) Havin-Mazya-Wolff potential ${\bf W}^{\mu}_{\beta, p}$ is defined by
$$
{\bf W}^{\mu}_{\beta, p}(x_0,r):= \int_0^r \left(\frac{|\mu|(B_{\varrho}(x_0))}{\varrho^{n-\beta p}}\right)^{1/(p-1)}\, \frac{d\varrho}{\varrho}\,, \qquad  \beta >0, \, p>1\,.
$$
\end{itemize}
\end{definition}
The relation between truncated Riesz potentials and the classical ones is very simple, i.e.,  
\eqn{recall2}
$$ {\bf I}_{\beta}^\mu(x_0,r)  \leq c(n) I_\beta(|\mu|)(x_0)=c(n)\int_{\er^n}\frac{d |\mu|(x)}{|x-x_0|^{n-\beta}}\quad \mbox{for every}\ r>0\;.$$ 
Riesz potentials naturally occur when dealing with linear equations via fundamental solutions.  On the other hand, starting by the fundamental work in \cite{maz, MH, HW}, Havin-Mazya-Wolff potentials naturally intervene when passing to nonlinear equations of $p$-Laplacean type, that is, equations as \rif{risolvi} under assumptions \rif{ass1}, equation \rif{chiefp} being a chief model example. They can be used to control the pointwise behaviour of Sobolev functions. Moreover, a $p$-Laplacean version of the Wiener criterion can be derived using such potentials \cite{maz, KiL0, KiL, MZ}. The fine existence theory for non-homogeneous equations makes use of nonlinear potentials \cite{PV1, PV2}. Of course Riesz potentials can be realised as nonlinear potentials with a suitable choice of the parameters, i.e.,$ {\bf W}^{\mu}_{1, 2}={\bf I}^{\mu}_{2}$ and ${\bf W}^{\mu}_{1/2, 2}={\bf I}^{\mu}_{1}$. 

Following \cite{piovra}, in the case $\mu$ is an integrable function we also have
\begin{definition} With $t,\delta>0$, $m, \theta\geq 0$, $\mu\in L^{1}(B_{r}(x_{0}))$ such that $|\mu|^{m} \in L^{1}(B_{r}(x_0))$, and with $B_{r}(x_0)\subset \er^n$, we define 
$$
{\bf P}_{t,\delta}^{m,\theta}(\mu;x_0,r) := \int_0^r \varrho^{\delta} \left(  \mint_{B_{\varrho}(x_0)}|\mu|^{m} \dx \right)^{\theta/t} \frac{\d\varrho}{\varrho} \,.
$$
\end{definition}
Again, these potentials incorporate Havin-Mazya-Wolff potentials, and therefore Riesz potentials, in the sense
$$
{\bf P}_{p-1,\frac{p\beta}{p-1}}^{1,1}(\mu;x_0,r)\equiv_{n,p}{\bf W}^{\mu}_{\beta, p}(x_0,r)\,.
$$
The behaviour of these nonlinear potentials with respect to Lorentz spaces is of interest here. The usual definition of Lorentz space $L(t,\gamma)(\Omega)$, with $t,\gamma\in (0,\infty)$, prescribes that, if $\mu\colon \Omega \to \er$ is a measurable function, then $\mu \in L(t,\gamma)(\Omega)$ iff 
\eqn{verydefi}
$$
\|\mu\|_{t, \gamma,\Omega} := \left(t\int_0^\infty (\lambda^t|\{x \in \Omega \, : \, |\mu(x)|> \lambda\}|)^{\gamma/t}\, \frac{d\lambda}{\lambda}\right)^{1/\gamma}<\infty \,.
$$
See \cite{gra, oneil, sw} for basic properties of such spaces. Lorentz spaces extend Lebesgue spaces and the second index tunes the first one in the sense that
 \eqn{lorentzbasic}
$$
\begin{cases}
\, \mbox{$L(t_{1}, \gamma_{1})(\Omega)\subset L(t_{2}, \gamma_{2})(\Omega)$ for all $0<t_{2}<t_{1}<\infty$, $\gamma_{1},\gamma_{2}\in (0,\infty]$}\\
\, \mbox{$L(t,\gamma_{1})(\Omega)\subset L(t,\gamma_{2})(\Omega)$ for all $t\in (0,\infty)$, $0<\gamma_{1}\leq \gamma_{2}\le\infty$}\\
\,  \mbox{$L(t,t)(\Omega)=L^{t}(\Omega)$ for all $t >0$}
 \end{cases}
$$
hold with continuous inclusions. The catch with nonlinear potentials is now given by the following lemma, whose proof can be found in \cite[Lemma 4.1]{piovra}. 
\begin{lemma}\label{crit}
Let $n \geq 2$, $t,\delta,\theta>0$ be numbers such that
\eqn{lo.2.1} 
$$
\frac{n\theta }{t\delta}>1\,.
$$
Let $B_{\tau_1}\Subset B_{\tau_1+r_{0}}\subset \mathbb{R}^{n}$ be two concentric balls with $\tau_1, r_0\leq 1$, and let $\mu\in L^{1}(B_{\tau_1+r_{0}})$ be such that $|\mu|^{m}\in L^{1}(B_{\tau_0+r_{0}})$, where $m>0$. Then
\begin{flalign}
\notag \|\mathbf{P}^{m,\theta}_{t,\delta}(\mu;\cdot,r_{0})\|_{L^{\infty}(B_{\tau_1})} &\le \tilde c\|\mu\|_{\frac{mn\theta}{t\delta},\frac{m\theta}{t};B_{\tau_1+r_{0}}}^{\frac{m\theta}{t}}\\
&  \leq c(\eps)\tilde c\|\mu\|_{L^{\frac{(1+\eps)mn\theta}{t\delta}}(B_{\tau_1+r_{0}})}^{\frac{m\theta}{t}}\label{11}
\end{flalign}
holds for every $\eps>0$, with $\tilde c\equiv \tilde c(n,t,\delta,\theta)$. 
\end{lemma}
By looking at \rif{11} the reader might be slightly surprised by the fact that \rif{11} covers the case $mn\theta/(t\delta)<1$. This looks somehow in contradiction with the usual embedding properties of potentials, maximal operators and so forth. This is only apparent. In fact, in Lemma \ref{crit} one can always reduce to the case  $m=1$, so that, by \rif{lo.2.1}, it is $mn\theta/(t\delta)=n\theta/(t\delta)>1$. Indeed, notice that 
$$
\mathbf{P}^{m,\theta}_{t,\delta}(\mu;\cdot,r_{0})=\mathbf{P}^{1,\theta}_{t,\delta}(|\mu|^m;\cdot,r_{0})\,.
$$
Then, using Lemma \ref{crit} with $m=1$, we find
\begin{flalign*}
\|\mathbf{P}^{m,\theta}_{t,\delta}(\mu;\cdot,r_{0})\|_{L^{\infty}(B_{\tau_1})} & = 
\|\mathbf{P}^{1,\theta}_{t,\delta}(|\mu|^m;\cdot,r_{0})\|_{L^{\infty}(B_{\tau_1})}\\
& \leq \tilde c \||\mu|^m\|_{\frac{n\theta}{t\delta},\frac{\theta}{t};B_{\tau_1+r_{0}}}^{\frac{\theta}{t}}\\
& =\tilde c \|\mu\|_{\frac{mn\theta}{t\delta},\frac{m\theta}{t};B_{\tau_1+r_{0}}}^{\frac{m\theta}{t}}\,.
\end{flalign*}
Note that in the last line we have used the very definition in \rif{verydefi}. 
As we shall see in Section \ref{rinormalizza} below, nonlinear potentials can be used to get a priori estimates for solutions to nonlinear elliptic equations and minima of variational integrals. A key is the following lemma. It contains a pointwise version De Giorgi's geometric iteration \cite{DG1}; its origins can be traced back in the seminal work of Kilpel\"ainen \& Mal\'y \cite{KiL} and the proof can again be found in \cite[Lemma 4.2]{piovra}. 
\begin{lemma}[Quantitative De Giorgi]\label{revlem}
Let $B_{r_{0}}(x_{0})\subset \mathbb{R}^{n}$ be a ball, $n\ge 2$, and consider functions $\mu_i$, 
$
|\mu_i|^{m_i} \in L^1(B_{2r_0}(x_{0}))$, and constants $$\chi >1, \quad t \geq 1, \quad \delta_i, m_i, \theta_i>0, \quad c_*,M_0 >0, \quad \kappa_0, M_i\geq 0,$$ 
where $i \in \{1, \ldots, h\}$, $h \in \en$. Assume that $v \in L^t(B_{r_0}(x_0))$ is such that for all $\kk\ge \kk_{0}$, and for every ball $B_{\rr}(x_{0})\subset B_{r_{0}}(x_{0})$, the inequality
\begin{flalign}\label{rev}
\notag \left(\mint_{B_{\rr/2}(x_{0})}(v-\kk)_{+}^{t\chi}  \dx\right)^{1/\chi}&\le c_{*}M_{0}^{t}\mint_{B_{\rr}(x_{0})}(v-\kk)_{+}^{t}  \dx\\  &
\quad\ \ +c_{*}\sum_{i=1}^{h} M_{i}^{t}\rr^{t\delta_{i}}\left(\mint_{B_{\rr}(x_{0})}|\mu_i|^{m_i}  \dx\right)^{\theta_i}
\end{flalign}
holds, where we denote, as usual, $ (v-\kk)_{+}:=\max\{v-\kk,0\}$. If $x_{0}$ is a Lebesgue point of $v$ in the sense that 
$$
\lim_{\rr\to 0} (v)_{B_{\rr}(x_0)}=v(x_0)\,,
$$ then
\begin{flalign}
\nonumber
 v(x_{0}) & \le\kk_{0}+cM_{0}^{\frac{\chi}{\chi-1}}\left(\mint_{B_{r_{0}}(x_{0})}(v-\kk_{0})_{+}^{t}  \dx\right)^{1/t}\\
 & \quad 
\ \ +cM_{0}^{\frac{1}{\chi-1}}\sum_{i=1}^{h} M_{i}\mathbf{P}^{m_i,\theta_i}_{t,\delta_{i}}(\mu_i;x_{0},2r_{0})\label{rev3388}
\end{flalign}
holds with $c\equiv c(n,\chi,\delta_i,\theta_i,c_{*})$.  
\end{lemma}
\begin{remark}
As we are going to see later on in Section \ref{rinormalizza}, the crucial point in the above lemma, when passing from \rif{rev} to \rif{rev3388}, is the quantified and explicit dependence on the constants $M_i$. We shall actually use Lemma \ref{revlem} with $h=1$, i.e., only one potential $\mathbf{P}^{m_i,\theta_i}_{t,\delta_{i}}$ will appear, but we prefer to report it in full generality. 
\end{remark}
Lemma \ref{revlem} has a Moser iteration counterpart, which is the next one. This time potentials are not involved. The emphasis is again on the precise dependence on the multiple $M_0$.
\begin{lemma}[Quantitative Moser]\label{lamoser} Let $B_{r_0}\subset \er^n$ be a ball and let $v\in L^{p/2}(B_{r_0})$, $p>1$, be a non-negative function such that
$$
\left(\int_{B_{\varrho_1}}v^{\left(\gamma+p/2\right)\chi}  \dx\right)^{1/\chi}
 \le \frac{M_0(1+\gamma)^t}{(\varrho_2-\varrho_1)^{t_*}}\int_{B_{\varrho_2}}v^{\gamma+p/2}  \dx
$$
holds for every $\gamma \geq 0$, where $M_0, t, t_*, \chi$ are positive constants with $\chi>1$, and where $B_{\tau_1} \Subset B_{\varrho_1} \Subset B_{\varrho_2} \Subset B_{\tau_2}\subset B_{r_0}$ are arbitrary concentric balls. Then it holds that
$$
\|v\|_{L^{\infty}(B_{\tau_1})} \lesssim_{\chi,t, t_*}\left[\frac{M_0}{(\tau_2-\tau_1)^{t_*}}\right]^{\frac{2}{p}\frac{\chi}{\chi-1}} \|v\|_{L^{p/2}(B_{\tau_2})}\;.
$$
\end{lemma} 
The proof can be found in \cite[Lemma 6.1]{dm} and it is nothing but the usual Moser iteration argument with a precise tracking of the dependence of the constants $M_0$.

\subsection{Gradient potential estimates in the uniformly elliptic case}\label{potsec2}
In the past decade there has been an intensive research activity in the field of Nonlinear Potential Theory. In particular, pointwise gradient estimates via Riesz potentials for solutions to $p$-Laplacean type equations
$$
-\divo\, a (Du)=\mu \ (\mbox{$=$ Borel measure with finite total mass})
$$
under assumptions \rif{ass1} have been discovered. It indeed holds
\eqn{rieszpot}
$$
|Du(x_0)|^{p-1} \lesssim_{\data} {\bf I}_{1}^\mu(x_0,r) + \left(\mint_{B_r(x_0)}(|Du|^2+s^2)^{t/2}\dx \right)^{\frac{p-1}{t}}
$$
\noindent for a.e. $x_0\in \Omega$ and $B_{r}(x_0)\Subset \Omega$, where $t$ depends on $p$ and $t=1$ for $p>2-1/n$. Note that letting $r \to \infty$ in \rif{rieszpot}, recalling \rif{recall2}, and assuming suitable decay properties of $Du$ at infinity yields
\begin{boxy}{The Gradient Potential Estimate}
$$
-\triangle_pu=\mu \Longrightarrow |Du(x_0)|^{p-1} \lesssim_{\data} \int_{\er^n}\frac{d |\mu|(x)}{|x-x_0|^{n-1}}\,.
$$
\end{boxy}
This is the same pointwise estimate that holds for the Poisson equation $-\Delta u=\mu$, apart from the scaling factor $p-1$.   
This result has been obtained first in \cite{mis3} when $p=2$ and then in \cite{DMjfa, KMli, nphuc2} for the case $p>1$. See \cite{DMam} for earlier gradient estimates via Wolff potentials. Inequality in \rif{rieszpot} holds once suitable notion of solutions to measure data problems are adopted as, in general, measure data problems do not support energy solutions (think of the fundamental solution of the Laplacean). The vectorial case for the $p$-Laplacean system and $p\geq 2$ is treated in \cite{KMvec}, while parabolic cases are in \cite{KMwolff, KMpisa, KMri, dong} after the initial results from \cite{DMam} valid in the case $p=2$. For an overview on nonlinear potential estimates we refer to \cite{KMguide}. Further results can be found in \cite{nphuc1, nphuc3}. For an interesting global counterpart of \rif{rieszpot} involving rearrangements see \cite{CMjems}. Estimates employing Havin-Mazya-Wolff potentials, this time for $u$ rather than $Du$, can be found in the groundbreaking work of Kilpel\"ainen \& Mal\'y \cite{KiL0, KiL}, see also \cite{truwa} and again \cite{DMam}. All such estimates can be also considered a part of what is called Nonlinear Calder\'on-Zygmund theory, for which we refer to the past CIME notes \cite{mincime}. For an overview of the basic facts concerning Nonlinear Potential Theory and use of nonlinear potentials in the fine analysis of solutions to elliptic PDEs we refer to the classical treatises \cite{AdHe, HKM, MZ}. 

Size bounds often come along with companion continuity results. In this case, as it is shown in \cite{KMli, KMvec, KMguide, nphuc2, nphuc3}, it holds that
\begin{boxy}{The Gradient Continuity Criterion}
\begin{flalign}
\notag &\lim_{r\to 0} {\bf I}_{1}^{\, \divo\, a (Du)}(x,r)=0 \ \mbox{uniformly w.r.t. $x$} \\
&\qquad \qquad\qquad \qquad  \qquad \Longrightarrow 
\mbox{$Du$ is continuous}\,. \label{seconda}
\end{flalign}
\end{boxy}
Estimate \rif{rieszpot} immediately leads to a Lorentz space criterion for Lipschitz continuity of solutions. Indeed, notice that
$
{\bf P}_{1,1}^{1,1}\equiv  {\bf I}^{\mu}_{1}, 
$
so that \rif{11} implies 
\eqn{prima0}
$$
\mu \in L(n,1) \Longrightarrow {\bf I}^{\mu}_{1}\in L^{\infty}
$$
locally. From this we find that $Du$ is locally bounded via \rif{rieszpot}.
\begin{boxy}{The key takeaway}
Pointwise gradient  bounds via Riesz potentials, typical of linear equations, actually hold for nonlinear problems too. No use of fundamental solutions is possible in this case and methods of proof are intrinsic and nonlinear in nature. Technical tools like De Giorgi type iterations can be framed in the setting of Nonlinear Potential Theory via the use of suitable nonlinear potentials. 
\end{boxy}
\section{Model Lipschitz results in the autonomous case} \label{gapbounds}
Here we report a few sample yet significant regularity results for minima of nonuniformly, autonomous elliptic integrals of the type 
\eqn{aut}
$$
 w \mapsto  \int_{\Omega}  F(Dw)\, dx
$$ 
under the assumptions
\eqn{assbase}
$$
\begin{cases}
\nu [H_{1}(z)]^{p/2}\le F(z)\le L [H_{1}(z)]^{q/2}\\
\nu [H_{1}(z)]^{(p-2)/2}|\xi|^{2}\le \partial_{zz}F(z)\, \xi\cdot\xi\\ 
|\partial_{zz}F(z)|\le L [H_{1}(z)]^{(q-2)/2} \;,
 \end{cases}\qquad 1 < p \leq q
 $$
satisfied for every $z, \xi \in \er^n$. The integrand $F \colon \er^n\to [0, \infty)$ is assumed to be $C^2$-regular. Integrands and functionals satisfying \rif{assbase}$_1$ are nowadays known as functionals with $(p,q)$-growth conditions, a terminology introduced by Marcellini in \cite{M2}. Conditions \rif{assbase} imply the following control or the ellipticity ratio of $F(\cdot)$:
\eqn{boundy}
$$
\mathcal R_{F}(z) \lesssim_{\data} |z|^{q-p}+1\,.
$$
Similarly to what we have seen in Section \ref{softy}, the crucial assumption to get regularity of minimizers is a bound of the type
\eqn{ilbound}
$$
\frac qp < 1 + {\rm o}(n)\,,
$$ 
with
\eqn{asymp}
$${\rm o}(n)\approx \frac 1n$$ for $n$ large. As already explained, bounds as \rif{ilbound} serve to limitate the growth of the ellipticity ratio $\mathcal R_{F}(\cdot)$ with respect to the gradient variable when proving gradient boundedness of minima. Such bounds are necessary already in the autonomous case, as shown by explicit counterexamples \cite{gia, M2}. A model regularity result is the following:
\vspace{3mm}\begin{mdframed}
\begin{theorem}[Marcellini \cite{M2}]\label{marcth}
Let $u\in W^{1,1}_{\loc}(\Omega)$ be a minimizer of the functional in \eqref{aut} under assumptions 
\trif{assbase} and 
\eqn{boundma}
$$
\frac qp < 1 + \frac{2}{n}\,.
$$
Then $Du$ is locally bounded in $\Omega$. 
\end{theorem}
\end{mdframed}\vspace{3mm}
While in nonautonomous case the optimal bound on $q/p$ is known to be \rif{bound1} by the counterexamples in \cite{sharp, FMM}, for autonomous functionals \rif{aut} the optimal bound is still unknown. An example of progress in this direction is in the following:
\vspace{3mm}\begin{mdframed}
\begin{theorem}[Bella \& Sch\"affner \cite{BS, BS2, schag}]\label{bstheorem} 
Let $u\in W^{1,1}_{\loc}(\Omega)$ be a minimizer of the functional in \eqref{aut} under assumptions 
\trif{assbase} and 
\eqn{boundbs}
$$
\frac qp < 1 + \frac{2}{n-1} \;.
$$
Then $Du$ is locally bounded in $\Omega$. 
\end{theorem}
\end{mdframed}\vspace{3mm}
The one in \rif{boundbs} is the best bound found on $q/p$ up to now. Advances are in \cite{CKP, legendo}, and are obtained under special structure conditions. The bound in \rif{boundbs} works in the vectorial case too, where it allows to prove partial regularity \cite{sha}. 

Bounds of the type in \rif{ilbound} also appear when proving lower order regularity. In this case the presence of coefficients plays no role. Indeed we have 
\begin{mdframed}
\begin{theorem}[Hirsch \& Sch\"affner \cite{HS}]\label{hstheorem} Let $u\in W^{1,1}_{\loc}(\Omega)$ be minimizer of the functional in \trif{genF} with $\mu=0$, where $F\colon\Omega \times \er^n \to \er$ is Carath\'eodory regular and satisfies
$$
\left\{
\begin{array}{c}
\nu|z|^p \leq  F(x,z) \leq L|z|^q +1\\[3 pt]
 F(x,2z) \leq LF(x, z)+L\,,
 \end{array}
 \right.
$$
for $(x, z)\in \Omega \times \er^n$, and where $1<p\leq q$.
Assume that
\eqn{hirsch}
$$
\frac 1p - \frac 1q \leq \frac 1{n-1}\;.
$$ 
Then $u$ is locally bounded in $\Omega$. 
\end{theorem} 
\end{mdframed}\vspace{3mm}
The bound in \rif{hirsch} is sharp for $L^\infty$-regularity in view of the available counterexamples \cite{gia, M2}. As in the standard case $p=q$ \cite{gg1}, this time no regularity on $F(\cdot)$ is required. In particular, no convexity of $z \mapsto F(\cdot,z)$ is needed and not even continuity is required on $x \mapsto F(x,\cdot)$. 

Finally, we mention that, as it already happened in \rif{twist2}, when starting by a solution which is already more regular, then bounds for Lipschitz continuity can be improved by using underlying interpolation effects. For instance, in perfect analogy with \rif{twist2} we have
\vspace{3mm}\begin{mdframed}
\begin{theorem}[Choe \cite{choe}]\label{choetheorem} 
Let $u\in W^{1,1}_{\loc}(\Omega)\cap L^{\infty}_{\loc}(\Omega)$ 
be a minimizer of the functional in \eqref{aut} under assumptions \trif{assbase} and 
\eqn{bouc1}
$$
1< p \leq q < p+1\,.
$$
Then $Du$ is locally bounded in $\Omega$. 
\end{theorem}
\end{mdframed}\vspace{3mm}
Further results using non-dimensional bounds as \rif{bouc1} can be found for instance in \cite{CKP0}.
 
Functionals and equations with $(p,q)$-growth nonuniform ellipticity have been the object of intensive investigation over the last decades. We refer to the surveys \cite{M6, M7, M8, Dark, mira} for an overview. More recently, nonstandard growth conditions have been examined in the setting of nonlocal and mixed operators, see for instance \cite{byun5, byun6, chaker1, chaker2, miscelato, depala, bonder, OK}. We forecast several further developments in this direction. 
\begin{boxy}{The key takeaway}
Solutions to nonuniformly elliptic problems with polynomial $(p,q)$-growth of the eigenvalues are still regular provided the gap $q/p$ is not to far from $1$. This means that the ellipticity ratio cannot grow too fast with respect to the gradient, as framed in \rif{boundy}. While the asymptotic in \rif{ilbound}-\rif{asymp} is known to be optimal, the sharp bound on $q/p$ for Lipschitz regularity in the autonomous case remains unknown. 
\end{boxy}
\section{Nonuniformly elliptic Schauder estimates \cite{piovra}}\label{piovra1}
In treating Schauder type estimates we shall first present the results in \cite{piovra}, in this section. We shall eventually consider the more recent and improved ones in \cite{gioiello} in the next section. The reason for this is that we shall eventually give a sketch of the proofs in Section \ref{fracschau}. This is easier in the case of the results in \cite{piovra} but becomes less straightforward in the case of those from \cite{gioiello}, which features more technically involved proofs. 
\subsection{Hopf, Caccioppoli and Schauder, reloaded} So-called Schauder estimates for linear elliptic equations
\eqn{lineare}
$$
-\divo\, (\texttt{A}(x)Du)=0\,, \qquad    \texttt{A}(\cdot) \approx  \mathbb{I}_{\rm d}
$$
are a basic tool in PDE theory and they play a role in a variety of situations. For instance, they are fundamental ingredient in the proof of higher regularity of solutions to nonlinear elliptic equations. When referring to \rif{lineare}, the basic assertion is
\eqn{variante}
$$
\texttt{A}(\cdot) \in C^{0, \alpha} \Longrightarrow Du \in C^{0, \alpha}
$$
in local and/or global fashion, eventually depending on boundary conditions. In other words, the gradient of solutions inherits the regularity of coefficients, as long as the operator without coefficients allows (this is the case of ``frozen coefficients"). Several variants
of \rif{variante} are possible, for instance when considering linear equations in non-divergence form, which is in fact the most classical case. Schauder estimates are actually a classic achievement of Hopf \cite{hopf}, Caccioppoli \cite{cacc1} and Schauder \cite{js, js2}; see also the work of Giraud \cite{gir1, gir2}. Modern proofs can be found in \cite{camp, gt, simon2, trusc}. 

After the classical, linear era, nonlinear versions of Schauder estimates, i.e., for equations of the type \rif{risolvi}, were established in the uniformly elliptic case -- see Giaquinta \& Giusti's \cite{gg2, gg3}, Manfredi's \cite{manth1, manth2} and Lieberman's \cite{liebe0} papers for full generality. A model result can be obtained by looking at solutions to
$$
-\divo\, (\ccc(x)|Du|^{p-2}Du)=0 \,, \quad 0< \nu \leq \ccc(\cdot) \in C^{0, \alpha}
$$
for which \rif{variante} holds locally as long as $\alpha < \alpha_0$ and solutions to \rif{chiefp}
are $C^{1, \alpha_0}$ (see for instance \cite{manth1, KMguide} for the precise meaning of this assertion). 

All the known proofs of Schauder estimates rely on a two-step argument based on comparing the solution to the original problem with the solutions of problems with constant coefficients, that typically enjoy good a priori estimates. Combining these with suitable comparison estimates, and iterating, gradient H\"older continuity follows. The crucial point in this procedure is that all the estimates involved are homogeneous, and this allows to combine and iterate them. In turn, this is a direct consequence the uniform ellipticity of the equations considered. Observe that this scheme is common to both the linear and the nonlinear case. 

On the other hand, when turning to the nonuniformly elliptic case, estimates are typically affected by a lack of homogeneity. An example is obviously given by functionals with $(p,q)$-growth, where two different exponents occur, generating an obvious lack of scaling. This is essentially the reason for which the problem of establishing nonlinear Schauder estimates has remained open for a long time, finally getting an answer in \cite{piovra}. Here we present a few facts from this last paper. We shall actually include two model results for the sake of readability; the interested reader will find more in \cite{piovra}. We start by the simplest possible variational case, that nevertheless already contains all the core difficulties and peculiarities of more general ones. We consider
\eqn{modello0}
$$
\begin{cases}
\, \displaystyle w \mapsto \int_{\Omega}\ccc(x)F(Dw) \dx, \quad 0< \nu \leq \ccc(\cdot) \leq L\\
\,  |\ccc(x_1)-\ccc(x_2)| \leq L|x_1-x_2|^{\alpha}, \quad  \alpha \in (0,1],
 \end{cases}
$$
for every choice of $ x_1, x_2\in \Omega$, where $F(\cdot)$ satisfies
\begin{flalign}\label{assFF}
\qquad\begin{cases}
\,F(\cdot) \in C^1(\er^n)\cap C^2(\er^n\setminus\{0\})\\
\,\nu[H_{s}(z)]^{p/2}\le F(z)\le L[H_{s}(z)]^{q/2}+L[H_{s}(z)]^{p/2}\\
\, \nu[H_{s}(z)]^{(p-2)/2}|\xi|^{2}\le \partial_{zz}F(z)\xi\cdot \xi \\
\, |\partial_{zz}F(z)|\le L[H_{s}(z)]^{(q-2)/2}+L[H_{s}(z)]^{(p-2)/2}
\end{cases}
\end{flalign}
for all $z, \xi \in \mathbb{R}^{n}$, $|z|\not=0$, $\xi\in \mathbb{R}^{n}$, and $H_{s}(\cdot)$ is defined in \rif{defiH}. 
\vspace{3mm}\begin{mdframed}
\begin{theorem}\label{t2}
Let $u\in W^{1,1}_{\loc}(\Omega)$ be a minimizer of the functional in \trif{modello0}, under assumptions \eqref{assFF}. If 
\eqn{bound333}
$$
\frac{q}{p}\leq  1+  \frac{1}{5}\left(\frac{\alpha}{n}\right)^2\,,
$$
then $Du$ is locally H\"older continuous in $\Omega$. Moreover
\eqn{stimascha}
$$
\|Du\|_{L^{\infty}(B_{r/2})}\lesssim_{\data, \alpha} \left( \mint_{B_r} F(Dw) \dx \right)^{\kk}+1
$$
holds whenever $  B_{r} \Subset \Omega$ is a ball with $r\leq 1$, where $ \kk\equiv \kk(n,p,q,\alpha)\geq 1$. 
\end{theorem}
\end{mdframed}\vspace{3mm}
When the problem becomes non-degenerate, that is, when $s >0$, we can then recover a result as in \rif{variante}. 
\vspace{3mm}\begin{mdframed}
\begin{theorem}\label{c1}
In the setting of Theorem \ref{t2} with $\alpha \in (0,1)$, assume also that $p\geq 2$, $s>0$, and that $\partial_{zz}F(\cdot)$ is continuous. Then $
u\in C^{1, \alpha}_{\loc}(\Omega).
$ 
\end{theorem}
\end{mdframed}\vspace{3mm}
The variational case is particularly natural to start with as it provides a natural setting for nonuniformly elliptic problems, allowing to put all the emphasis on a priori regularity estimates. In fact, it dispenses us to distinguish from various definitions of solutions. Minimality automatically selects the right solutions as it involves a definition that works against all possible competitors (variations). The situation changes when considering equations, and already in the uniformly elliptic case. Indeed, by considering the distributional version of \rif{chiefp} 
\eqn{ppllaa}
$$
\int_{\Omega} |Du|^{p-2}Du\cdot D\varphi \dx =0 \qquad \forall\ \varphi \in C^{\infty}_0(\Omega)
$$
we notice that this makes actually sense whenever $u\in W^{1,p-1}_{\loc}(\Omega)$; we are therefore in the realm of Definition \ref{defisol0}. Such solutions, not necessarily belonging to the natural space $u\in W^{1,p}_{\loc}(\Omega)$ as the standard energy ones from Definition \ref{defisol}, are called very weak solutions. Although some authors firmly believed of the contrary \cite{is}, in general very weak solutions to \rif{chiefp} do not belong $u\in W^{1,p}_{\loc}(\Omega)$  \cite{colombotione}. Note that energy solutions are those allowing to test \rif{ppllaa} with functions $\varphi$ that are proportional to the solution itself. This is crucial in order to get regularity estimates, even basic ones such as Caccioppoli inequalities.   

This very basic issue naturally reappears in the nonuniformly elliptic setting. As an example, we consider an equation of the type
\eqn{risolvi0}
$$
-\divo\, a(x,Du) =0 \qquad \mbox{in}\  \Omega \subset \er^n\,.
$$
Here the vector field $a\colon \Omega \times \er^n \to \er^n$ is of class $C^1(\er^n\setminus\{0_{\er^n}\})$ with respect to gradient variable, and satisfies
\begin{flalign}\label{assAA}
\qquad\begin{cases}
 \,|a(x,  z)|+|\partial_z a(x,  z)| [H_{s}(z)]^{1/2}\\
 \quad  \le L[H_{s}(z)]^{(q-1)/2}+L[H_{s}(z)]^{(p-1)/2}\\
\, \nu [H_{s}(z)]^{(p-2)/2}|\xi|^{2}\le \partial_{z}a(x,z)\xi\cdot \xi \\
\, |a(x_{1},z)- a(x_{2}, z)| \\
\quad   \le L |x_{1}-x_{2}|^{\alpha }([H_{s}(z)]^{(q-1)/2}+[H_{s}(z)]^{(p-1)/2})
\end{cases}
\end{flalign}
whenever $x_1, x_2 \in \Omega$ and $z, \xi\in \er^n$, $|z|\not=0$. A natural ambiguity on the notion of weak solution now arises. Is this assumed to belong to $W^{1,p}$, as in the case of minima of functionals with $(p,q)$-growth? Or it is rather to be taken from $W^{1,q}$? In this last case, in the distributional form of \rif{risolvi0}, that is
$$
\int_{\Omega} a(x,Du)\cdot D\varphi \dx =0, \qquad \forall\ \varphi \in C^{\infty}_0(\Omega)
$$
we can take test functions $\varphi$ that are proportional to the solution $u$ and again derive basic energy estimates. Therefore considering $W^{1,q}$-solutions corresponds to take energy solutions in this setting. Issues concerning the various notions of weak solutions to problems with $(p,q)$-growth conditions are examined in depth in Zhikov's papers \cite{Z1, Z2, Z3, Z4}. In view of this discussion, there are now two possible approaches usually pursued in the literature
\begin{itemize}
\item Prove a priori estimates for more regular, i.e. $W^{1,q}$-solutions \cite{LU, M2, simon2}, thereby starting by energy solutions. 
\item Simultaneously proving both the existence and regularity for solutions to assigned boundary value problems, say for instance Dirichlet problems \cite{seminale, ivanov2, M2}. 
\end{itemize}
No distinction between the two approaches takes place when $p=q$. 
The second point from the above couple is the one we are now proposing, as already done in \cite{piovra} (and in most of the literature).  For this we consider the Dirichlet problem
\eqn{dir1}
$$
\left\{
\begin{array}{c}
-\divo\, a (x,Du)=0\quad \mbox{in}\ \Omega\\ [3 pt]
u \equiv u_0  \quad \mbox{on}\ \partial \Omega\;,
 \end{array}\right.
\qquad\ \ 
u_0 \in W^{1,\frac{p(q-1)}{p-1}}(\Omega)\;,
$$
where $\Omega\subset \er^n$ is a bounded and Lipschitz domain. 
\vspace{3mm}\begin{mdframed}
\begin{theorem}\label{t6}
Assume that the vector field $a(\cdot)$ satisfies \eqref{assAA}. 
If 
\eqn{deter}
$$
\frac{q}{p}\leq  1+  \frac{p-1}{10 p}\left(\frac{\alpha}{n}\right)^2\,,
$$
then there exists a solution $u\in W^{1,p}(\Omega)$ to the Dirichlet problem \eqref{dir1}, such that $Du$ is locally H\"older continuous in $\Omega$.
 Moreover, the estimate
$$
\|Du\|_{L^{\infty}(\Omega_0)}\lesssim_{\data, \alpha} \frac{1}{[\dist(\Omega_0, \partial \Omega)]^{\kk}} \left(\int_{\Omega} (|Du_0|+1)^{\frac{p(q-1)}{p-1}}\dx +1  \right)^{ \kk}$$
holds whenever $\Omega_0\Subset \Omega $ is an open subset, with $\kk \equiv \kk(n,p,q,\alpha)\geq 1$.
\end{theorem} 
\end{mdframed}\vspace{3mm}
\vspace{3mm}\begin{mdframed}
\begin{theorem}\label{c4} 
In the setting of Theorem \ref{t6} with $\alpha \in (0,1)$, assume also that $p\geq 2$, $s>0$, and that $\partial_z A(\cdot)$ is continuous on $\Omega \times \er^n$. Then $
u\in C^{1, \alpha}_{\loc}(\Omega).
$\end{theorem}
\end{mdframed}\vspace{3mm}
In the case $p=q=2$, Theorem \ref{c4} is a classical result from the uniformly elliptic theory \cite{cacc1, hopf, js, gg2,gg3}. The main point in the proof of Theorems \ref{t2}-\ref{c4} is that plain perturbation arguments, of the type used in the nonuniformly elliptic setting, do not work due to the aforementioned lack of homogeneous estimates (estimates that are not invariant under scaling). To overcome this point, in \cite{piovra} a new approach via renormalized Caccioppoli inequalities in fractional Sobolev spaces is considered. After this, a fractional version of De Giorgi's geometric iteration technique, and the use of nonlinear potentials via Lemma \ref{revlem}, leads to establish gradient boundedness. Once this is achieved, one can then then adapt more standard perturbation techniques. We shall expand on this in Section \ref{fracschau} below. 
\begin{remark} On the contrary of \rif{bound333}, the bound in \rif{deter} deteriorates as $p\to1$, forcing $q/p\to 1$. In both cases the factors $1/5$ and $1/10$ have been introduced in order to simplify the presentation and can be replaced by slightly larger quantities. For instance, \rif{bound333} can be replaced by 
$$
\frac{q}{p}<  1+  \texttt{k}(n,p, q,\alpha)\left(\frac{\alpha}{n}\right)^2\,, \qquad \frac 15 < \texttt{k}(n,p, q, \alpha) <1
$$
where $\texttt{k}(n,p,q, \alpha)$ is an involved function of its parameters. See \cite[Proposition 7.1]{piovra}. 
\end{remark}
\subsection{Relaxation}\label{lavgap} 
For general functionals of the type
\eqn{FMx}
$$
w \mapsto \mathcal{F}_{\texttt{x}}(w,\Omega) := \int_{\Omega}  F(x,Dw)\, dx
$$
Theorem \ref{t2} continues to hold provided 
\eqn{xx.0}
$$
\begin{cases}
\, z \mapsto F(x, z) \ \mbox{satisfies \rif{assFF} uniformly with respect to $x \in \Omega$}\\
\, |\partial_{z}F(x_{1},z)-\partial_{z}F(x_{2},z)|\\
\, \  \le L|x_{1}-x_{2}|^{\alpha}([H_{\mu}(z)]^{(q-1)/2}+[H_{\mu}(z)]^{(p-1)/2})
\end{cases}
$$
and the approximation-in-energy property \rif{energyconv} holds for the minimizer in question $u$. This is a in a sense a necessary and natural assumption in view of the potential occurrence of the Lavrentiev phenomenon \rif{lavvi}. In fact, more can be said. The first step is to define the so-called relaxed functional \cite{BFM, FMa, FMal, M0, ma5, M1, sc}
\eqn{rilassato}
$$
\overline{\mathcal F_{\texttt{x}}}(w,U)  \, := \inf_{\{w_k\}\subset W^{1,q}(U)}  \left\{ \liminf_k  \mathcal F_{\texttt{x}}(w_k,U)  \, \colon \,  w_k \deb w\  \mbox{in} \ W^{1,p}(U)  \right\} 
$$
for every open subset $U \subset \Omega$ and $w \in W^{1,1}_{\loc}(\Omega)$. 
Considering this type of lower semicontinuous envelope goes back to Lebesgue, Caccioppoli, Serrin and De Giorgi. In the nonuniformly elliptic setting, it appears for the first time in the work of Marcellini \cite{M0, ma5}. Accordingly, the Lavrentiev gap functional is defined by
\eqn{lav1}
$$
\mathcal {L}_{\mathcal{F}_{\texttt{x}}}(w,U):=
\begin{cases}
\, \displaystyle  \overline{\mathcal F_{\texttt{x}}}(w,U)- \mathcal{F}_{\texttt{x}}(w,U) \quad \mbox{if $\mathcal{F}_{\texttt{x}}(w,U)< \infty$}\\
0 \quad \mbox{if $\mathcal{F}_{\texttt{x}}(w,U)= \infty$\,.}
\end{cases}
$$
for every $w\in W^{1,1}(U)$. We have
\begin{itemize}
\item By $W^{1,p}$-weak lower semicontinuity of $\mathcal{F}_{\texttt{x}}(\cdot,U)$ we have $\mathcal{F}_{\texttt{x}}(\cdot,U)\leq \overline{\mathcal{F}_{\texttt{x}}}(\cdot,U)$ and therefore $\mathcal {L}_{\mathcal{F}_{\texttt{x}}}(\cdot,U)\geq 0$. It follows that $w \in W^{1, p}(U)$ whenever $ \overline{\mathcal F_{\texttt{x}}}(w,U)$ is finite. 
\item A function $u \in W^{1, p}(\Omega)$ is a minimizer  of $\overline{\mathcal F_{\texttt{x}}}(\cdot, \Omega)$ iff $\overline{\mathcal F_{\texttt{x}}}(u, \Omega)$ is finite and $\overline{\mathcal F_{\texttt{x}}}(u, \Omega) \leq \overline{\mathcal F_{\texttt{x}}}(w, \Omega)$ holds whenever $w \in u+ W^{1,1}_0(\Omega)$. 
\item 
Assume that $u\in W^{1,1}_{\loc}(\Omega)$ is a minimizer of the original functional $\mathcal F_{\texttt{x}}$ in \rif{FMx} such that $\mathcal {L}_{\mathcal{F}_{\texttt{x}}}(u,B)\equiv 0$ for every ball $B\Subset \Omega$. Then
\eqn{minimal}
$$
 \overline{\mathcal F_{\texttt{x}}}(u,B) =\mathcal{F}_{\texttt{x}}(u,B)  \leq \mathcal{F}_{\texttt{x}}(w,B) \leq \overline{\mathcal F_{\texttt{x}}}(w,B)  
$$ holds whenever $w \in u+W^{1,1}_0(B)$. Therefore $u$ also minimize $\overline{\mathcal F_{\texttt{x}}}(\cdot, B)$ for every ball $B \Subset \Omega$. 
\end{itemize}
This last point allows to follow a natural strategy, devised in \cite{sharp}: first, proving regularity of minima of the relaxed functional, and, from this, deducing the regularity of minima when the Lavrentiev gap vanishes. This is in the next results, still taken from \cite{piovra}. 
\vspace{3mm}\begin{mdframed}
\begin{theorem}\label{t4}
Let $u\in W^{1,p}(\Omega)$ be a minimizer of the functional $\overline{\mathcal F_{\texttt{x}}}(\cdot, \Omega)$, where $\Omega$ is a Lipschitz regular domain, and under assumptions \trif{bound333} and \eqref{xx.0}. Then $Du$ is locally H\"older continuous in $\Omega$. Moreover
\eqn{stima3}
$$
\|Du\|_{L^{\infty}(\Omega_0)} \lesssim_{\data, \alpha}  \frac{1}{[\dist(\Omega_0, \partial \Omega)]^{\kk}} \left[ \overline{\mathcal F_{\texttt{x}}}(u, \Omega) +1\right]^{\kk}
$$ 
holds whenever $\Omega_0\Subset \Omega $ is an open subset, where $\kk \equiv \kk(n,p,q,\alpha)\geq 1$. 
\end{theorem}
\end{mdframed}
\newpage\begin{mdframed}
 \begin{corollary}\label{c3}
Let $u\in W^{1,1}_{\loc}(\Omega)$ be a minimizer of the functional $\mathcal F_{\texttt{x}}$ in \eqref{FMx}, under assumptions \trif{bound333} and \eqref{xx.0}. Assume that 
\eqn{annulla}
$$\mathcal {L}_{\mathcal{F}_{\texttt{x}}}(u,B)=0 \quad \mbox{holds for every ball $B \Subset \Omega$}\,.$$ 
Then $Du$ is locally H\"older continuous in $\Omega$. Moreover, \eqn{stima3m}
$$
\|Du\|_{L^{\infty}(B_{t})} \lesssim_{\data, \alpha} \frac{1}{(r-t)^{\kk}}\left[ \mathcal F_{\texttt{x}}(u, B_{r})+1\right]^{\kk}
$$ holds whenever $  B_t \Subset B_{r} \Subset \Omega$ are concentric balls  with $r\leq 1$, where $ \kk\equiv \kk(n,p,q,\alpha)\geq 1$. 
\end{corollary}
\end{mdframed}\vspace{3mm}
To make the Theorem \ref{t4} and Corollary \ref{c3} really effective, it remains to understand when \rif{annulla} holds. This is in fact equivalent to require that the approximation-in-energy \rif{energyconv} takes place for $u$ on every ball $B \Subset \Omega$. This in fact happens in several common situations. For instance, via a simple convolution argument it is easy to see \cite{sharp} that if there exists a convex function $G\colon \er^n \to [0, \infty)$ such that 
$$
G(z) \lesssim F(x,z) \lesssim G(z) +1\,,$$ 
then $\mathcal {L}_{\mathcal{F}_{\texttt{x}}}(\cdot,B)\equiv 0$ holds for every ball $B \Subset \Omega$. This is the case of Theorem \ref{t2}, that in fact follows from Corollary \ref{c3}. The property is in a sense self-reproductive, as if $F(\cdot)$ is an integrand such that  \rif{energyconv} holds, then \rif{energyconv} also holds for any other Crarath\'eodory integrand $F_1(\cdot)$ for which
$$
F(x,z) \lesssim F_1(x,z) \lesssim F(x,z) +1$$ 
holds true. For instance \rif{energyconv} holds whenever 
$$
|z|^p+\aaa(x)|z|^q \lesssim F(x,z) \lesssim |z|^p+\aaa(x)|z|^q +1$$ 
is satisfied, as described before \rif{consider}. See Section \ref{softy} and \cite{sharp} for further cases in which \rif{energyconv} holds. 
\subsection{Non-differentiable functionals}\label{noeln}
A classical series of results in the Calculus of Variations are concerned with functionals of the type
\begin{flalign}\label{ggg2}
w\mapsto 
\int_{\Omega}[F(Dw)+ \hhh(x, w)] \dx
\end{flalign}
treated for instance in \cite{Ce1, gg2, hast, KM, phil, stamp}. The situation is as follows. When the Carath\'eodory function $\hhh(x,y)$ is differentiable with respect to the second variable, regularity of minima can be proved using the Euler-Lagrange equation, that is
 \eqn{eln}
$$
\divo\, \partial_{z}F(Du)= \partial_{y}\hhh(x, u)\,.
$$
Things change when $\hhh(\cdot)$ is not differentiable, let's say it only satisfies
\eqn{asshh}
$$
 |\hhh(x,y_{1})-\hhh(x,y_{2})|\le L|y_{1}-y_{2}|^{\alpha}, \quad \alpha\in (0,1]
$$
and equation \rif{eln} cannot be used simply because the right-hand side does not exist (unless $\alpha=1$). In the uniformly elliptic case a breakthrough was achieved in \cite{gg2} where it was proved that, if $\partial_zF(\cdot)$ satisfies \rif{ass1} with $p=2$ and $\hhh(\cdot)$ is a bounded and H\"older continuous functions, then $Du$ is locally H\"older continuous in $\Omega$. Later on, the literature reports several extensions of this kind of result in the case $p\not =2$, starting by \cite{manth1, manth2}. All such methods fail when $F(\cdot)$ is not uniformly elliptic, for the same reasons explained in Section \ref{piovra1} (lack of homogeneous estimates). Here we report the first gradient regularity result for minima of non-differentiable functionals of the type in \rif{ggg2}, obtained in \cite{piovra}. 
\vspace{3mm}\begin{mdframed}
\begin{theorem}\label{t1}
Let $u\in W^{1,1}_{\loc}(\Omega)$ be a minimizer of the functional in \trif{ggg2}, under assumptions \eqref{assFF} and \eqref{asshh} and 
\eqn{boundlolo}
$$
\frac{q}{p}\leq  1+\frac 15\left(1-\frac{\alpha}{p}\right)\frac{\alpha}{n}\,.
$$
Then $Du$ is locally H\"older continuous in $\Omega$. Moreover, 
$$
\|Du\|_{L^{\infty}(B_{t})}\lesssim_{\data, \alpha} \frac{1}{(r-t)^{\kk}}\left[ \|F(Du)\|_{L^1(B_{r})} +\|\hhh(\cdot,u)\|_{L^1(B_{r})}+1\right]^{\kk}
$$
holds whenever $B_{t}\Subset B_{r} \Subset \Omega$ are concentric balls with $r\leq 1$, where $ \kk\equiv \kk(n, p, q, \alpha)\geq 1$.
\end{theorem}
\end{mdframed}\vspace{3mm} 
\begin{boxy}{The key takeaway}
For uniformly elliptic problems Schauder estimates follow via perturbation methods. This is not possible in the nonuniformly elliptic case, due to basic lack of homogeneous estimates. They actually fail when not assuming a bound of the type $q/p\leq 1+ \texttt{o}(\alpha/n)$, encoding a direct and delicate interaction between growth of the ellipticity ratio and H\"older continuity of coefficients (see Section \ref{softy2}). On the positive side, Schauder estimates can be proved in the nonuniformly elliptic case by means of fractional spaces based techniques (see Section \ref{fracschau}). 
\end{boxy}

\newpage
\section{The sharp growth rate in nonuniformly elliptic\\ Schauder theory \cite{gioiello}}\label{gioi}
When looking at Theorem \ref{t2}, and in particular at the gap bound \rif{bound333}, it is natural to wonder whether this can be improved in at least something of the form 
\eqn{prebound}
$$
\frac{q}{p} < 1 + {\rm o}(n,\alpha)\,, \qquad  {\rm o}(\alpha, n) \approx \frac{\alpha}{n}
$$
that would be in line with \rif{boundma} and \rif{boundbs}. In other words, in something that approaches $1$ linearly, rather that quadratically, with respect to $\alpha/n$. 
In fact, we have 
\vspace{3mm}\begin{mdframed}
\begin{theorem}\label{maingioi}
Let $u\in W^{1,p}(\Omega)$ be a minimizer of the functional $\overline{\mathcal F_{\texttt{x}}}(\cdot, B)$ in \eqref{rilassato} for every ball $B\Subset \Omega$ under assumptions \eqref{xx.0}. If 
\eqn{bound333d}
$$
\frac{q}{p}<  1+  \frac{\alpha}{n}\,,
$$
then with $B_{\rr}\Subset \Omega$ being a ball, $\rr \in (0,1]$, and  $\theta\in (0,1)$, 
\begin{itemize} 
\item The local Lipschitz estimate
$$
\|Du\|_{L^{\infty}(B_{\rr/2})} \lesssim_{\data, \alpha}  \rr^{-n\kappa}\overline{\mathcal F_{\texttt{x}}}(u,B_{\rr})^{\kappa}+1
$$
holds with $\kappa\equiv \kappa(n,p,q,\alpha)\geq 1$. 
\item The local H\"older estimate
\eqn{m.2}
$$
[Du]_{0,\alpha_*;B_{\theta \rr}}\le c
$$
holds with both $\alpha_*\in (0,1)$ and $c\geq 1$ depending on $n,p,q,\alpha,L,\theta$ and $\overline{\mathcal F_{\texttt{x}}}(u,B_{\rr})$.
\item  If, in addition, $s >0$ and $\partial_{zz}F(\cdot)$ is continuous, then $u\in C^{1,\alpha}_{\loc}(\Omega)$ when $\alpha <1$. Specifically, 
$$
[Du]_{0,\alpha;B_{\theta \rr}}\le c_{*}
$$
holds as in \eqref{m.2}, where $c_{*}$ depends also on the modulus of continuity of $\partial_{zz}F(\cdot)$ and on $s$.  
\end{itemize}\end{theorem}
\end{mdframed}\vspace{3mm}
Needless to say, in the case the Lavrentiev gap functional vanishes, we recover the result for the original functional.
\vspace{3mm}\begin{mdframed}
\begin{corollary}\label{mfmf} 
Theorem \ref{maingioi} extends to any minimizer $u$ of the original functional in \eqref{FMx} provided $\mathcal {L}_{\mathcal{F}_{\texttt{x}}}(u,B)=0$ holds for every ball $B\Subset \Omega$. \end{corollary}
\end{mdframed}\vspace{3mm}
Theorem \ref{maingioi} is sharp, in particular with respect to \rif{bound333d}, in view of the counterexamples in \cite{sharp, FMM}. Comments on the proof of Theorem \ref{maingioi} can be found in Section \ref{maingioi2}. 
\section{Universal Lorentz conditions for Lipschitz regularity}
\label{steinuni}

\subsection{Stein type theorems} An essentially equivalent formulation of a classical theorem of Stein \cite{St} asserts\begin{boxy}{Stein's Theorem}
\eqn{solutio}
$$
\triangle u\in L(n,1) \Longrightarrow  
\mbox{$Du$ is continuous}\,.
$$
\end{boxy} 
\noindent This is sharp \cite{cianchijga}. Surprisingly enough, this condition is in a sense universal, and fo instance holds for energy solutions of the $p$-Laplacean system 
\begin{boxy}{Nonlinear Stein Theorem \cite{KMstein}}
\eqn{solutiop}
$$
\triangle_p u\in L(n,1) \Longrightarrow  
\mbox{$Du$ is continuous}\,.
$$
\end{boxy}
More in general, once the proper functional setting and the definition of solution is settled, the result continues to hold - at $C^1$ or $C^{0,1}$ regularity level - when replacing equation \rif{solutio} with \rif{risolvi}
in the case of
\begin{itemize}
\item A $p$-Laplacean type operator as in \rif{ass1} with Dini-continuous coefficients $x \mapsto a(x, \cdot)$ \cite{DMpoinc, KMli, KMri, KMguide}. 
\item The $p$-Laplacean system (again with Dini coefficients) \cite{DMpoinc, KMstein, KMvec} and degenerate systems involving differential forms \cite{sil}. See also \cite{thefive, thefive2}. 
\item General uniformly elliptic operators of the type in \rif{uni-ell0}-\rif{uni-ell}, \cite{Baroni}. 
\item  General nonuniformly elliptic operators, also with fast exponential growth, especially in the variational case \cite{seminale, ciccio, BS2, depicci}. See Sections \ref{expsec} and \ref{nonuniuni}. 
\item Fully nonlinear operators \cite{toti, pime} and renormalized $p$-Laplacean operators (via viscosity solutions) \cite{bamu}. 
\item General nonlinear elliptic systems, without Uhlenbeck type structure, and in the context of partial regularity \cite{KMpol, byun4}. This continues to hold in the case of minimizers of quasiconvex functionals \cite{camel, destro}. 
\end{itemize}
Moreover 
\begin{itemize}
\item Global Lipschitz regularity results for solutions to Dirichlet and Neumann problems, under suitable regularity on the boundary, can be obtained \cite{CMPDE, CMarma, CMjems}. 
\item Many of the results above extend to parabolic equations; see \cite{KMpisa, KMwolff, KMri}. 
\end{itemize}

The gradient continuity of solutions is essentially linked to \rif{prima0}. In fact, it holds that 
$$
\mu\in L(n,1) \Longrightarrow \lim_{r\to 0} {\bf I}_{1}^\mu(x,r)=0 \ \mbox{uniformly w.r.t. $x$}
$$
and therefore the continuity of $Du$ follows by the criterion displayed in \rif{seconda}. This is the general scheme of the various nonlinear versions of the original Stein's theorem, from degenerate to fully nonlinear elliptic problems. We refer to \cite{KMguide} for an overview. 

\subsection{Non-differentiable Stein \cite{piovra}} Of course solutions to \rif{solutio} locally minimize the functional in \rif{ggg2} when $F(Dw)\equiv |Dw|^2/2$ and $\hhh(x,w)=-\mu(x)w$. 
It is therefore natural to wonder what happens in more general cases. For this we consider a general integrand $F(\cdot)$ and 
\eqn{asshh2}
$$
\begin{cases}
\,  |\hhh(x,y_{1})-\hhh(x,y_{2})|\le \mu(x)|y_{1}-y_{2}|^{\alpha}, \quad \alpha\in (0,1]\\
\, \displaystyle \mu \in L\left(\frac n\alpha , \frac{1}{2-\alpha}\right)\,.
\end{cases}
$$
\vspace{3mm}\begin{mdframed}
\begin{theorem}\label{t5}
Let $u\in W^{1,1}_{\loc}(\Omega)$ be a minimizer of the functional in \trif{ggg2}, under assumptions \eqref{assFF} with $p=2$, \trif{boundlolo} and \eqref{asshh2}. Then $Du\in L^{\infty}_{\loc}(\Omega,\er^n)$ and moreover 
\begin{flalign*}
& \|Du\|_{L^{\infty}(B_{t})} \\
& \quad \le \frac{c}{(r-t)^{\kk}}\left[\|F(Du)\|_{L^1(B_{r})}+\|\hhh(\cdot,u)\|_{L^1(B_{r})}+\|\mu\|_{n/\alpha,1/(2-\alpha);B_{r}}+1\right]^{\kk}
\end{flalign*}
holds whenever $B_{t}\Subset B_{r} \Subset \Omega$ are concentric balls with $r\leq 1$, where $c\equiv c(\data, \alpha)\geq 1$ and $ \kk\equiv \kk(n, p,q,\alpha)\geq 1$.
\end{theorem}
\end{mdframed}\vspace{3mm}
Notice that, when $\alpha=1$, the Lorentz condition \rif{asshh2}$_2$ gives, as expected, $\mu \in L(n,1)$. Moreover, following the monotonicity properties in \rif{lorentzbasic}$_{1,2}$, condition \rif{asshh2}$_2$ gets weaker when $\alpha$ increases, again as expected. A more general result can be found in \cite{piovra}. 

\subsection{Stein type theorems in the fast growth case \cite{seminale, ciccio}}\label{expsec}
There are integrands whose nonuniform ellipticity cannot be controlled via polynomial growth of the eigenvalues. These can be of the type 
\eqn{veryfast} 
$$
w \mapsto \int_{\Omega} \exp(\exp(\dots \exp(|Dw|^p)\ldots))\, dx \;,  \qquad p \geq 1\,,$$
treated for instance in \cite{DE, evans, liebe2, M4, M5}.
Following the computations in \cite[(6.13)]{seminale}, in the case of \rif{veryfast} we find the following optimal upper bound on the ellipticity ratio
\eqn{veryfast2} 
$$
\mathcal R (z)\lesssim t^{p-1} \exp(\exp(\dots \exp(|z|^p)\ldots))+1
$$
where, if $k \geq 1$ is the number of nested exponentials involved in \rif{veryfast}, the number in \rif{veryfast2} is $k-1$ (it is zero when \rif{veryfast} involves one $\exp$ only and in this case $\mathcal R(z)$ grows polynomially in $|z|$). We report a sample result, concerning a full nonautonomous version of \rif{veryfast}. A main point here is the dependence on coefficients, that are assumed to be Sobolev functions, as first considered in the $(p,q)$-case in \cite{elema}. For this, we fix functions $\{p_k(\cdot)\}$ and $\{\ccc_k(\cdot)\}$, all defined on the open subset $\Omega \subset \er^n$, such that
\eqn{assexp}
$$
\begin{cases}
1 < p_{\mathfrak{m}} \leq p_0(\cdot) \leq p_M \,, \quad  0 < p_m \leq p_k(\cdot)\leq p_M\,, \quad \mbox{for $k \geq 1$} \\
 0 < \nu \leq \ccc_k(\cdot)\leq L \,, \quad p_k(\cdot), \ccc_k(\cdot) \in W^{1,d}(\Omega)\,, \ \  d>n\,, \quad \mbox{for $k \geq 0$}\;.
\end{cases}
$$
For every $k\in \en$, we next inductively define ${\bf e}_{k}\colon \Omega \times  [ 0, \infty) \to \er$ by
 $$\left\{
 \begin{array}{ccc}
 {\bf e}_{k+1}(x,t) & := & \exp \left(\ccc_{k+1}(x)\left[{\bf e}_{k}(t)\right]^{p_{k+1}(x)}\right) \\ [6 pt]{\bf e}_{0}(x,t) &:= & \exp\left(\ccc_0(x)t^{p_0(x)}\right)\;,
 \end{array}
 \right.
 $$
and consider functionals
\eqn{espik}
$$
w \mapsto \mathcal E_k(w, \Omega):=   \int_{\Omega}  \left[ {\bf e}_{k}(x,|Dw|)-w \mu \right]\dx\;.
$$
For this we have the following result, taken from \cite{ciccio}, and valid in the vectorial case too $N\geq 1$:
\vspace{3mm}\begin{mdframed}
\begin{theorem}\label{sample5} Let $u \in W^{1,1}_{\loc}(\Omega, \er^N)$ be a minimizer of the functional $\mathcal E_k$ in 
\trif{espik} for some $k \in \en$, under assumptions \trif{assexp} and such that $\mu\in L(n,1)$ with $n\geq 3$. Then $Du\in L^{\infty}_{\loc}(\Omega,\er^{N\times n})$. 
\end{theorem}
\end{mdframed}\vspace{3mm}
For more general cases and the two-dimensional situation $n=2$, we refer to \cite{seminale, ciccio}. Amongst the various applications of the methods in \cite{ciccio} we mention the possibility of proving sharp conditions for Lipschitz regularity of solutions to obstacle problems involving nonuniformly elliptic functionals. For this we consider a measurable function $\psi \colon \Omega \to \er$  
and the convex set $$\mathcal{K}_{\psi}(\Omega):=\{w\in W^{1,1}_{\loc}(\Omega)\colon  w(x)\ge \psi(x) \  \mbox{for a.e.} \ x\in \Omega \}\,.$$ We then say that a function $u \in W^{1,1}_{\loc}(\Omega)\cap \mathcal{K}_{\psi}(\Omega)$ is a constrained local minimizer of $\mathcal E_k$ if, for every open subset $\tilde \Omega\Subset \Omega$, we have $\mathcal E_k(u;\tilde \Omega) <\infty$ and if $\mathcal E_k(u;\tilde \Omega)\leq \mathcal E_k(w;\tilde \Omega)$ holds for every competitor $w \in u + W^{1,1}_0(\tilde \Omega)$ such that $w \in \mathcal{K}_{\psi}(\tilde \Omega)$. 
\vspace{3mm}\begin{mdframed}
\begin{theorem}\label{sample9}
Let  $u\in W^{1,1}_{\loc}(\Omega)\cap \mathcal{K}_{\psi}(\Omega)$ be a constrained local minimizer of $\mathcal E_k$ in \trif{espik} for some $k \in \en$, with $\mu\equiv 0$ and $n\geq 3$. If $\psi\in W^{2,1}_{\loc}(\Omega)$ and $|D^{2}\psi|\in L(n,1)$, then $Du\in L^{\infty}_{\loc}(\Omega,\er^n)$. \end{theorem}
\end{mdframed}\vspace{3mm}
Theorem \ref{sample9} continues to hold when, instead of functionals as in \rif{espik}, we consider general functionals with $(p,q)$-growth as for instance the one in Theorem \ref{marcth}. Details are again in \cite{ciccio}. 

Commenting on the growth of the ellipticity ratio \rif{veryfast2} is probably useful here. In the case of functionals with $(p,q)$-growth conditions a quantified rate of H\"older continuity of coefficients 
\eqn{sobby3}
$$x \mapsto \frac{\partial_{z} F(x, z)}{(|z|^2+s^2)^{(q-1)/2}}\,, \qquad |z|+s\not=0$$ is necessary to guarantee Lipschitz continuity of minima; see Section \ref{softy} and the gap bound in \rif{bound1}. The same holds even assuming Sobolev differentiability on coefficients, that is 
\eqn{sobby}
$$
\frac{\left|\partial_{zx} F(x, z)\right|}{(|z|^2+s^2)^{(q-1)/2}} \lesssim g(x) \in L^d(\Omega)\,,\qquad d>n\,.
$$ 
This is for instance the approach from \cite{ciccio,elema}, where \rif{sobby} is used with $d$ large enough to satisfy 
\eqn{sobby2}
$$
\frac{q}{p} < 1+ \frac 1n - \frac 1d \Longleftrightarrow d > \frac{np}{p-(q-p)n}\,.
$$
Note that Sobolev-Morrey embedding theorem implies that the vector field in \rif{sobby3} belongs to $C^{0,\alpha}$ with $\alpha := 1-n/d$, so that \rif{sobby2} and \rif{bound1} actually coincide. This means that a quantified rate of integrability $d$ is needed anyway. On the other hand, in Theorem \ref{sample5} no lower bound on $d$ is needed. This is in striking contrast with the fact that, while in the case of \rif{espik} the ellipticity ratio grows exponentially, it just grows polynomially fast in $(p,q)$-growth functionals. This apparently counterintuitive situation is explained as follows. A relevant quantity implicitly appearing in the estimates is in fact a renormalized ratio 
$$
\sup_{x\in B}\, \frac{\mathcal R_{F}(x,z)}{F(x,z)}
$$
for a fixed ball $B \subset \Omega$. While in the $(p,q)$ case the above quantity decays to zero (as $|z|\to \infty$) polynomially, in the exponential case the convergence to zero is exponentially fast (compare \rif{veryfast2}).  This explains the stronger requirement on the integrability of coefficients in the polynomial growth case and also reflects the fact that exponential type functionals force Lipschitz estimates more strongly due to their fast growth (although they are more delicate to treat from the technical view point). The case of functionals with exponential growth therefore reveals to be closer to that of uniformly elliptic functionals and operators, for which we have the following borderline result, once again from \cite{ciccio}, that deals with elliptic systems:
\vspace{3mm}\begin{mdframed}
\begin{theorem}\label{sample6} Let $u \in W^{1,1}_{\loc}(\Omega, \er^N)$ be a distributional solution to 
$$
-\divo\, (\ccc(x)\tilde a(|Du|)Du)=\mu  \,, \quad 0 < \nu \leq \ccc(\cdot) \leq L\;,
$$
such that $A(\cdot, |Du|)\in L^1(\Omega)$, where $A(\cdot)$ is defined in \trif{defiA}, 
and under assumptions \trif{uni-ell}. If $\mu, |D\ccc|\in   L(n,1)$ with $n\geq 3$, then $Du\in L^{\infty}_{\loc}(\Omega,\er^{N\times n})$. 
\end{theorem}
\end{mdframed}\vspace{3mm}
The borderline nature of the Theorem \ref{sample6} with respect to Theorem \ref{sample5} relies in that $L^d \subset L(n,1)\subset L^n$ holds for every $d >n$; see \rif{lorentzbasic}. 
\vspace{2mm}

\begin{boxy}{The key takeaway}
Condition $\mu\in L(n,1)$ in \rif{genF}-\rif{risolvi} is universal for Lipschitz and/or $C^1$-regularity of minima and solutions, as it works both in the uniformly and in the nonuniformly elliptic case. It does not depend on the specific operator/integrand considered but it is ultimately an embedding type effect linked to any form of ellipticity. New Lorentz conditions appear when considering non-differentiable functionals. 
\end{boxy}

\section{Schauder estimates at nearly linear growth \cite{loggy}}\label{logsec}
The functionals of the type \rif{genF} considered in the previous sections have superlinear polynomial growth in the sense that 
$$
|z|^p \lesssim F(x,z)\,, \qquad p>1\,. 
$$
Therefore integrands with so-called nearly linear growth are excluded. In fact, these are characterized by 
$$
\lim_{|z|\to \infty} \, \frac{F(x,z)}{|z|}=\infty\,,\qquad \lim_{|z|\to \infty} \, \frac{F(x,z)}{|z|^p}=0 \quad \ \forall \ p>1\,.
$$
This happens for integrals of the type 
\eqn{modellol}
$$
 w\mapsto \int_{\Omega}\ccc(x)|Dw|\log(1+|Dw|)\dx\,, \qquad 0 < \nu \leq \ccc(\cdot) \leq L
$$
that we already encountered in Section \ref{esempi}. Going to the realm of $(p,q)$-growth problems when $p=1$, a borderline case of the double phase functional in \rif{doppio} is
\eqn{modello}
$$
 w\mapsto \int_{\Omega}\left[\ccc(x)|Dw|\log(1+|Dw|)+\aaa(x)|Dw|^{q}\right] \dx\,.
$$
This one was first considered in \cite{dm} under Sobolev differentiability assumption on the coefficients $\ccc(\cdot), \aaa(\cdot)$. A main point here is that, although the double phase integral in \rif{doppio} is softly nonuniformly elliptic in the sense of \rif{limitatasoft}, the functional in \rif{modello} is nonuniformly elliptic in the classical sense of \rif{classe}. Indeed, the best upper bound we can find for the ellipticity ratio is
$$
\mathcal{R}(x, z) \lesssim \log(1+|z|) +1\,.
$$
For this very reason, the extension to the functional \rif{modello} of the techniques from \cite{BCM, CMa1, CMa2, HO1, HO2}, devised to treat softly nonuniformly elliptic integrals, is impossible. A different approach has been introduced in \cite{loggy} and gives
\vspace{3mm}\begin{mdframed}
\begin{theorem}\label{t1log}
Let $u\in W^{1,1}_{\loc}(\Omega)$ be a local minimizer of the functional in \trif{modello}, with
\eqn{bound}
\begin{cases}
$$
\displaystyle 0\le \aaa(\cdot)\in C^{0,\alpha}(\Omega), & 1<q<1+\alpha/n\\
 \ccc(\cdot)\in C^{0, \ao}_{\loc}(\Omega),&\displaystyle 1/\Lambda \leq  \ccc(\cdot) \leq \Lambda\,,
\end{cases}
$$ where $\alpha, \ao\in (0,1)$ and $\Lambda\geq 1$. 
Then $Du$ is locally H\"older continuous in $\Omega$ and moreover, for every ball $B_r\Subset \Omega$, $r\leq1$, there holds
$$
\|Du\|_{L^{\infty}(B_{r/2})}\le c\left(\mint_{B_r}\left[|Du|\log(1+|Du|)+\aaa(x)|Du|^{q}\right] \dx\right)^{\kk} +c
$$
with $$\begin{cases}
c\equiv c(n,q, \Lambda,\alpha, \alpha_0,\|a\|_{C^{0,\alpha}}, \|\ccc\|_{C^{0,\ao}})\geq 1\\
\kk\equiv \kk(n,q,\alpha, \alpha_0)\geq 1\,.
\end{cases}
$$
\end{theorem}
\end{mdframed}\vspace{3mm}
The main assumption in Theorem \ref{t1log} is
\eqn{ilbb}
$$
q < 1+ \frac \alpha n\,,
$$
which is the obvious borderline version of \rif{bound1} (let $p\to1$). The equality case in \rif{doppiob}, which is lacking in 
 \rif{ilbb}, is based on Gehring type results \cite{dm0} that typically miss in the nearly linear growth case. In the non-degenerate case, the gradient H\"older continuity exponent can be quantified and we have an analog of Theorems \ref{c1} and \ref{c4}, that is 
\vspace{3mm}\begin{mdframed}
\begin{theorem}\label{t2log}
Under assumptions \trif{bound}, local minimizers of the functional 
$$
 w\mapsto \int_{\Omega}\big[\ccc(x)|Dw|\log(1+|Dw|)+\aaa(x)(|Dw|^{2}+1)^{q/2}\big] \dx
$$
are locally $C^{1,\tia/2}$-regular in $\Omega$, where $\tia:=\min\{\alpha_0, \alpha\}$. In particular, local minimizers of the functional in \trif{modellol} are locally $ C^{1,\alpha_0/2}$-regular provided \trif{bound}$_2$ is assumed. 
\end{theorem}
\end{mdframed}\vspace{3mm}
One might of course wonder whether or not it is possible to approach linear growth conditions from below more closely, for instance considering arbitrary compositions of logarithms as follows:
\eqn{generale}
$$
 w\mapsto \int_{\Omega}\big[F(x,Dw)+\aaa(x)(|Dw|^{2}+s^2)^{q/2}\big] \dx\,,
$$
where $s \in [0,1]$ and 
\eqn{verylinear}
$$
\begin{cases}
\ F(x,z) \equiv \ccc(x)|z|L_{k+1}(|z|)  \ \  \mbox{for $k\geq 0$} \\
 \ L_{k+1}(|z|)=\log(1+L_{k}(|z|))  \ \  \mbox{for $k\geq 0$}\,, \quad 
L_{0}(|z|)= |z|\,.
\end{cases}
$$
Here $\ccc(\cdot)$ is as in \rif{bound}$_2$. For this we have to cook up a more technical setting. We consider continuous integrands $F\colon \Omega \times \er^n\to [0, \infty)$ such that $z \mapsto F(x, z) \in C^{2}(\er^n)$ for every $x\in \Omega$ and satisfying
\eqn{assif}
$$
\begin{cases}
\ \nu|z|g(z) \leq F(x,z) \leq L(|z|g(z) +1)\\
\ \displaystyle \frac{\nu |\xi|^2}{(|z|^2+1)^{\gamma/2}}  \leq \partial_{zz}F(x,z)\xi \cdot \xi \, , \quad  |\partial_{zz}F(x,z)| \leq  \frac{Lg(|z|)}{(|z|^2+1)^{1/2}} \\
\ |\partial_zF(x,z)-\partial_z F(y,z)| \leq L |x-y|^{\ao}g(|z|), 
\end{cases}
$$
for every choice of $x,y\in \Omega$, $z, \xi \in \er^n$, where $\ao \in (0,1)$ and $1\leq \gamma <3/2$ being fixed constants. Here $g\colon [0, \infty) \to [1, \infty)$ is a non-decreasing, concave and unbounded function such that $ t\mapsto tg(t) $ is convex. Moreover, we assume that for every $\eps >0$ there exists a constant $c_{g}(\eps)$ such that
 \eqn{0.1}
 $$
 g(t) \leq c_{g}(\eps) t^{\eps}\quad \mbox{holds for every $t\ge 1$}\,.
$$
Conditions \rif{assif}-\rif{0.1} imply the local Lipschitz continuity of minima. Local H\"older continuity needs one more condition, that is
\eqn{assi3}
$$
|\partial_z F(x, z)| \leq L|z| \quad \mbox{holds for every $|z|\leq 1$\,.}
$$
We then have
\vspace{3mm}\begin{mdframed}
\begin{theorem}\label{t3log}
Let $u\in W^{1,1}_{\loc}(\Omega)$ be a local minimizer of the functional in \trif{generale} under assumptions \trif{bound}$_1$ and \trif{assif}-\trif{0.1}. There exists $\gamma_m \in (1,2)$, depending only $n,q, \alpha_0, \alpha$, such that, if $1\leq \gamma < \gamma_m$, then 
$$
\|Du\|_{L^{\infty}(B_{r/2})}\le c\left(\mint_{B_r}\big[F(x, Du)+\aaa(x)(|Du|^2+s^2)^{q/2}\big] \dx\right)^{\kk} +c
$$
holds whenever $B_r\Subset \Omega$, $r\leq 1$. Here it is $c\equiv c(n,q, \nu, L, \Lambda,\alpha, \alpha_0, \linebreak\|a\|_{C^{0,\alpha}},\|\ccc\|_{C^{0,\ao}})\geq 1$ and $\kk \equiv \kk (n,q, \alpha_0, \alpha)\geq 1$. Moreover, assuming also \trif{assi3} implies that $Du$ is locally H\"older continuous is $\Omega$. Finally, if also $s >0$ and $\partial_{zz}F(\cdot)$ is continuous, then $u \in C^{1,\tilde \alpha/2}_{\loc}(\Omega)$-regular with $\tilde \alpha = \min\{\alpha_0, \alpha\}$.   \end{theorem}
\end{mdframed}\vspace{3mm}
As a matter of fact, Theorems \ref{t1log} and \ref{t2log} follow from Theorem \ref{t3log} with $\gamma=1$. In order to treat functionals in \rif{generale} with the choice in \rif{verylinear} it is sufficient to apply Theorem \ref{t3log} with any $\gamma>1$. For more details the reader is referred to \cite{loggy}. The proof of Theorem \ref{t3log} is again based on the use of fractional energy estimates. This time a careful use of the double phase structure of the integral \rif{generale} is necessary, and involves an ad hoc anisotropic version of the Bernstein technique. The final outcome is the Lipschitz regularity of minima under the optimal bound \rif{ilbb}. Once again, once local Lipschitz regularity is gained, one can (carfully) adapt more standard perturbation methods. 
 
The methods introduced in \cite{loggy} are general enough to provide a blueprint to treat more general situations. One of these is in  \cite{depiccif}, where the authors considered multi-phase models with nearly linear growth. In this case (three phases case) a model is provided by the functional 
\eqn{multiphase}
$$
 w\mapsto \int_{\Omega}\left[\ccc(x)|Dw|\log(1+|Dw|)+\aaa(x)|Dw|^{q}+\mf{b}(x)|Dw|^{\gamma}\right] \dx
$$
where
$$
\begin{cases}
\displaystyle 1 < q < 1+\frac{\alpha}{n}\,, \qquad 1 < \gamma < 1+\frac{\beta}{n}\\ \\
0\leq \aaa(\cdot) \in C^{0, \alpha}(\Omega)\,,\qquad 0\leq \mf{b}(\cdot) \in C^{0, \beta}(\Omega)
\end{cases}
$$
with $\alpha, \beta \in (0,1]$ and $\ccc(\cdot)$ as in \rif{bound}$_2$. In this situation, as it has been proved in \cite{depiccif},  minimizers of the functional in \rif{multiphase} have locally H\"older continuous first derivatives. Needless to say, the result extends to functionals of the type
$$
 w\mapsto \int_{\Omega}\left[\ccc(x)|Dw|\log(1+|Dw|)+\sum_{k=1}^m \mf{a}_k(x)|Dw|^{q_k}\right] \dx
$$
with 
\eqn{condille}
$$
\displaystyle 1 < q_k < 1+\frac{\alpha_k}{n}\,, \quad 0\leq \aaa_k(\cdot) \in C^{0, \alpha_k}(\Omega)\,, \quad \alpha_k \in (0,1]\,.
$$ 
The case of multi-phase functionals with super-linear growth, i.e., the one displayed in \rif{multino}, was originally considered in \cite{crismulti0}. Observe that conditions \rif{condille} are the sharp analogue of those in \rif{boundulti} when approaching nearly linear growth conditions, that is, when formally letting $p\to 1$ in \rif{multino} and \rif{boundulti}.

Functionals with nearly linear growth have been the object of intensive investigation over the years. For this we mention, for instance, \cite{bildhauer, BF, BF1, elemape, FM, FS2, MP}. Their nonuniform ellipticity stems from their closeness to linear growth functionals, like for instance the minimal surface one, that are always nonuniformly elliptic; for this see \cite{BeS, BeS2, gme1, gme2, gmek1, gmek2} and related references for recent developments. Beside their theoretical interest, nearly linear growth integrals also appeared in relation to applications in Plasticity and Non-Newtonian Fluid-dynamics \cite{FrS1, FrS2, FS2, FS1}. 

\begin{boxy}{The key takeaway}
Schauder type estimates still hold for minima of nearly linear growth integrals. These are in general nonuniformly elliptic in the classical, pointwise sense \rif{classe} and therefore perturbation methods are ruled out. Again, a direct proof of Lipschitz continuity is needed and goes via fractional Caccioppoli inequalities and nonlinear potential theoretic methods. 
\end{boxy}

\section{Renormalized Caccioppoli inequalities and nonlinear potentials}\label{rinormalizza}
In this section we shall try to present and synthesize some of the methods from \cite{seminale, ciccio, piovra, loggy, gioiello}, where the use of nonlinear potential theoretic methods in the setting of nonuniformly elliptic problems was introduced. As mentioned in Section \ref{sezione1}, in order to give a unified treatment and to emphasize the main ideas, we shall pick a few special cases, often not the most general ones. These can be anyway treated via more specific variants of the same leading principles, sometimes discussed in the commentaries. We shall demonstrate four different applications of the general ideas in order to emphasize the flexibility of the approach presented and provide a gradual, gentle introduction to the most technical cases. We shall start by a toy model in Section \ref{toy}, that actually deals with classical linear uniformly elliptic case, in order to highlight the connections between uniformly and nonuniformly elliptic problems. Of course we shall confine ourselves to give sketches of the proofs as a complete treatment would just (more than) double the length of these notes. 
\subsection{Linear elliptic equations}\label{toy} We consider a linear elliptic equation
\eqn{lineare2}
$$
-\divo\, (\texttt{A}(x)Dv)=\mu \,, \qquad     \nu\mathbb{I}_{\rm d} \leq  \texttt{A}(\cdot) \leq L\mathbb{I}_{\rm d}
$$
where $\texttt{A}(\cdot)$ has measurable entries. 
We want to prove that 
\eqn{imply1}
$$
\mu \in L\left(n/2, 1\right) \Longrightarrow v \in L^{\infty}_{\loc}(\Omega)
$$
provided $n> 4$. Of course we start assuming that 
$v\in W^{1,2}_{\loc}(\Omega)$, i.e., it is a standard energy solution (according to Definition \ref{defisol}). 
Notice that assuming $\mu \in L(n/2, 1)$ implies that $\mu \in L^2$ holds provided $n \geq 4$ by \rif{lorentzbasic}$_2$. 
\begin{boxy}{The Classical Caccioppoli inequality}
\eqn{para20}
$$
    \int_{B_{\rr/2}} |D(v-\kk)_+|^2 \, dx  \lesssim_{\nu, L}  \frac{1}{\rr^2} \int_{B_{\rr}} (v-\kk)_+^2\, dx
 +   \rr^2\int_{B_{\rr}} |\mu|^2\, dx
 $$
 \end{boxy}
 \noindent now holds for every ball $B_{\rr} \Subset \Omega$. As usual, here it is 
$ (v-\kk)_{+}:=\max\{v-\kk,0\}$ and we are taking any $\kk \geq 0$. The proof is very standard and can be obtained testing \rif{lineare2} by $\eta^2(v-\kk)_{+}$, where $\eta\in C^1_0(B_{\rr})$ is such that $\eta\equiv 1$ on $B_{\rr/2}$ and $|D\eta|\lesssim 1/\rr$ (see for instance \cite[Chapter 7]{G}). Sobolev embedding theorem now implies the following:
\begin{boxy}{Reverse type inequality with remainder}
\eqn{reversa}
    $$
   \left( \mint_{B_{\rr/2}} (v-\kk)_+^{2\chi} \, dx \right)^{1/\chi} \lesssim_{n,\nu, L} 
      \mint_{B_{\rr}} (v-\kk)_+^{2} \, dx +  \rr^4\mint_{B_{\rr}} |\mu|^2\, dx
    $$ 
    \end{boxy}
\noindent where $\chi := 2^*/2>1$ and $2^*$ has been defined in \rif{sobby0}. Applying Lemma \ref{revlem} with the choices $h=1$, $m_1=t=\delta_1=2$, $\theta_1=1$, $M_0=M_1\approx 1$, $\kk_0=0$, yields 
$$
 v(x_{0}) \lesssim_{n,\nu, L}  \left(\mint_{B_{r_0}(x_{0})} v^2  \dx\right)^{1/2}+\mathbf{P}^{2,1}_{2,2}(\mu;x_{0},2r_0)\,.
$$
This holds whenever $B_{2r_0}(x_{0}) \Subset \Omega$ such that $x_0$ is Lebesgue point of $v$. 
By observing that $-v$ solves \rif{lineare2} with $\mu$ replaced by $-\mu$, we arrive at the pointwise nonlinear potential estimate
\eqn{para1}
$$
 |v(x_{0})| \lesssim_{n,\nu, L}   \left(\mint_{B_{r_0}(x_{0})} v^2  \dx\right)^{1/2}+\mathbf{P}^{2,1}_{2,2}(\mu;x_{0},2r_0)
$$
that again holds whenever $B_{2r_0}(x_0)\Subset \Omega$, for a.e. $x_0$. We now consider a ball $B_{r}\Subset \Omega$; applying \rif{para1} to $B_{r_0}(x_0)$, $r_0=r/4$, for all $x_0\in B_{r/2}$ for which \rif{para1} holds, we get 
$$
\|v\|_{L^\infty(B_{r/2})} \lesssim_{n,\nu, L} \left(\mint_{B_{r}} v^2  \dx\right)^{1/2}+\|\mathbf{P}^{2,1}_{2,2}(\mu;\cdot,r/2)\|_{L^\infty(B_{r/2})}\,.
$$
It is now time to use Lemma \ref{crit}; this is possible since \rif{lo.2.1} is verified by our assumption $n>4$. We therefore conclude with
\eqn{para3}
$$
\|v\|_{L^\infty(B_{r/2})} \lesssim_{n,\nu, L} \left(\mint_{B_{r}} v^2  \dx\right)^{1/2}+\|\mu\|_{n/2,1;B_{r}}
$$
that holds for every ball $B_r\Subset \Omega$. 
\begin{commentary}
{\rm 
Both the result in \rif{imply1} and estimate \rif{para3} are sharp with respect to the function space $L(n/2,1)$ considered for the right-hand side. The result continues to hold in the cases $n=3,4$, that are missing here due to the use of the nonlinear potential $\mathbf{P}^{2,1}_{2,2}$. This is in fact a consequence of another, more basic nonlinear potential estimate, due to Kilpelainen \& Maly \cite{KiL} and valid for general quasilinear elliptic equations with measure data right-hand side, that in the case $p=2$ reduces to 
$$
|v(x_0)| \lesssim_{n,\nu, L}  {\bf I}_{2}^\mu(x_0,r) + \mint_{B_r(x_0)}v \dx \approx_{n}  
\mathbf{P}^{1,1}_{1,2}(\mu;x_{0},r) + \mint_{B_r(x_0)}v \dx \,.
$$
See also \cite{KMguide}. This again implies \rif{imply1} by Lemma \ref{crit} as this time \rif{lo.2.1} is satisfied for $n>2$ (see also \cite{cianchipisa}). In fact, more general results are available for equations of the type \rif{risolvi} under assumptions \rif{assumi}. 
In this case, for suitably defined solutions to \rif{risolvi}, and restricting for simplicity to the case $p\geq 2$, we have the
\begin{boxy}{Wolff potential estimate \cite{KiL, KMguide}}
\begin{flalign} |v(x_0)| &\lesssim_{\data}  {\bf W}_{1,p}^{\mu}(x_0,r) +  \left(\mint_{B_r(x_0)} |v|^{p-1} \, dx \right)^{1/(p-1)}\notag \\
&\approx_{n,p}    \mathbf{P}^{1,1}_{p-1,p/(p-1)}(\mu;x_0,r) +  \left(\mint_{B_r(x_0)} |v|^{p-1} \, dx \right)^{1/(p-1)}
\label{primadonna-mw}
\end{flalign} 
for a.e. $x_0 \in \Omega$ and balls $B_r(x_0)\Subset \Omega$. 
\end{boxy} 
In this case, again as a consequence of Lemma \ref{crit}, applied with $m=\theta=1$, $t=p-1$ and $\delta=p/(p-1)$, we find
\eqn{imply2}
$$
\mu \in L\left(\frac{n}{p}, \frac{1}{p-1}\right) \Longrightarrow v \in L^{\infty}_{\loc}(\Omega)
$$ 
provided 
\eqn{dime1}
$$\frac{n}{t\delta}=\frac{n}{p}>1\,,$$ 
with the a priori estimate 
\eqn{para5}
$$
\|v\|_{L^\infty(B_{r/2})}\lesssim_{\data} \left(\mint_{B_{r}} |v|^{p-1}  \dx\right)^{1/(p-1)}+\|\mu\|_{n/p,1/(p-1);B_{r}}^{1/(p-1)}
$$
that holds for every ball $B_{r}\Subset \Omega$. In fact, this last one reduces to \rif{para3} when $p=2$.  
Let us briefly see how to recover \rif{para5} without appealing to \rif{primadonna-mw}, but rather using more standard energy estimates as done for \rif{para3}. For equations of the type in \rif{risolvi}, under assumptions \rif{assumi}, the following Caccioppoli type inequality holds:
\eqn{para11}
$$
    \int_{B_{\rr/2}} |D(v-\kk)_+|^p \, dx  \lesssim_{\nu, L,p} \frac{1}{\rr^p} \int_{B_{\rr}} (v-\kk)_+^p\, dx
 +   \rr^{\frac{p}{p-1}}\int_{B_{\rr}} |\mu|^{\frac{p}{p-1}}\, dx
 $$ 
 for every ball $B_\rr \Subset \Omega$; see \cite[Chapter 7]{G}. 
 Then, again via Sobolev embedding, we come to the analog of \rif{reversa}, that is
 \eqn{reversap}
     $$ 
   \left( \mint_{B_{\rr/2}} (v-\kk)_+^{p\chi} \, dx \right)^{1/\chi} \lesssim_{\data} 
      \mint_{B_{\rr}} (v-\kk)_+^{p} \, dx +  \rr^{\frac{p^2}{p-1}}\mint_{B_{\rr}} |\mu|^{\frac{p}{p-1}}\, dx
    $$
    for some $\chi\equiv \chi(n,p)>1$. We can now proceed as after \rif{reversa}, this time the relevant nonlinear potential being  
    $\mathbf{P}^{p/(p-1),1}_{p,p/(p-1)}$ (Lemma \ref{revlem}). This yields \rif{imply2} upon checking that (Lemma \ref{crit})
    \eqn{dime2}
    $$
\frac{n }{t\delta}>1 \Longleftrightarrow n > \frac{p^2}{p-1} \,.
$$
Therefore we again come to \rif{para5}, and again with a restricted range of dimensions $n$. The difference between \rif{dime1} and \rif{dime2} stems from the kind of {\em underlying energy estimates} they rely on. In the case of \rif{para11} one is using estimates at the natural energy level $W^{1,p}$. These are natural estimates for energy solutions, but yield weaker results. In order to get estimates {\em below the natural growth exponent} like \rif{primadonna-mw}, more sophisticated techniques are needed \cite{KiL, KMguide}. They relate to so-called measure data problems \cite{elencoo, BG1, BG2, dalmy, LSW} and usually work at the energy level fixed by the space $W^{1, p-1}$.  
}
\end{commentary}
\begin{boxy}{The key takeaway}
Energy inequalities of the type in \rif{para20} and \rif{para11}, with related reverse H\"older type inequalities \rif{reversa} and \rif{reversap}, automatically encode $L^\infty$-bounds on solutions under optimal regularity assumptions on $\mu$,  via nonlinear potentials. Similar De Giorgi type iteration schemes work below the natural growth exponent allowing for slightly larger range of parameters \cite{KiL, KMguide}. The crucial point here is that, as the proof in \cite{seminale, piovra} of Lemma \ref{revlem} reveals, the basic estimates involved \rif{reversa} and \rif{reversap} must be homogeneous with respect to the function $v$. This happens here as a consequence of the assumed uniform ellipticity in \rif{lineare2} and \rif{assumi}. 
 \end{boxy}
\subsection{Nonuniformly elliptic problems \cite{seminale, ciccio}}\label{nonuniuni}
Here we switch to gradient estimates and present a few basic methods from \cite{seminale, ciccio}. 
For the sake of simplicity we shall confine ourselves to the case $n\geq 3$, again referring to \cite{seminale, ciccio} for the two dimensional case $n=2$. In this section we consider the simplified case of functionals of the type
\eqn{semplificato}
$$
w \mapsto   \int_{\Omega}  \left[ F(Dw)-\mu w \right]\, dx\,,
$$ 
where $F(\cdot)$ is a $C^2$-regular integrand satisfying \rif{assbase} as in Section \ref{gapbounds}. Moreover, we assume that $\mu \in L(n,1)$.  
Here a disclaimer is necessary: \cite{seminale, ciccio} are long papers and reporting here a self-contained proof of the results would be impossible. In particular, a lengthy and nonstandard procedure is necessary in \cite{seminale, ciccio} in order to approximate the original nonuniformly elliptic problem with uniformly elliptic ones and this is aimed at making all the computations necessary for the a priori estimates legal. The complexity of such a procedure stems from the fact that this is tailored to cover elliptic functionals whose growth of the eigenvalues is faster than polynomial, therefore going much beyond the $(p,q)$ case considered here. We shall therefore follow the classical approach of delivering the core a priori estimates for more regular solutions and data (see \rif{appiass}). Of course the shape of the estimates will not quantitatively incorporate such a priori assumption that can be eventually removed by the aforementioned approximation procedure. We now fix a local minimizer $u$ of the functional in \rif{semplificato} and assume that 
\eqn{appiass}
$$u\in W^{1,\infty}_{\loc}(\Omega)\cap 
W^{2,2}_{\loc}(\Omega)\,, \qquad \mu \in L^\infty_{\loc}(\Omega)\,.$$
In particular, this will guarantee that all the forthcoming integrals will be finite. 
We now proceed via a variant of the so-called Bernstein technique. In the uniformly elliptic case 
this classical method prescribes that, if $u$ solves an equation of the type $\divo\, a(Du)=0$, then, under appropriate  differentiability assumptions, functions of the type $|Du|^{\gamma}$ are subsolutions of certain linear elliptic equations, for suitable powers $\gamma>0$. As such, they are locally bounded functions and therefore original solutions $u$ are locally Lipschitz regular. In turn, upper bounds for subsolutions typically goes via Caccioppoli type inequalities on level sets and De Giorgi type iterations. In our setting we are using the Euler-Lagrange equation 
\eqn{elf}
$$-\divo \, \partial_zF(Du)=\mu$$ 
of the functional in \rif{semplificato}, and this poses two problems. First, it is nonuniformly elliptic. Second, the right-hand side does not vanish. Therefore the subsolution type scheme described above falls short of deriving estimates. Nevertheless, by differentiating \rif{elf} we can use a raw version of Bernstein's method. We derive a Caccioppoli type inequality for a certain convex function of $|Du|$, as at this stage we do not have any information concerning any type of subsolution property due to the presence of a general right-hand side $\mu$. From this Caccioppoli inequality we directly derive local gradient bounds with the aid of nonlinear potentials. Specifically, we start defining 
\eqn{triala}
$$
\begin{cases}
\displaystyle G(t):= \int_0^t (y^2+s^2)^{(p-2)/2}y\dy=\frac{1}{p}\left[(t^{2}+s^{2})^{p/2}-s^{p}\right]\\
 v := G(|Du|)\,.
 \end{cases}
$$ 
We next consider concentric balls 
\eqn{palline}
$$
B_{\rr}\equiv B_{\rr}(x_0) \subset  B_{r_0}(x_0)\subset  B_{2r_0}(x_0)\Subset \Omega\,.$$ Following \cite[Lemma 4.5]{seminale}, the Caccioppoli type inequality
\begin{flalign}
    \int_{B_{\rr/2}} |D(v-\kk)_+|^2 \, dx  & \lesssim_{\data}  \frac{1}{\rr^2} \int_{B_{\rr}} [\mathcal R_{F}(Du)+1](v-\kk)_+^2\, dx \notag  \\
    &\qquad  \ \ \ \ +   \int_{B_{\rr}}(|Du|^2+1) |\mu|^2\, dx \label{para202}
\end{flalign} 
 holds for every ball $B_{\rr}\Subset \Omega$ and number $\kk \geq 0$, where $\mathcal R_{F}(\cdot)$ denotes the ellipticity ratio of the integrand $F(\cdot)$. Note that this is a point where \rif{appiass} is needed. In fact \rif{para202} is obtained differentiating the Euler-Lagrange equation \rif{elf} so that  \rif{appiass} guarantees that this is possible and that all the terms generated by such a procedure make sense and are finite. By \rif{assbase} the best upper bound we can get on $\mathcal R_{F}(Du)$ is
 \eqn{para2022}
 $$
 \mathcal R_{F}(Du)\lesssim_{L/\nu} |Du|^{q-p}+1\,.
 $$  
In comparison to \rif{para20}, inequality \rif{para202} fails to be homogeneous with respect to $v$ due to the appearance of $ \mathcal R_{F}(Du)$. This makes proceeding as after \rif{para20} impossible. Using \rif{para2022} we renomalize \rif{para202} getting the following  
\begin{boxy}{Renormalized Caccioppoli inequality}
\eqn{parallela} 
$$
\begin{cases}
 \, \displaystyle    \int_{B_{\rr/2}} |D(v-\kk)_+|^2 \, dx \\
  \, \displaystyle \qquad   \lesssim_{\data}   \frac{M^{q-p}}{\rr^2} \int_{B_{\rr}}(v-\kk)_+^2\, dx+   M^2\int_{B_{\rr}} |\mu|^2\, dx \\ \\
 \, \displaystyle  \|Du\|_{L^\infty(B_{\rr})} +1 \leq M\,,\quad v \approx |Du|^p\,.
 \end{cases}
$$ 
\end{boxy}
By paying the prize of the appearance of $M$, we now have an inequality which is homogeneous with respect to $v$. Recalling \rif{palline}, the number $M$ is from on chosen to be such that
\eqn{soddisfaM}
$$
\|Du\|_{L^\infty(B_{r_0}(x_0))} +1 \leq M
$$
so that \rif{parallela} holds for every choice of balls $B_{\rr}$ in \rif{palline}. 
 Sobolev embedding theorem now yields 
 \begin{boxy}{Renormalized reverse type inequality with remainder}
\eqn{reversanorm}
    $$
   \left( \mint_{B_{\rr/2}} (v-\kk)_+^{2\chi} \, dx \right)^{1/\chi} \lesssim_{\data} 
      M^{q-p}\mint_{B_{\rr}} (v-\kk)_+^{2} \, dx +  M^2\rr^2\mint_{B_{\rr}} |\mu|^2\, dx
    $$ 
    \end{boxy}
\noindent where (recall here it is $n>2$)
\eqn{classico} 
$$\chi=\frac{n}{n-2}>1\,.$$ 
By \rif{soddisfaM} inequality \rif{reversanorm} holds for every choice of balls as in \rif{palline} and therefore we are allowed to apply Lemma \ref{revlem} with $h=1$, $m_1=t= 2$, $\delta_1=\theta_1= 1$, $M_0=M^{(q-p)/2}$, $M_1=M$ and $\kk_0=0$. This yields
$$
 v(x_{0})  \lesssim_{\data} M^{\frac{n(q-p)}{4}}\left(\mint_{B_{r_{0}}(x_{0})}v^2  \dx\right)^{1/2}
  +M^{\frac{(q-p)(n-2)}{4}+1} \mathbf{P}^{2,1}_{2,1}(\mu;x_{0},2r_{0})\,.
$$
Finally, recalling \rif{triala}, we further estimate 
\begin{flalign}
 \notag |Du(x_{0})| &  \lesssim_{\data} M^{\frac{n}{4}\left(\frac qp-1\right)+\frac 12}\left(\mint_{B_{r_{0}}(x_{0})}[F(Du) +1] \dx\right)^{\frac 1{2p}}\\
 & \qquad \quad 
  +M^{\frac{n-2}{4}\left(\frac qp-1\right)+\frac 1p}\left[ \mathbf{P}^{2,1}_{2,1}(\mu;x_{0},2r_{0})\right]^{\frac 1p} + 1\,.\label{brahms0}
\end{flalign}
All in all, this last estimate holds whenever $B_{2r_0}(x_0)\Subset \Omega$ and $M$ satisfies \rif{soddisfaM}. 
\begin{boxy}{The pointwise-to-global trick}
A general scheme to turn the pointwise estimate \rif{brahms0} into a global one on $B_{r/2}$ getting rid of the presence on $\|Du\|_{L^\infty}$ on the right-hand side, which is implicit in the choice of $M$ in \rif{soddisfaM}.
\end{boxy}
We consider a ball $B_{r}\Subset \Omega$, $r\leq 1$, and concentric balls 
\eqn{concentriche}
$$B_{r/2}\Subset B_{\tau_1}\Subset B_{\tau_2}\Subset B_r\,,$$ and apply \rif{brahms} at the generic point $x_0\in B_{\tau_1}$, with $r_0:=(\tau_2-\tau_1)/2$, and with the choice $M:=\|Du\|_{L^\infty(B_{\tau_2})}+1$, so that \rif{soddisfaM} is met as $B_{2r_0}(x_0)\subset B_{\tau_2}$. We gain
\begin{flalign}
 \notag \texttt{h}(\tau_1) &  \lesssim_{\data}  \frac{[ \texttt{h}(\tau_2)]^{\frac{n}{4}\left(\frac qp-1\right)+\frac 12}}{(\tau_2-\tau_1)^{\frac n{2p}}}\left(\int_{B_r}[F(Du)+1]  \dx\right)^{\frac 1{2p}}\\
 & \qquad \quad 
  + [\texttt{h}(\tau_2)]^{\frac{n-2}{4}\left(\frac qp-1\right)+\frac 1p} \|\mathbf{P}^{2,1}_{2,1}(\mu;\cdot,\tau_2-\tau_1) \|_{L^\infty(B_{\tau_1})}^{1/p}+ 1\,,\label{brahms}
\end{flalign}
where we have set 
\eqn{laHH}
$$
 \texttt{h}(t):=  \|Du\|_{L^\infty(B_{t})} +1\,, \quad r/2 \leq  t \leq r\,.
$$
We now require that
\eqn{condy1}
$$
\frac n4\left(\frac qp-1\right)+\frac 12<1 \Longleftrightarrow \frac qp < 1+\frac 2n
$$
and 
\eqn{condy2}
$$
\frac{n-2}{4}\left(\frac qp-1\right)+\frac 1p <1 \Longleftrightarrow \frac qp < 1+\frac{4(p-1)}{p(n-2)}\,.
$$
These allow to use Young's inequality in \rif{brahms} thereby getting
\begin{flalign*}
 \notag \texttt{h}(\tau_1) &  \leq \frac{\texttt{h}(\tau_2)}{2}+ \frac{c}{(\tau_2-\tau_1)^{\frac{2n}{(n+2)p-nq}}}\left(\int_{B_r}[F(Du) +1]  \dx\right)^{\frac{2}{(n+2)p-nq}}\\
 & \qquad 
  +c\|\mathbf{P}^{2,1}_{2,1}(\mu;\cdot,\tau_2-\tau_1) \|_{L^\infty(B_{\tau_1})}^{\frac{4}{4(p-1)-(n-2)(q-p)}}
\end{flalign*}
with $c\equiv  c (\data)$.  Using Lemma \ref{crit} (recall we are assuming $n>2$ and this allows to verify \rif{lo.2.1}) we estimate 
$$
\|\mathbf{P}^{2,1}_{2,1}(\mu;\cdot,\tau_2-\tau_1) \|_{L^\infty(B_{\tau_1})}\lesssim_{n} 
\|\mu\|_{n,1;B_{\tau_2}} \leq 
\|\mu\|_{n,1;B_{r}}
$$
so that we gain 
\begin{flalign*}
 \notag \texttt{h}(\tau_1) &  \leq \frac{\texttt{h}(\tau_2)}{2}+ \frac{c}{(\tau_2-\tau_1)^{\frac{2n}{(n+2)p-nq}}}\left(\int_{B_r}[F(Du) +1]  \dx\right)^{\frac{2}{(n+2)p-nq}}\\
 & \qquad 
  +c\|\mu\|_{n,1;B_{r}}^{\frac{4}{4(p-1)-(n-2)(q-p)}}\,.
\end{flalign*}
Applying Lemma \ref{l5} below we finally arrive at 
\begin{boxy}{The Lorentz a priori estimate}
\begin{flalign}
 \notag
 \|Du\|_{L^\infty(B_{r/2})}  &\lesssim_{\data}  \left(\mint_{B_r}F(Du) \dx\right)^{\frac{2}{(n+2)p-nq}} 
\\  
& \qquad +\|\mu\|_{n,1;B_{r}}^{\frac{4}{4(p-1)-(n-2)(q-p)}}+1\, .\label{milano}
\end{flalign}
\end{boxy}
Summarising, by \rif{condy1}-\rif{condy2} estimate \rif{milano} holds provided
\eqn{boundfinale}
$$ 
\frac qp < 1+\min\left\{\frac 2n, \frac{4(p-1)}{p(n-2)} \right\}\,.
$$
The one in \rif{milano} is a Lipschitz a priori bound that does not depend on the regularity of $u$ assumed in \rif{appiass} but just of the energy and the right-hand datum $\mu$. It is therefore suitable for a match with approximation arguments bound to remove \rif{appiass}. 
The final outcome is that the implication 
\eqn{implicalo}
$$
\mu \in L(n,1) \Longrightarrow Du \in L^{\infty}_{\loc}(\Omega)
$$ 
and estimate \rif{milano} hold for any general local minimizer $u$. 
\begin{lemma} \label{l5}\, \hspace{-2.5mm} {\em \cite[Lemma 1.1]{gg1}}
\, Let $\textnormal{\texttt{h}}\colon [r/2,r]\to \mathbb{R}$ be a non-negative and bounded function, and let $a,b, \gamma$ be non-negative numbers. Assume that the inequality in 
$$
\textnormal{\texttt{h}}(\tau_1)\le  \frac {\textnormal{\texttt{h}}(\tau_2)}{2}+\frac{a}{(\tau_2-\tau_1)^{\gamma}}+b,
$$
holds whenever $r/2\le \tau_1<\tau_2\le r$. Then 
$$
\textnormal{\texttt{h}}(r/2)\le \frac{c( \gamma)a}{r^{\gamma}}+c(\gamma)b
$$ holds too. 
\end{lemma}
\begin{commentary} {\rm The method exposed applies to general nonuniformly elliptic problems, including those with fast growth conditions as in \rif{veryfast}, see \cite{seminale, ciccio}. We have chosen the $(p,q)$-setting in order to give a streamlined presentation of some of the main ideas, especially those concerned with a priori estimates. We also stress that, although we confined ourselves to the scalar case, the approach presented here works in the vectorial case too when considering functionals with radial structure, i.e., $F(z)\equiv \tilde F(|z|)$. This last assumption is necessary already in the uniformly elliptic case, otherwise counterexamples to Lipschitz regularity emerge; see \cite{Dark, SY} and related references. For this and more general situations we again refer to \cite{seminale, ciccio}. Here the essence of the matter is better understood when $\mu\equiv 0$. In this case \rif{para202} becomes 
\eqn{similar}
$$
    \int_{B_{\rr/2}} |D(v-\kk)_+|^2 \, dx  \lesssim_{\data}  \frac{1}{\rr^2} \int_{B_{\rr}} \mathcal R_{F}(Du)(v-\kk)_+^2\, dx\,.
$$
See again \cite[Lemma 4.5]{seminale}. 
In the uniformly elliptic regime $\mathcal R_{F}(Du)\lesssim 1$ this last inequality reduces to \rif{para20} (with $\mu\equiv 0$) and the standard De Giorgi's iteration applies. In the nonuniform case we pass to the renormalized Caccioppoli inequality
\eqn{renora}
$$
    \int_{B_{\rr/2}} |D(v-\kk)_+|^2 \, dx   \lesssim_{\data}  \frac{M^{q-p}}{\rr^2} \int_{B_{\rr}} (v-\kk)_+^2\, dx\,. 
$$
The nonuniform ellipticity has in a sense disappeared and is now encoded in the presence of $M$. From \rif{renora} we arrive at the following version of \rif{milano}:
\begin{boxy}{Reference Lipschitz estimate}
\eqn{milano2}
$$
\|Du\|_{L^{\infty}(B_{r/2})} \lesssim_{\data}  \left(\mint_{B_r}F(Du)\dx\right)^{\ssf/q}+1
$$
\end{boxy}
\noindent that holds provided 
\eqn{autobb}
$$
\frac qp < 1 + \frac 2n\,.
$$
In \rif{milano2} we find 
\begin{boxy}{The exponent $\ssf$}
\eqn{sceltissime2}
$$
\ssf := 
\begin{cases}
 \displaystyle \frac{2q}{(n+2)p-nq} \quad \mbox{if $n\geq 3$}\\ \\ 
  \displaystyle \mbox{any number $\displaystyle  > \frac{q}{2p-q}$ \quad if $n=2$ and $q>p$} \\ \\ 
  1  \quad \mbox{if $n=2$ and $q=p$}
\,.
\end{cases}
$$
\end{boxy}
\vspace{3mm}
This is a very streamlined and short proof of Theorem \ref{marcth}, and \rif{milano2} is exactly the a priori estimate found by Marcellini in his classic paper \cite{M2} (our proof works in this case also when $n>2$ as $\mu \equiv 0$). In fact, by passing to renormalized energy inequalities like \rif{renora}, we are in position to use uniformly elliptic iteration procedures. Notice that, at this stage, it is crucial to get a precise control of the constants $M_i$ in Lemma \ref{revlem} in order to conclude, an information that is not needed in the uniformly elliptic case \rif{para20}-\rif{reversa}, as seen in Section \ref{toy}. Surprisingly, this single fact allows to reduce the analysis of nonuniformly elliptic problems to that of uniformly elliptic ones. This idea is general and can be implemented also via Moser's iteration method; see Section \ref{mosersec}. Moreover, the method above is general enough to provide another proof of Theorem \ref{bstheorem}. This has been done by Bella \& Sch\"affer \cite{BS2}, who cleverly combined the methods of \cite{BS0, BS} and those from \cite{seminale} to prove \rif{implicalo} under the improved bound 
$$
\frac qp < 1+\min\left\{\frac 2{n-1}, \frac{4(p-1)}{p(n-3)} \right\}
$$
for $n>3$; this improves \rif{boundfinale}. When $\mu =0$ this bound reduces to \rif{boundbs}. We finally remark that in the case $p=q$ (uniform ellipticity) estimate \rif{milano} reduces to
$$
 \|Du\|_{L^\infty(B_{r/2})}  \lesssim_{\data}   \left(\mint_{B_r}|Du|^p  \dx\right)^{1/p} +\|\mu\|_{n,1;B_{r}}^{1/(p-1)}+1
$$
which the standard a priori estimate for $p$-Laplacean type problems \cite{DMpoinc, KMstein, KMvec}. 
\begin{boxy}{The key takeaway}
The analysis of nonuniformly elliptic problems can be reduced to that of uniform ones via suitable renormalized Caccioppoli inequalities \rif{parallela}. These are aimed at overcoming the lack of homogeneous estimates which are typical of the nonuniform case. The mechanism works thanks to a precise quantitative control of the constants along the relevant iterative schemes (Lemma \ref{revlem}). In this setting nonlinear potentials naturally appear along the iterations and encode in a sharp way the regularity properties of the data. 
 \end{boxy}
}\end{commentary}

\subsection{Fractional approach to nonuniform ellipticity \cite{piovra}}\label{fracschau} As already mentioned, Lipschitz regularity is the focal point in nonuniformly elliptic problems. The classical approach to $C^{0,1}$-regularity is to start differentiating the Euler-Lagrange equation \rif{elf} (in the $x_s$ direction, $1\leq s \leq n$). The resulting differentiated equation formally looks like \rif{lineare2}, with $v\equiv D_su$, $\texttt{A}(x)\equiv \partial_{zz}F(Du(x))$ and $\mu$ replaced by $\partial_{x_s}\mu$ (in the distributional sense). This leads to the basic energy inequality \rif{para202}, which is the starting result for the subsequent developments. Here we want to address what happens in the case of functionals with H\"older continuous coefficients, for simplicity looking at the model case \rif{modello0} and Theorem \ref{t2}, whose proof we are indeed going to sketch. Therefore from now on we shall assume that 
\eqn{bound333-dopo}
$$
\frac{q}{p}\leq  1+  \frac{1}{5}\left(\frac{\alpha}{n}\right)^2\,,
$$
holds. Moreover, as explained in \cite{piovra}, we can confine ourselves on the case
\eqn{alffa}
$$
0<\alpha<1\,.
$$
On the other hand, the differentiable case $\alpha=1$ - Lipschitz continuous coefficients - can also be dealt with by more standard methods, see also Section \ref{mosersec}. 
In the case of Theorem \ref{t2} the Euler-Lagrange equation of the functional is therefore
\eqn{elf2}
$$-\divo \, [\ccc(x)\partial_zF(Du)]=\mu\,.$$ 
This time we cannot differentiate \rif{elf2} as $\ccc(\cdot)$ is only H\"older continuous, so that the previous approach fails immediately unless $\ccc(\cdot)$ is differentiable, i.e., $\alpha=1$. This is in fact the point leading to consider indirect, perturbation methods to prove gradient regularity, and this is the classical basic strategy pursued in the proof of Schauder estimates. As mentioned before, this approach is succesful in the uniformly elliptic case and falls short of working in nonuniform problems. We therefore go back to the direct approach, but we perform a ``fractional differentiation" of \rif{elf2}. In other words, while Lipschitz continuous coefficients lead to energy estimates involving second order derivatives of solutions as \rif{para202}, H\"older coefficients lead to energy estimates featuring higher order fractional derivatives of solutions. This approach, developed in \cite{piovra}, is really too long to reproduce here and we confine ourselves to give a sketch of it; moreover, we shall confine for simplicity to the case $p\geq 2$. For this, we start recalling the basic definition of fractional Sobolev spaces. In this case, we use so-called Sobolev-Slobodevsky spaces, equipped with Gagliardo norms \cite{AF, guide}.  
\begin{definition}\label{fra1def}
Let $\beta \in (0,1)$, $N \in \en$, $n \geq 2$, and let $\Omega \subset \er^n$ be an open subset. The space $W^{\beta ,2}(\Omega,\er^N )$ consists of maps $w \colon \Omega \to \er^N$ such that
$$
\| w \|_{W^{\beta ,2}(\Omega )}  := \|w\|_{L^2(\Omega)}+ \left(\int_{\Omega} \int_{\Omega}  
\frac{|w(x)- w(y)|^{2}}{|x-y|^{n+2\beta}} \dx \dy \right)^{1/2}\,.
$$
The local variant $W^{\beta ,2}_{\loc}(\Omega,\er^N )$ is defined by requiring that $w \in W^{\beta ,2}_{\loc}(\Omega,\er^N )$ iff $w \in W^{\beta ,2}(\tilde{\Omega},\er^N)$ for every open subset $\tilde{\Omega} \Subset \Omega$. 
\end{definition}
Fractional spaces come along with their own embedding; see for instance \cite[Theorem 6.7]{guide}.  We shall use it on balls $B_{\rr}\subset \er^n$ in the form
\begin{flalign}
\notag \left(\mint_{B_{\rr}}|w|^{\frac{2n}{n-2\beta}}\dx \right)^{\frac{n-2\beta}{2n}} &\lesssim_{n,\beta}  
\rr^{\beta} \left( \int_{B_{\varrho}} \mint_{B_{\varrho}}  
\frac{|w(x)- w(y)|^{2}}{|x-y|^{n+2\beta}} \dx \dy\right)^{1/2}\\
& \qquad \quad +  \left(\mint_{B_{\rr}}|w|^{2}\dx \right)^{1/2} \,, \label{immersione}
\end{flalign}
that holds whenever $w\in W^{\beta,2}(B_{\rr})$. This follows via a simple scaling argument.   

Proceeding as in the scheme of a priori estimates for more regular minimizers, we shall now consider a local minimizer $u$
of the functional which is such that
\eqn{assuntoc1} 
$$u \in C^1_{\loc}(\Omega)\,.$$
As usual, such an additional regularity assumption can be removed by an approximation procedure (see \cite{piovra}).  With the function $v$ being defined as in \rif{triala}, we have
\begin{boxy}{Renormalized fractional Caccioppoli inequality}
\begin{flalign}\label{xx.10}
& \int_{B_{\varrho/2}} \mint_{B_{\varrho/2}}  
\frac{|(v-\kk)_{+}(x)- (v-\kk)_{+}(y)|^{2}}{|x-y|^{n+2\beta}} \dx \dy
\notag \\ & \quad  \lesssim_{\data,\alpha, \beta} \frac{M^{\ssf (q-p)}}{\varrho^{2\beta}}\mint_{B_{\rr}}(v-\kk)_{+}^{2} \dx \notag  \\
& \qquad  \quad +\frac{M^{\ssf q+p-\mathrm{b}}}{\varrho^{2\beta}}\rr^{2\alpha}\mint_{B_{\rr}}(|Du|+1)^{2q-2p+\mathrm{b}} \dx\,.
\end{flalign}
\end{boxy}
This holds whenever $p\geq 2$, $B_{\rr}\Subset \Omega$, $\kk \geq 0$, and for any choice of numbers $M, \mathrm{b},  \beta $ obeying
\eqn{sceltissime}
$$
 \|Du\|_{L^{\infty}(B_{\rr})}+1\leq M,
\qquad 0 \leq  \mathrm{b} \leq p,\qquad  0< \beta< \frac{\alpha}{1+\alpha}\,,
$$
where $\ssf$ has been defined in \rif{sceltissime2}. 
\begin{remark}[The exponent $\ssf$ in  \rif{sceltissime2}] Note that definition of $\ssf$ makes sense as $(n+2)p-nq>0$, and this is a consequence of \rif{autobb} (which is of course implied by \rif{bound333-dopo}). 
As one could guess from its shape, the derivation of \rif{xx.10}, as an a priori estimate, does not yet require the bound in \rif{bound333-dopo}. At this stage we only need to assume
\rif{autobb}. Also note that
\eqn{sopra}
$$
\ssf \geq 1  \Longleftrightarrow 1 \leq \frac qp <  1 + \frac 2n $$
and that $\ssf=1$ iff $p=q$. The exponent $\ssf$ plays a special role here, as, in fact,  estimate \rif{milano2} is a basic building block in the proof of \rif{xx.10}. See the Commentary \ref{commentami} below. 
\end{remark}
Inequality \rif{xx.10} is an analog of \rif{parallela}. In both cases we have a Caccioppoli type estimate where the integral of a higher derivative of $(v-\kk)_{+}$ is controlled by the integral of $(v-\kk)_{+}$. In the case of \rif{parallela} we have a first derivative of $(v-\kk)_{+}$, while in \rif{xx.10} we have a fractional derivative. The last terms are concerned wit the so-called {\em external ingredients}: right-hand side data $\mu$ when looking at \rif{parallela}, presence of coefficients when dealing with \rif{xx.10}. Both inequalities are homogeneous with respect to $(v-\kk)_{+}$, i.e., they have been renomalized hiding the nonuniform ellipticity using the constant $M$. A main difference between the two is the presence of the free renormalization parameter $\mathrm{b}$ in \rif{xx.10}, that will be used later on. We can now apply \rif{immersione} with $w\equiv (v-\kk)_{+}$ in \rif{xx.10} and the result is
 \begin{boxy}{Renormalized fractional reverse type inequality with remainder}
\begin{flalign}\label{xx.101}
\notag  \left(\mint_{B_{\rr/2}}(v-\kk)_{+}^{2\chi}  \dx\right)^{1/\chi}  &\lesssim_{\data,\alpha, \beta}  M^{\ssf (q-p)} \mint_{B_{\rr}}(v-\kk)_{+}^{2} \dx \\ & \ \quad +M^{\ssf q+p-\mathrm{b}}\rr^{2\alpha}\mint_{B_{\rr}}(|Du|+1)^{2q-2p+\mathrm{b}} \dx
\end{flalign}
\end{boxy}
\noindent where
\eqn{classico2} 
$$
\chi\equiv \chi(\beta):= \frac{n}{n-2\beta}>1\,.
$$
Again, \rif{xx.101} parallels \rif{reversanorm}, with a different value of $\chi$ (compare \rif{classico} with \rif{classico2}). To proceed, we fix balls as in \rif{concentriche}, we take $$M:=\|Du\|_{L^\infty(B_{\tau_2})}+1$$ and use \rif{brahms} at the generic point $x_0\in B_{\tau_1}$, with $r_0:=(\tau_2-\tau_1)/2$, so that $B_{2r_0}(x_0)\subset B_{\tau_2}$. This leads us to apply Lemma \ref{revlem}, with $h=1$, $m_1\equiv 2q-2p+\mathrm{b}$, $\theta_1\equiv 1$, $t\equiv 2$, $\delta_1\equiv \alpha$, $M_0\equiv M^{\ssf(q-p)/2}$, $M_1\equiv M^{(\ssf q+p-\mathrm{b})/2}$, $\kk_0\equiv 0$, thereby getting 
\begin{flalign*}
\nonumber
 v(x_{0}) & \lesssim_{\data,\alpha, \beta}  M^{\frac{\chi\ssf(q-p)}{2(\chi-1)}}\left(\mint_{B_{r_{0}}(x_{0})}v^{2}  \dx\right)^{1/2}\\
 & \quad 
\ \ +M^{\frac{\ssf(q-p)}{2(\chi-1)}+\frac{\ssf q+p-\mathrm{b}}{2}}\mathbf{P}^{2q-2p+\mathrm{b},1}_{2,
\alpha}(|Du|+1;x_{0},2r_{0})
\end{flalign*}
so that, recalling \rif{triala} 
\begin{flalign*}
\nonumber
|Du(x_{0})| & \lesssim_{\data,\alpha, \beta}  M^{\frac{\chi\ssf}{2(\chi-1)}\left(\frac qp-1\right)+\frac 12}\left(\mint_{B_{r_{0}}(x_{0})}(|Du|^p  +1)\dx\right)^{\frac 1{2p}}\\
 & \quad 
\ \ +M^{\frac{\ssf}{2(\chi-1)}\left(\frac qp-1\right)+\frac{\ssf q+p-\mathrm{b}}{2p}}\left[\mathbf{P}^{2q-2p+\mathrm{b},1}_{2,
\alpha}(|Du|+1;x_{0},2r_{0})\right]^{1/p}+1\,.
\end{flalign*}
The inequalities in the last two displays hold for every ball $B_{2r_0}(x_0)\Subset \Omega$ as we are assuming that $Du$ is continuous \rif{assuntoc1}. We can now use the covering argument in concentric balls \rif{concentriche} as after \rif{brahms}, thereby obtaining 
\begin{eqnarray}
&& \label{brahms2} \qquad \texttt{h}(\tau_1)  \\
 && \notag   \lesssim_{\data,\alpha, \beta}  \frac{[ \texttt{h}(\tau_2)]^{\frac{\chi\ssf}{2(\chi-1)}\left(\frac qp-1\right)+\frac 12}}{(\tau_2-\tau_1)^{\frac{n}{2p}}}\left(\int_{B_r}[F(Du)+1]  \dx\right)^{\frac 1{2p}}\\
  && \notag 
  + [\texttt{h}(\tau_2)]^{\frac{\ssf}{2(\chi-1)}\left(\frac qp-1\right)+\frac{\ssf q+p-\mathrm{b}}{2p}} \|\mathbf{P}^{2q-2p+\mathrm{b},1}_{2,
\alpha}(|Du|+1;\cdot,\tau_2-\tau_1) \|_{L^\infty(B_{\tau_1})}^{1/p}+ 1
\end{eqnarray}
where $\texttt{h}(\cdot)$ is again defined as in \rif{laHH} and remark that $\mathrm{b}$ is still a free parameter.  
Now we want to proceed as after \rif{brahms}. For this, we first want to check that that the $L^\infty$-norm of 
$\mathbf{P}^{2q-2p+\mathrm{b},1}_{2,
\alpha}$ is finite. This means applying Lemma \ref{crit} and therefore verifying that
\eqn{condif1}
$$
\frac{n\theta }{t\delta}\equiv \frac{n}{2\alpha}>1\quad \mbox{and}\quad
\frac{mn\theta}{t\delta}\equiv \frac{n(2q-2p+\mathrm{b})}{2\alpha}< p\,.
$$
Note that the first inequality in the above display is automatically verified by \rif{alffa} and that the second one is equivalent to require
\eqn{condif1-dopo}
$$
\frac qp < 1+ \frac{\alpha}{n} -\frac{\mathrm{b}}{2p}\,.
$$
These ensure that
\begin{flalign}
\notag  &\|\mathbf{P}^{2q-2p+\mathrm{b},1}_{2,
\alpha}(|Du|+1;\cdot,\tau_2-\tau_1) \|_{L^\infty(B_{\tau_1})}  \leq c\|Du\|_{L^{p}(B_{r})}^{q-p+\mathrm{b}/2}+c\\
&\qquad \qquad\qquad \qquad\qquad\qquad\leq c \left(\int_{B_{r}} \ccc(x)F(Du)\dx \right)^{\frac qp-1+\frac{\mathrm{b}}{2p}} +c\notag \\
& \qquad \qquad\qquad\qquad\qquad\qquad\leq c \left(\int_{B_{r}} F(Du)\dx \right)^{\frac qp-1+\frac{\mathrm{b}}{2p}} +c\label{condif11}
\end{flalign}
holds with $c\equiv c(\data,\mathrm{b})$ as a consequence of \rif{11}. Next, we have to make sure that the exponents of $ \texttt{h}(\tau_2)$ in \rif{brahms2} are smaller than one, i.e., 
\eqn{prima}
$$
\begin{cases}
\displaystyle \frac{\chi(\beta)\ssf}{2(\chi(\beta)-1)}\left(\frac qp-1\right)+\frac 12<1\\
\displaystyle  \frac{\ssf}{2(\chi(\beta)-1)}\left(\frac qp-1\right)+\frac{\ssf q+p-\mathrm{b}}{2p}<1\,.
\end{cases}
$$
In fact, the second inequality in the above display implies the first and reduces to
\eqn{condif2}
$$
 \left(\frac qp-1\right)\frac{\ssf n}{4\beta}+\frac{\ssf}{2}-\frac{\mathrm{b}}{2p}
 < \frac 12\,.
$$
In terms of $q/p$, condition \rif{condif2} translates into
\eqn{condif2-dopo}
$$
\frac qp < 1 + \left(1-\ssf+\frac{{\rm b}}{p}\right) \frac{2\beta}{\ssf n}\,.
$$
Notice that in \rif{condif1-dopo} and \rif{condif2-dopo} we can still choose the values of $\beta$ and $\mathrm{b}$ appearing in \rif{sceltissime}. To do this, note that while the right-hand side of \rif{condif1-dopo} is a decreasing function of ${\rm b}$, the right-hand side of \rif{condif2-dopo} is increasing. We then optimize the choice of ${\rm b}$ by equalizing the two expressions and this leads to choose 
\eqn{ilb}
$$
\mathrm{b}\equiv \mathrm{b}(\beta)=\frac{2p[\ssf \alpha  +2\beta(\ssf-1)]}{\ssf n+4\beta}>0\,.
$$
Notice that this is an admissible value provided ${\rm b}\leq p$. Recalling the definition of $\ssf$ in \rif{sceltissime2}, we note that 
\eqn{ammissibile}
$$
 \mathrm{b}  < \frac{2p\alpha}{n} \Longleftrightarrow \mathfrak{s} < 1 +\frac{2\alpha}{n}\,.
$$
On the other hand we have
\eqn{notabase}
$$
\trif{bound333-dopo} \Longrightarrow 
 \ssf < 1 + \frac{\alpha}{2n} 
$$
and this proves the admissibility of ${\rm b}$ in \rif{ammissibile}. 
Notice that condition \rif{ammissibile} is independent of $\beta$. In fact, 
\eqn{decrescente}
$$ \mbox{$t \in [0,p]\mapsto {\rm b}(t)$ is decreasing $\displaystyle \Longleftrightarrow \mathfrak{s} < 1 +\frac{2\alpha}{n}$}$$
and therefore
$p \geq {\rm b}(0)= 2p\alpha/n> {\rm b}(\beta)$, which is again \rif{ammissibile}. By inserting the value of ${\rm b}$ found in \rif{ilb} in the right-hand side of \rif{condif1-dopo} (which is the same that inserting ${\rm b}$ in \rif{condif2-dopo}), gives the unified condition
\eqn{bound3veroprepre}
$$
\frac qp < 1 + \frac{2\beta}{n}\left[\frac{2\alpha-n(\ssf -1)}{\ssf n+4\beta}
\right]\,.
$$
Notice that the right-hand side is an increasing function of $\beta$ as consequence of \rif{decrescente}. This leads us to consider the (forbidden) limiting case 
$$
\beta = \frac{\alpha}{1+\alpha}
$$
in \rif{bound3veroprepre}, thereby getting 
$$
\frac qp < 1 + \frac{2\alpha}{(1+\alpha)n}\left[\frac{2\alpha-n(\ssf -1)}{\ssf n+4\alpha/(1+\alpha)}
\right]=: 1+\mathcal{R}_1(n, p,q,\alpha) \frac{\alpha^2}{n^2}
$$
with
$$
 \mathcal{R}_1(n, p,q,\alpha)= \frac{2}{1+\alpha}\left[\frac{2-n(\ssf -1)/\alpha}{\ssf +\frac{4\alpha}{n(1+\alpha)}}
\right]\,.
$$
Using \rif{notabase} we observe that $ \mathcal{R}_1> 2/3$. It follows that we can choose $\beta < \alpha/(1+\alpha)$ close enough to $\alpha/(1+\alpha)$ to guarantee that \rif{bound3veroprepre} holds so that both \rif{condif1-dopo}
 and \rif{prima} hold too, with the corresponding choice of ${\rm b}\equiv {\rm b(\beta)}$ in \rif{ilb}. We can apply Young's inequality in \rif{brahms2}, that, together with \rif{condif11}, yields
$$
 \notag \texttt{h}(\tau_1)   \le  \frac{\texttt{h}(\tau_2)}{2}+ \frac{c}{(\tau_2-\tau_1)^{n\kk }}\left(\int_{B_{r}}[F(Du) +1]  \dx\right)^{\kk}
$$
with $c\equiv c (\data, \alpha)$ and $\kk\equiv \kk(n,p,q,\alpha)>1$, and with $\texttt{h}(\cdot)$ being defined as in \rif{laHH}. Applying Lemma \ref{l5} we conclude with  
$$
\|Du\|_{L^\infty(B_{r/2})}   \lesssim_{\data, \alpha}  \left(\mint_{B_{r}} F(Du)  \dx\right)^{\kk} +1 
$$
that is \rif{stimascha} under the additional a priori regularity assumption \rif{assuntoc1}. As already mentioned, this can be removed via an approximation argument \cite[Corollary 3]{piovra} so that we obtain the full statement of Theorem \ref{t2}. Full proofs, including the case $1<p<2$, can be found in \cite{piovra}. 
\begin{commentary}\label{commentami}{\rm A crucial point in the proof of Theorem \ref{t2} is the fractional energy estimate \rif{xx.10}, that, as already mentioned, cannot be obtained via direct differentiation of \rif{elf2}, on the contrary of \rif{para202}. Inequality \rif{xx.10} is obtained via a sort of nonlinear dyadic/atomic decomposition technique, finding its roots in \cite{KM}, and that resembles the one used for Besov spaces \cite{AdHe, triebel} in the setting of Littlewood-Paley theory. This roughly goes as follows (see \cite{piovra} for the full details). First we rescale the problem to $B_1(0)$ by passing to 
\eqn{riscale} 
$$x\to \frac{u(x_0+\rr x)}{\rr}\,.$$ For simplicity we shall still denote this function by $u$. We select a scale of differentiation $0<|h|\ll1$, where $h \in \er^n$, and consider small cubes and balls with radius $|h|^{1/(1+\alpha)}\gg |h|$. 
The idea is to cover $B_{1/2}$ with a lattice of disjoint dyadic cubes $\{Q_k\equiv Q_h(x_k)\}_{k \leq \mathfrak{n}}$ with (equal) sidelength $\approx |h|^{1/(1+\alpha)}$, centred at points $x_k \in B_{1/2}$ and with sides parallel to the coordinate axes. The number $\mathfrak{n}$ of such cubes is comparable to $|h|^{-n/(1+\alpha)}$. It follows that any dilated family $\{tQ_k\}$, $t\geq 1$ has the finite intersection property. This means that each cube of the type $tQ_k$ does not touch others of the same type more that a finite number $\texttt{n}\equiv \texttt{n}(n,t)$ of times. The same happens for the outer balls $\{B_k\equiv B_h(x_k)\}$, defined as the smallest concentric balls to $Q_k$ such that $Q_k \subset B_k$ (for this it is sufficient to consider the corresponding outer cubes, i.e., the smallest concentric cubes containing $B_k$; these are still of the type $tQ_k$, with $t= \sqrt{n}$). We take $|h|$ small enough to guarantee that $8B_k\Subset B_1$ (recall that $B_k$ is centred in $B_{1/2}(0)$). All in all, we have 
$$
\Big| \ B_{1/2}(0)\setminus \bigcup_{k\le \mathfrak{n}}Q_{k} \ \Big|=0,\qquad Q_{i}\cap Q_{j}=\emptyset \ \Leftrightarrow \ i\not =j
$$
and
\eqn{finite}
$$
\lambda(B_{1/2}(0))\leq \sum_{k\leq  \mathfrak n}\lambda(B_k) \leq \sum_{k\leq  \mathfrak n}\lambda(8B_k) \lesssim_n \lambda(B_1(0))\,,
$$
that holds for any Borel measure $\lambda(\cdot)$ with finite total mass defined on $B_1(0)$. 
Next, on each ball $8B_k$, one defines $v_k\equiv v(B_k)\in W^{1,q}(8B_k)$ solving  
\eqn{elf2x}
$$
\begin{cases}
\, -\divo \, [\ccc(x_k)\partial_zF(Dv_k)]=0\\
\, v_k\equiv u\  \mbox{on $\partial 8B_k$}\,.
\end{cases}
$$ 
The function $v_k$ plays the role of an ``atom" in the sense of the dyadic decompositions used to decompose functions Besov spaces \cite{AdHe, triebel}. In a very rough sense, we shall see that the atomic decomposition in question looks like
$$
u(x) \approx \sum_k v_k(x) 1_{B_k}(x) + \texttt{o}(|h|)
 $$
see \rif{cccc1} below, where $\texttt{o}(|h|)$ is measured in a suitable sense. This time we employ atoms that are more special functions than those considered in the usual Besov spaces decompositions (see for instance \cite[Definition 4.6.1]{AdHe}). Indeed, they  themselves are solutions to nonlinear equations \rif{elf2x}. This allows to such a decomposition to adapt more closely to the original problem, which is nonlinear. To proceed, we apply to $v_k$ a variant of the classical Bernstein technique, in the sense that we get the following energy inequality
\begin{flalign}
 \notag & \int_{2B_k} |D(G(|Dv_k|) - \kk)_+|^2 \dx\\
 & \qquad 
  \lesssim 
  \frac{1}{|h|^{\frac{2}{1+\alpha}}} \int_{4B_k} \mathcal R_{F}(Dv_k) (G(|Dv_k|) - \kk)_+^2  \dx
  \label{lavecchia}
\end{flalign}
for every $\kk \geq 0$, which is in fact \rif{similar}; see \cite[Lemma 5.1]{piovra} and recall that $G(\cdot)$ is defined in \rif{triala}. Here and until the end of this section all the constants implied in the symbol $\lesssim$  will depend at most on $\data$ and $\alpha$. By standard properties of finite differences we get  
\eqn{frompi0}
$$
\int_{B_k} |\tau_{h}(G(|Dv_k|) - \kk)_+|^2 \dx \lesssim |h|^2\int_{2B_k} |D(G(|Dv_k|) - \kk)_+|^2 \dx
$$ 
where, as usual, 
  $
  \tau_{h}w(x):= w(x+h)-w(x)
  $  denotes the standard finite difference operator of a function $w$ in direction $h\in \er^n$. Note that here we have used that 
  \eqn{pppalle}
  $$|h|\leq |h|^{1/(1+\alpha)} \Longrightarrow B_k+h \subset 2B_k\,.$$
Combining \rif{frompi0} and \rif{lavecchia} we gain
\begin{flalign}
 \notag &\int_{B_k} |\tau_{h}(G(|Dv_k|) - \kk)_+|^2 \dx
  \\
  &\quad \lesssim |h|^\frac{2\alpha}{{1+\alpha}}  \int_{4B_k} \mathcal R_{F}(Dv_k) (G(|Dv_k|) - \kk)_+^2  \dx\, . \label{frompi}
\end{flalign}
Using additional a priori estimates for solutions to \rif{elf2x}, that is, applying \rif{milano2} to $v_k$, then gives
$$
\|Dv_k\|_{L^\infty(4B_k)} \lesssim \|Du\|_{L^\infty(8B_k)}^{\ssf}+1 \lesssim M^{\ssf}
$$
where $\ssf$ is in \rif{condif2} and 
provided $M$ is as in \rif{sceltissime}, \cite[Proposition 5.3]{piovra}. Recalling \rif{para2022}, and using the last inequality in \rif{frompi},  we find 
\begin{flalign}
\notag & \int_{B_k} |\tau_{h}(G(|Dv_k|) - \kk)_+|^2 \dx
 \\
 &\quad  \lesssim |h|^\frac{2\alpha}{{1+\alpha}} M^{(q-p)\ssf} \int_{4B_k} (G(|Dv_k|) - \kk)_+^2  \dx\,.\label{cccc}
\end{flalign}
The one in \rif{cccc} is a sort of Nikolski space estimate at the microlocal level $B_k$. Indeed, recall that a function $w\in L^{2}(\Omega)$ belongs to Nikolski space $\mathcal N^{s,2}(\Omega)$ iff 
\eqn{nikki}
$$
 \int_{\Omega_{|h|}} |\tau_{h}w|^2 \dx \lesssim |h|^{2s}\,,\qquad \forall \ h \in \er^n\,,
$$
where $\Omega_{h}:= \{x \in \Omega \colon \dist(x, \partial \Omega)>|h|\}$. Such spaces relate  well with Sobolev-Slobodevski spaces (see Lemma \ref{l4} below and \cite{AKM} for more details); all of them are special cases of a larger family of fractional spaces called Besov spaces \cite{triebel}. Here the word microlocal refers to the fact that in \rif{cccc} the differentiability gain $|h|^{2\alpha/(1+\alpha)}$ is obtained on a domain $B_k$ whose diameter is proportional to $|h|^{1/(1+\alpha)}$, and therefore depends on $|h|$, on the contrary of the classical definition of Nikolski space in \rif{nikki}. A point worth remarking here is that the gain of the factor $|h|^{2\alpha/(1+\alpha)}$ comes from the combination of \rif{lavecchia} and \rif{frompi0}; in this respect choosing the size of the sidelength of $Q_k$ to be comparable to $|h|^{1/(1+\alpha)}$, and therefore larger to the differentiation scale $|h|$, is crucial. The special exponent $1/(1+\alpha)$ optimizes the procedure.  

Estimate \rif{cccc} is eventually transferred to to $u$ via the comparison estimate
\begin{flalign}
\notag & \int_{4B_k}|(G(|Du|)-\kk)_{+}-(G(|Dv_k|)-\kk)_{+}|^{2}\dx
\\
&\qquad \quad  \lesssim |h|^\frac{2\alpha}{{1+\alpha}}M^{\mathfrak{s}p+p-\mathrm{b}}\rr^{2\alpha}\int_{8B_k}(|Du|+1)^{2q-2p+\mathrm{b}} \dx\,.\label{cccc1}
\end{flalign} 
We transfer the microlocal information in \rif{cccc} from $v_k$ to $u$ on $B_k$, as follows 
\begin{flalign*}
\notag  & \int_{B_k}|\tau_{h}(G(Du)-\kk)_{+}|^{2}  \dx    \lesssim  \int_{B_k}|\tau_{h}(G(Dv_k)-\kk)_{+}|^{2}\dx \\
&  \qquad\qquad\qquad+ \int_{B_k}|(G(Du(\cdot+h))-\kk)_{+}-(G(Dv_k(\cdot+h))-\kk)_{+}|^{2}\dx\\
&\qquad\qquad\qquad+ \int_{B_k}|(G(Du)-\kk)_{+}-(G(Dv_k)-\kk)_{+}|^{2}\dx\\
& \qquad\qquad\quad \lesssim \int_{B_k}|\tau_{h}(G(Dv_k)-\kk)_{+}|^{2}\dx \\
&  \qquad\qquad\qquad + \int_{2B_k}|(G(Du)-\kk)_{+}-(G(Dv_k)-\kk)_{+}|^{2}\dx\,.
\end{flalign*}
Note that here we have again used \rif{pppalle}. Using \rif{cccc} and \rif{cccc1} we then find 
\begin{flalign*}
\notag   \int_{B_k}|\tau_{h}(G(Du)-\kk)_{+}|^{2}  \dx &   \lesssim |h|^\frac{2\alpha}{{1+\alpha}} M^{(q-p)\ssf} \int_{4B_k} (G(|Dv_k|) - \kk)_+^2  \dx\\
& \quad  + |h|^\frac{2\alpha}{{1+\alpha}}M^{\mathfrak{s}p+p-\mathrm{b}}\rr^{2\alpha}\int_{8B_k}(|Du|+1)^{2q-2p+\mathrm{b}} \dx
\end{flalign*}
and finally applying \rif{cccc1} once again, we conclude with 
\begin{flalign*}
\notag  \int_{B_k}|\tau_{h}(G(Du)-\kk)_{+}|^{2}  \dx 
 & \lesssim |h|^\frac{2\alpha}{{1+\alpha}}M^{\mathfrak{s}(q-p)}\int_{8B_k}(G(Du)-\kk)_{+}^{2}  \dx\\
 &\quad +|h|^\frac{2\alpha}{{1+\alpha}}M^{\mathfrak{s}q+p-\mathrm{b}}\rr^{2\alpha}\int_{8B_k}(|Du|+1)^{2q-2p+\mathrm{b}} \dx.
\end{flalign*}
Adding up the above inequalities over the decomposition $\{B_k\}$, and using \rif{finite}, we obtain
\begin{flalign*}
\notag & \int_{B_{1/2}(0)}|\tau_{h}(G(Du)-\kk)_{+}|^{2}  \dx 
 \lesssim |h|^\frac{2\alpha}{{1+\alpha}}M^{\mathfrak{s}(q-p)}\int_{B_1(0)}(G(Du)-\kk)_{+}^{2}  \dx\\
 &\qquad \qquad \quad \qquad  \ \qquad \qquad +|h|^\frac{2\alpha}{{1+\alpha}}M^{\mathfrak{s}q+p-\mathrm{b}}\rr^{2\alpha}\int_{B_1(0)}(|Du|+1)^{2q-2p+\mathrm{b}} \dx\,.
\end{flalign*}
We are now in the position to use Lemma \ref{l4} below, that gives 
\begin{flalign*}
& \int_{B_{1/2}(0)} \mint_{B_{1/2}(0)}  
\frac{|(G(Du)-\kk)_{+}(x)- (G(Du)-\kk)_{+}(y)|^{2}}{|x-y|^{n+2\beta}} \dx \dy
\notag \\ & \qquad  \qquad  \lesssim M^{\ssf (q-p)}\mint_{B_{1}(0)}(G(Du)-\kk)_{+}^{2} \dx \\
& \qquad  \qquad   \qquad+M^{\ssf q+p-\mathrm{b}}\rr^{2\alpha}\mint_{B_{1}(0)}(|Du|+1)^{2q-2p+\mathrm{b}} \dx\,.
\end{flalign*}
This holds whenever $0<\beta< \alpha/(1+\alpha)$. From this \rif{xx.10} follows rescaling as in \rif{riscale}. 
\begin{lemma}\, \hspace{-2.5mm} {\em \cite[Lemma 1]{dmaamp}}\label{l4}
Let $B_{\varrho} \Subset B_{r}\subset \er^n$ be concentric balls with $r\leq 1$, $w\in L^{2}(B_{r},\mathbb{R}^{k})$, $s\geq 1$ and assume that, for $\beta \in (0,1]$, $\mathcal H\ge 1$, there holds
$$
\|\tau_{h}w\|_{L^{2}(B_{\rr})}\le \mathcal H|h|^{\beta } \quad \mbox{
for every $h\in \mathbb{R}^{n}$ with $0<|h|\le \frac{r-\rr}{K}$, where $K \geq 1$}.
$$
Then it holds that 
\begin{flalign*}
&\|w\|_{W^{\alpha_{0},2}(B_{\rr})}\\
& \qquad \lesssim_{n,s}\frac{1}{(\beta -\alpha_{0})^{1/2}}
\left(\frac{r-\rr}{K}\right)^{\beta -\alpha_{0}}\mathcal H+\left(\frac{K}{r-\rr}\right)^{n/2+\alpha_{0}} \|w\|_{L^{2}(B_{\rr})}\,.
\end{flalign*}
\end{lemma} 

The approach outlined here is robust and flexible enough to allow for several applications.   The point is to choose the right form of the gradient function $v$ playing the role of subsolution in the setting of Bernstein technique. An example is given in \cite{loggy}, where, in order to prove Theorems \ref{t1log} and \ref{t3log}, a fractional Caccioppoli inequality has been derived for the function
\eqn{intrinseca}
$$
v = |Du|^{2-\gamma} + \aaa(\cdot)|Du|^q\,,
$$
where $\gamma \in [1,2)$ is suitably close to $1$ as described in Theorem \ref{t3log}. One takes $\gamma=1$ in the case of Theorem \ref{t1log} and of the model functional in \rif{modello}.  
The choice \rif{intrinseca} reflects the structure of the functional considered in \rif{generale}. Specifically, the scheme outlined above works with an inequality of the type
\begin{eqnarray}
&& \int_{B_{\varrho/2}} \mint_{B_{\varrho/2}}  
\frac{|(v-\kk)_{+}(x)- (v-\kk)_{+}(y)|^{2}}{|x-y|^{n+2\beta}} \dx \dy
\notag \\ 
&& \qquad  \lesssim \frac{M^{2\ssss_{1}}}{\rr^{2\beta}}\mint_{B_{\rr}}(v-\kk)_{+}^{2} \dx\nonumber\\
&& \qquad \quad \notag +\frac{M^{2\ssss_{2}}}{\rr^{2\beta}}\rr^{2\alpha}\mint_{B_{\rr}}(|Du|+1)^{2(q-1+\delta_{2})} \dx\nonumber \\
&&  \qquad \quad  +\frac{M^{2\ssss_{3}}}{\rr^{2\beta}}\rr^{\alpha}\mint_{B_{\rr}}(|Du|+1)^{q-1+\delta_{2}} \dx \notag\\
&& \quad \qquad  +\frac{M^{2\ssss_{3}}}{\rr^{2\beta}}\rr^{\alpha_0}\mint_{B_{\rr}}(|Du|+1)^{3\delta_{2}} \dx \label{16}
\end{eqnarray} 
where $v$ is as in \rif{intrinseca} and this time $M$ satisfies the intrinsic relation 
$$
\sup_{x \in B_{\rr}}\,  |Du(x)|^{2-\gamma}+\aaa(x)|Du(x)|^{q}  + 1 \leq M\,.
$$ 
Here the free parameters $\delta_{1,2,3}$ can be chosen small at will and the parameters $\ssss_{1,2,3}$ are suitable numbers related to $\delta_{1,2,3}$ and $\gamma,q $. Inequality \rif{16} provides useful estimates provided the sharp bound \rif{ilbb} 
is satisfied and $\gamma$ is suitably close to $1$. The crucial point in \cite{loggy} is the intrinsic approach fully using the structure in \rif{generale} for which the choice in \rif{intrinseca} becomes effective.

We finally observe that the approach to $L^\infty$-bounds via fractional De Giorgi's iterations has been introduced, independently, in \cite{caff} and \cite{mis3}. In both cases the starting point is a fractional Caccioppoli inequality. In case of \cite{caff} this is natural as the problem considered involves a nonlocal operator (like a fractional power of the Laplacean) and therefore the starting energy inequality comes from simple testing. In \cite{mis3} the equation considered is of the type in \rif{risolvi}, with $p=2$ and $\mu$ being a Borel measure with finite total mass, and the decomposition argument explained in this section is required. In this last case, fractional Caccioppoli's inequalities are obtained for the function $v=|Du|$. An approach using fractional De Giorgi's classes has been used later in \cite{cozzi, naka}. See \cite{minpre, Dark} for fractional estimates for nonlinear elliptic systems and their applications to singular sets estimates. 
} 
\end{commentary}
\begin{boxy}{The key takeaway}
When full higher derivatives of solutions are not available, try fractional ones. The scheme [renormalized  energy estimates $\oplus$ De Giorgi type iterations $\oplus$ nonlinear potentials] is robust enough to resist  the lack of both full second order derivatives and uniform ellipticity. This allows to replace differentiability of coefficients with H\"older continuity and to give a unified approach to gradient boundedness of solutions, avoiding perturbation methods already in classical cases. 
\end{boxy}
\subsection{Stein type theorems \cite{piovra}} Here we briefly sketch the proof of the a priori estimate of Theorem \ref{t5}, and again in the case $p\geq 2$; we keep the notation introduced in the previous sections. Full details of the proofs can be found in \cite{piovra}.  As described in Section \ref{noeln}, this time we cannot even use \rif{eln} since $y \mapsto h(\cdot, y)$ is only H\"older continuous. On the other hand, it is still possible to get a fractional Caccioppoli inequality \cite[Proposition 5.2]{piovra}, that now looks like
\begin{flalign}\label{xx.1011} 
& \int_{B_{\varrho/2}} \mint_{B_{\varrho/2}}  
\frac{|(v-\kk)_{+}(x)- (v-\kk)_{+}(y)|^{2}}{|x-y|^{n+2\beta}} \dx \dy
\notag \\ & \quad \ \ \lesssim_{\data,\alpha, \beta}  \frac{M^{\ssf (q-2)}}{\varrho^{2\beta}}\mint_{B_{\rr}}(v-\kk)_{+}^{2} \dx  +\frac{M^{\ssf  q}}{\rr^{2\beta}}\rr^{\frac{2\alpha}{2-\alpha}}\mint_{B_{\rr}}\mu^{\frac{2}{2-\alpha}}  \dx\,.
\end{flalign}
This holds whenever $B_{\rr} \Subset \Omega$ and $\kk\geq 0$, where $M$ is as in \rif{parallela}, $\ssf$ as in \rif{sceltissime2}, $v$ as in \rif{triala} and this time $\beta$ satisfies 
\eqn{newrange}
$$
0<\beta<\frac{\alpha}{2+\alpha}\,.
$$
Of course, we still assume the a propri regularity in \rif{assuntoc1}; this additional property can be removed via a (delicate) approximation procedure as detailed in \cite{piovra}. Using \rif{immersione} with $w\equiv (v-\kk)_{+}$ in \rif{xx.1011} yields 
\begin{flalign*}
\notag   & \left(\mint_{B_{\rr/2}}(v-\kk)_{+}^{2\chi}  \dx\right)^{1/\chi}  \\
&\qquad \lesssim_{\data,\alpha, \beta} M^{\ssf (q-2)} \mint_{B_{\rr}}(v-\kk)_{+}^{2} \dx + M^{\ssf  q}\rr^{\frac{2\alpha}{2-\alpha}}\mint_{B_{\rr}}\mu^{\frac{2}{2-\alpha}}  \dx
\end{flalign*}
where $\chi\equiv \chi(\beta)>1$ is as in \rif{classico2}. 
Again, compare this last inequality with \rif{reversa} and with \rif{reversanorm}. In all cases the last integrals account for the external ingredients: coefficients or data. We are now in business to apply Lemma \ref{revlem}, that gives 
\begin{flalign*}
 v(x_{0})   &\lesssim_{\data,\alpha, \beta} M^{\frac{\chi\ssf(q-2)}{2(\chi-1)}}\left(\mint_{B_{r_{0}}(x_{0})}v^{2}  \dx\right)^{1/2}\\
 & \quad +M^{\frac{\ssf(q-2)}{2(\chi-1)}+\frac{\ssf q}{2}}\mathbf{P}^{2/(2-\alpha),1}_{2,
\alpha/(2-\alpha)}(\mu;x_{0},2r_{0})
\end{flalign*}
so that, after a few elementary manipulations, we get 
\begin{flalign*}
\nonumber
|Du(x_{0})| & \lesssim_{\data,\alpha, \beta} M^{\frac{\chi\ssf}{2(\chi-1)}\left(\frac q2-1\right)+\frac 12}\left(\mint_{B_{r_{0}}(x_{0})}(|Du|^2  +1)\dx\right)^{1/4}\\
 & \quad 
\ \ +M^{\frac{\ssf}{2(\chi-1)}\left(\frac q2 -1\right)+\frac{\ssf q}{4}}\left[\mathbf{P}^{2/(2-\alpha),1}_{2,
\alpha/(2-\alpha)}(\mu;x_{0},2r_{0})\right]^{1/2}\,.
\end{flalign*}
As after \rif{concentriche}, we find 
\begin{flalign*}
 \notag  \texttt{h}(\tau_1)   & \lesssim_{\data,\alpha, \beta} \frac{[ \texttt{h}(\tau_2)]^{\frac{\chi\ssf}{2(\chi-1)}\left(\frac q2-1\right)+\frac 12}}{(\tau_2-\tau_1)^{\frac{n}{2p}}}\left(\int_{B_r}[F(Du)+1]  \dx\right)^{1/4}\\
  & \quad 
  + [\texttt{h}(\tau_2)]^{\frac{\ssf}{2(\chi-1)}\left(\frac q2 -1\right)+\frac{\ssf q}{4}} \|\mathbf{P}^{2/(2-\alpha),1}_{2,
\alpha/(2-\alpha)}(\mu;\cdot,\tau_2-\tau_1) \|_{L^\infty(B_{\tau_1})}^{1/2}+ 1\,.
\end{flalign*}
From this inequality we can now proceed as after \rif{brahms}. Lemma \ref{crit} now gives
$$
\|\mathbf{P}^{2/(2-\alpha),1}_{2,
\alpha/(2-\alpha)}(\mu;\cdot,\tau_2-\tau_1) \|_{L^\infty(B_{\tau_1})}\lesssim_{n,\alpha}  \|f\|^{1/(2-\alpha)}_{n/\alpha,1/(2-\alpha);B_{r}}  \,.
$$
Then we require 
$$
\begin{cases}
\displaystyle \frac{\chi(\beta)\ssf}{2(\chi(\beta)-1)}\left(\frac q2-1\right)+\frac 12<1\\
\displaystyle  \frac{\ssf}{2(\chi(\beta)-1)}\left(\frac q2 -1\right)+\frac{\ssf q}{4}<1\,.
\end{cases}
$$
Also in this case the second inequality implies the first and therefore we finally ask that
$$
\left(\frac q2 -1\right)\frac{\ssf(n-2\beta)}{4\beta} +\frac{\ssf q}{4} <1\,.
$$
In turn, by using \rif{boundlolo}, we can find $\beta$ within the range \rif{newrange} such that the above inequality is true. The rest of the proof, leading to the a priori estimate of Theorem \ref{t5}, follows as in Section \ref{fracschau}. 
\subsection{Optimal gap bounds \cite{gioiello}}\label{maingioi2} A key to the proof of Theorem \ref{maingioi} is a preliminary higher integrability result, i.e., 
\eqn{giogio}
$$
Du \in L^{\mf{q}}_{\loc}(\Omega)\,, \quad \mbox{for every $\mf{q} <\infty$}\,.
$$
This is based on a careful use of delicate Besov spaces techniques. Besov spaces methods are sometimes employed in uniformly elliptic problems with a certain degree of lack of ellipticity like for instance operators in the Heisenberg group \cite{d}, or in nonlocal problems \cite{bl, bls, gl}. The proof of \rif{giogio}  relies on estimates of the type
\eqn{derive}
$$
\left\|\frac{\tau_{h}(\tau_{h}u)}{|h|^{s_i}}\right\|_{L^{\mf{q}_i}} \leq c_i
$$
for $h \in \er^n$ and $|h|$ suitably small and with sequences
$
\mf{q}_i \to \infty$, $1< s_i  \to 1$ an $c_i\to \infty$. Such estimates eventually imply \rif{giogio}. The proof of \rif{giogio} also yields an intermediate result, connected to Theorem \ref{doubleschauder}, \rif{twist2}, and Theorem \ref{choetheorem}, in the sense that the local boundedness of minima makes the gap bounds independent of $n$. 
\vspace{3mm}\begin{mdframed}
\begin{theorem}[Non-dimensional gap bound]\label{main3} Let $u\in W^{1,p}_{\loc}(\Omega)\cap L^{\infty}_{\loc}(\Omega)$  be a minimizer of the functional $\overline{\mathcal F_{\texttt{x}}}(\cdot, B)$ in \eqref{rilassato} for every ball $B\Subset \Omega$, under assumptions \eqref{xx.0} with $p\leq n$ and 
\eqn{pqbb}
$$
 q<p+\alpha\min\left\{1,p/2\right\}\,.
$$
Then \eqref{giogio} holds. Moreover, 
$$
\|Du\|_{L^{\mf{q}}(B_{\rr/2})}^{\mf{q}}\lesssim_{\data}  \left(\frac{\|u- (u)_{B_{\rr}}\|_{L^{\infty}(B_{\rr})}}{\rr} +1\right)^{\mf{b}_{\mf{q}}}\left[\overline{\mathcal F_{\texttt{x}}}(u,B_{\rr})+\rr^n\right]
$$
holds for every $\mf{q}<\infty$, and for every ball $B_{\rr}\Subset \Omega$, $\rr \leq 1$, where $\mf{b}_{\mf{q}}\equiv \mf{b}_{\mf{q}}(n,p,q,\alpha, \mf{q})$. Finally, if $\mathcal {L}_{\mathcal{F}_{\texttt{x}}}(u,B)=0$ holds for every ball $B\Subset \Omega$,  then the  same result holds for $W^{1,p}_{\loc}(\Omega)\cap L^{\infty}_{\loc}(\Omega)$-regular minimizers $u$ of $\mathcal F_{\texttt{x}}$ in \trif{FMx}. 
\end{theorem}
\end{mdframed}\vspace{3mm}
When $p\geq 2$ condition \rif{pqbb} is sharp, again by the counterexamples in \cite{sharp, FMM}. Note that in Theorem \ref{main3} we are restricting to the case $p\leq n$ otherwise $u \in  L^{\infty}_{\loc}(\Omega)$ is automatically verified and indeed condition \rif{pqbb} is more restrictive than $q/p<1+\alpha/n$.  

Going back to the proof of Theorem \ref{maingioi}, the use of higher integrability \rif{giogio} allows to upgrade the fractional technique explained in Section \ref{fracschau} and improve the exponents via a new fractional Caccioppoli type inequality of the type in \rif{xx.10}. A crucial role in obtaining the sharp condition $q/p<1+\alpha/n$ is at this stage played by the use of Theorem \ref{bstheorem} (and the related a priori estimate).  In other words,  estimate \rif{milano2} holds with
\begin{boxy}{The new exponent $\ssf$}
\eqn{nuovino}
$$
\mf{s}:=\begin{cases}
\displaystyle 
\, \frac{2q}{(n+1)p-(n-1)q}\qquad \mbox{if either $n\ge 4 $ or $n=2$}\\
\\
\displaystyle
\, \mbox{any number} > \frac{q}{2p-q}\qquad \mbox{if $n=3$ and $q>p$}
\\ \\
\displaystyle
\, 1 \qquad \mbox{if $n=3$ and $q=p$}\,.
\end{cases}
$$
\end{boxy}
This follows using the results in \cite{BS, BS2, schag}. This new exponent is smaller than the one appearing in \rif{sceltissime2}. Summarizing, passing from the first gap bound in \rif{bound333} to the sharp one $q/p<1+\alpha/n$ goes in two steps:
\begin{itemize}
\item The structural improvement. Using \rif{giogio} allows to improve the quadratic decay gap bound in \rif{bound333} to a  linear one, that is something of the form \rif{prebound}. 
\item The fine improvement. The use of estimate \rif{milano2} with the new exponent $\sss$ in \rif{nuovino}, allowing to finally reach $q/p<1+\alpha/n$. In this respect note that the use of the results in  \cite{BS, BS2, schag} is really essential in that the use of the previous exponent in \rif{sceltissime2} would at this stage lead to yet the suboptimal bound 
$$
\frac qp < 1 +\frac{\alpha}{n+1}\,.
$$
\end{itemize}

\section{Moser strikes again \cite{dm}}\label{mosersec}
In the above section we have seen the general three-step scheme 
\eqn{itre}
$$
\begin{cases}
\mbox{Renormalized Caccioppoli type inequalities}\\
\mbox{De Giorgi type iterations}\\
\mbox{Use of Nonlinear potentials}.
\end{cases}
$$
As observed in \cite{dm} a similar, parallel approach using the first two points in \rif{itre} can be employed via Moser iteration. We indeed sketch how to obtain an a priori estimate for minima of functionals of the type \rif{genF} with $\mu\equiv 0$. The assumptions 
we are considering here are 
\eqn{assF2fine}
$$
\left\{
\begin{array}{c}
[H_{1}(z)]^{p/2}\le F(x,z)\le L [H_{1}(z)]^{q/2}\\ [6 pt]
[H_{1}(z)]^{(p-2)/2}|\xi|^{2}\le \partial_{zz}F(x,z)\, \xi\cdot\xi\\ [6 pt]
\displaystyle
|\partial_{zz}F(x,z)| [H_{1}(z)]^{1/2} +|\partial_{zx}F(x,z)|\le L [H_{1}(z)]^{(q-1)/2} \;,
 \end{array}\right.
 $$
for every $x\in \Omega$, $z, \xi \in \er^n$, where $L\geq 1$; the integrand $F \colon \Omega \times \er^n\to [0, \infty)$ is assumed to be locally $C^2$-regular with respect to the gradient variable and is continuous; finally its derivatives are Carathe\'odory-regular. We recall that, according to \rif{defiH}, it is $H_{1}(z)=|z|^{2}+1$. 
Assumptions \rif{assF2fine} guarantee that minimizers $u$ satisfy the a priori estimate 
\eqn{toylip0}
$$
\|Du\|_{L^{\infty}(B/2)} \lesssim_{\data}\left(\mint_{B}F(x, Du)\dx \right)^{\frac{1}{(n\mathfrak{b}+1)p-n\mathfrak{b}q}}+1\,,
$$
whenever $B\equiv B_{\rr}\Subset \Omega$ is such that $\rr\leq 1$, where $\data=(n,p,q, L)$ and
\eqn{bboo}
$$
\mathfrak{b}:=
\begin{cases}
\displaystyle 1 \ \mbox{if $n>2$}\\
\displaystyle >1 \ \mbox{such that $2\mathfrak{b}(q-p)<p$ if $n=2$}\,,
\end{cases}
$$
provided
\eqn{boundx}
$$
\frac{q}{p}< 1 +\frac1n
$$
is satisfied and the approximation in energy property \rif{energyconv} holds for $w\equiv u$ (absence of Lavrentiev phenomenon), i.e., 
\eqn{basedpon}
$$
\begin{cases}
\, \mbox{there exists $\{u_k\}_{k}\subset C^{\infty}(B)$ such that $u_k\to u$ strongly in $W^{1,p}(B)$}\\
\, \mbox{and $F(\cdot,Du_k)\to F(\cdot, Du)$ in $L^1(B)$}\,.
\end{cases}
$$ 
Notice that condition \rif{boundx} is optimal by the counterexamples in \cite{sharp,FMM}; see Section \ref{softy2}. Notice also that in \rif{basedpon} it is equivalent to require $\{u_k\}_{k}\subset W^{1,q}(B)$ due to the growth assumptions in \rif{assF2fine}$_1$.

We shall give a rapid proof of \rif{toylip0} as an a priori estimate, that is, assuming the a priori higher integrability $u \in W^{1,\infty}(B)$; the general case can then be treated via an approximation argument, based on \rif{basedpon}, as for instance shown in \cite{sharp}. The point we want to stress now is the following: similarly to the method demonstrated in Section \ref{rinormalizza}, renormalizing the relevant energy inequalities, and taking advantage of the precise dependence of the constants along the iterations (in this case Lemma \ref{lamoser}) rapidly leads to \rif{toylip0}. The outcome is a very streamlined proof that in fact literally reduces estimates for $(p,q)$-growth functionals to those valid when $p=q$. For the proof of \rif{toylip0} we preliminary observe that we can reduce to the case $B\equiv \mathcal B_1:= B_1(0)$ by using the same scaling argument in \rif{riscale} (in the following lines we denote $\mathcal B_{\tau}\equiv B_\tau(0)$). The function in \rif{riscale} locally minimizes the same functional endowed with new (rescaled) integrand $(x,z)\mapsto F(x_0+\rr x, z)$, which satisfies \rif{assF2fine} as we are assuming that $\rr\leq 1$. We still denote by $u$ the blown-up function in \rif{riscale}. Following \cite{dm},  
we first recall that in the uniformly elliptic case $p=q$ the inequality
\eqn{revvi1}
$$
\displaystyle \left(\int_{\mathcal B_{\varrho_1}}v^{\left(\gamma+p/2\right)\chi} \dx \right)^{1/\chi}
 \lesssim_{n,p} \frac{L^2(1+\gamma)^2}{(\varrho_2-\varrho_1)^{2}}\int_{\mathcal B_{\varrho_2}}v^{\gamma+p/2}  \dx\,, \quad \forall \ \gamma\geq 0
 $$
 with 
$$
 v := H_1(Du) = |Du|^2+1
$$
holds whenever $ \mathcal B_{\tau_1}\Subset \mathcal B_{\varrho_1} \Subset \mathcal B_{\varrho_2} \Subset \mathcal B_{\tau_2}\Subset \mathcal B_1$, and $\chi =2^{*}/2$. As usual, $2^*>2$ is the Sobolev embedding exponent defined as in \rif{sobby0}. Note that
\eqn{chichi}
$$
\frac{\chi}{\chi-1}= \frac{2^*}{2^*-2}
$$
and that a crucial point here is that in \rif{revvi1} the dependence on $L$ is made explicit. 
Via Lemma \ref{lamoser}, and recalling \rif{chichi}, this yields
\eqn{localp}
$$
\|v\|_{L^{\infty}(\mathcal B_{\tau_1})} \lesssim_{n,p} \frac{L^{\frac{4}{p}\frac{2^*}{2^*-2}}}{(\tau_2-\tau_1)^{\frac{4}{p}\frac{2^*}{2^*-2}}} \|v\|_{L^{p/2}(\mathcal B_{\tau_2})}\,.
$$
We now show how to use the above estimate, that holds when $p=q$, to treat the case $p<q$ straightaway. For this, observe that, when handling integral quantities evaluated on the ball $B_{\tau_2}$ we have
$$
|z| \equiv |Du| \leq \|Du\|_{L^{\infty}(\mathcal B_{\tau_2})}
$$
and therefore we can always replace \rif{assF2fine}$_3$ by
\eqn{revvi4} 
$$
\begin{cases}
\, \displaystyle |\partial_{zz}F(x,z)|[H_1(z)]^{1/2}+ |\partial_{zx}F(x,z)|\le  \mathfrak{L}[H_1(z)]^{(p-1)/2} \ \mbox{on $B_{\tau_2}$}\,,\\
\, \mathfrak{L}:=L\|H_1(Du)\|_{L^{\infty}(B_{\tau_2})}^{(q-p)/2}\equiv L\|v\|_{L^{\infty}(B_{\tau_2})}^{(q-p)/2}\,.
\end{cases}
$$
Indeed, all we need to prove \rif{revvi1} is the Euler-Lagrange equation
$$
\int_{\mathcal B_1} \partial_zF(x, Du)\cdot D\varphi\dx =0
$$
with sutable test functions $\varphi$ supported $B_{\tau_2}$; this means that we can use \rif{revvi4}. Note that the only part of \rif{assF2fine} which is used to prove \rif{revvi1} is actually \rif{assF2fine}$_{2,3}$. 
Applying \rif{localp} with $L$, replaced by $\mathfrak L$, we now gain
\eqn{revvi30}
$$
\|v\|_{L^{\infty}(\mathcal B_{\tau_1})} \lesssim_{n,p,L} \frac{ \|v\|_{L^{\infty}(\mathcal B_{\tau_2})}^{\frac{2(q-p)}{p}\frac{2^*}{2^*-2}}}{(\tau_2-\tau_1)^{\frac{4}{p}\frac{2^*}{2^*-2}}} \|v\|_{L^{p/2}(\mathcal B_{\tau_2})}\;.
$$
When $n>2$ we note that 
$$\frac{2(q-p)2^*}{p(2^*-2)} = \frac{(q-p)n}{p} \stackrel{\eqref{boundx}}{<}1\,.$$
Young's inequality in \rif{revvi30} then yields
\eqn{revvi3}
$$
\texttt{h}(\tau_1) \leq \frac {\texttt{h}(\tau_2)}{2}
+ \frac{c\|v\|_{L^{p/2}(\mathcal B_{\tau_2})}^{\frac{p}{(n+1)p-nq}}}{(\tau_2-\tau_1)^{\frac{2n}{(n+1)p-nq}}}
$$
where 
$
\texttt{h}(t):= \|v\|_{L^{\infty}(B_{t})}$ and $c\equiv c (n,p,)$.
Lemma \ref{l5} finally gives
\begin{flalign}
\notag \|v\|_{L^{\infty}(\mathcal B_{1/2})}  &\lesssim_{\data} \|v\|_{L^{p/2}(\mathcal B_{1})}^{\frac{p}{(n+1)p-nq}}\lesssim_{\data}  \left(\int_{\mathcal B_1}(|Du|^p+1)\dx \right)^{\frac{2}{(n+1)p-nq}}
\label{toylip}
\end{flalign}
from which \rif{toylip0} follows using \rif{assF2fine}$_1$. This proves \rif{toylip0} in the case $n>2$ and $\rr=1$; scaling back we get the general case. 
When $n=2$ note that $q/p<3/2$, that is the assumed bound \rif{boundx}, implies $2(q-p)/p<1$, therefore we take $2^*$ large enough (recall \rif{sobby0}) in order to achieve 
$$
\frac{2(q-p)2^*}{p(2^*-2)} \leq \frac{2\mathfrak{b} (q-p)}{p} \stackrel{\eqref{bboo}}{<}1$$
so that \rif{revvi30} is replaced by 
$$
\|v\|_{L^{\infty}(\mathcal B_{\tau_1})} \lesssim_{p} \frac{ \|v\|_{L^{\infty}(\mathcal B_{\tau_2})}^{2\mathfrak{b}(q-p)/p}}{(\tau_2-\tau_1)^{4\mathfrak{b}/p}} \|v\|_{L^{p/2}(\mathcal B_{\tau_2})}
$$
and we can proceed as in the case $n>2$. 
\begin{boxy}{The key takeaway}
Renormalized energy inequalities can be used also in the setting of Moser iteration technique. This, together with quantitative knowledge of the constants along the interation schemes, allows in several cases to reduce the treatment of nonuniformly elliptic problems with polynomial growth to that of uniformly elliptic ones.  
\end{boxy}

\renewcommand{\refname}{}    

\end{document}